\documentclass{article}

     \PassOptionsToPackage{round}{natbib}

     \usepackage[preprint]{neurips_2023}

\usepackage[utf8]{inputenc} 
\usepackage[T1]{fontenc}    
\RequirePackage{titletoc}
\usepackage[colorlinks,citecolor=blue]{hyperref}       
\usepackage{hyperref}       
\usepackage{url}            
\usepackage{booktabs}       
\usepackage{amsfonts}       
\usepackage{nicefrac}       
\usepackage{graphicx}
\usepackage{microtype}      
\usepackage{algorithm}
\usepackage{algpseudocode}
\usepackage{subcaption}
\usepackage{multirow}
\usepackage{xcolor}         

\title{Speeding up Monte Carlo Integration: \\ Control Neighbors for Optimal Convergence}

\author{%
 Rémi Leluc \\
 CMAP, \'Ecole Polytechnique,\\
 Institut Polytechnique de Paris, France\\
  \texttt{remi.leluc@gmail.com} \\
   \And
	Fran\c{c}ois Portier \\
	CREST \\
	ENSAI, France \\
	\texttt{francois.portier@gmail.com} \\
	\And
	Aigerim Zhuman \\
	LIDAM, ISBA\\
	UCLouvain, Belgium \\
	\texttt{aigerim.zhuman@uclouvain.be} \\
	\And
	Johan Segers \\
	LIDAM, ISBA\\
	UCLouvain, Belgium \\
	\texttt{johan.segers@uclouvain.be} \\
}

\include{preambule}
\begin{document}

\maketitle

\begin{abstract}
A novel linear integration rule called \textit{control neighbors} is proposed in which nearest neighbor estimates act as control variates to speed up the convergence rate of the Monte Carlo procedure on metric spaces. The main result is the $\Oh(n^{-1/2} n^{-s/d})$ convergence rate -- where $n$ stands for the number of evaluations of the integrand and $d$ for the dimension of the domain -- of this estimate for Hölder functions with regularity $s \in (0,1]$, a rate which, in some sense, is optimal. Several numerical experiments validate the complexity bound and highlight the good performance of the proposed estimator.
\end{abstract}

\section{Introduction} \label{sec:1_intro}

Consider the classical numerical integration problem of approximating the value of an integral $\meas(\func) = \int \func \, \diff \meas$ where $\meas$ is a probability measure on a metric space $(\ms, \rho)$ and the integrand $\func$ is a real-valued function on the support of $\mu$. Suppose that random draws from the measure $\meas$ are available and calls to the function $\func$ are possible. The standard Monte Carlo estimate consists in averaging $\func(X_i)$ over $i=1,\ldots,n$, where the particles $X_i$ are drawn independently from $\meas$. For square-integrable integrands, the Monte Carlo estimate has convergence rate $\Oh(n^{-1/2})$ as $n \to \infty$, whatever the dimension of the domain. In some applications, calls to the integrand may be expensive \citep{sacks1989design,doucet2001sequential}
and one may only have access to a small number of evaluations $\func(X_i)$, such as in Bayesian inference on complex models \citep{higdon2015bayesian,toscano2022bayesian}. The $n^{-1/2}$-rate of the standard Monte Carlo estimate then becomes too slow and leads to highly variable estimates.

As detailed in \citet{novak2016some}, the complexity of integration algorithms may be analyzed through the convergence rate of the error. Any randomized procedure based on $n$ particles yields an estimate $\hat{\meas}_n(\func)$ of the integral $\meas(\func)$ and the root mean-square error of the procedure is $\Exp[|\hat{\meas}_n(\func) - \meas(\func)|^2]^{1/2}$. For the specific problem of integration with respect to the uniform measure over the unit cube $[0,1]^d$, the complexity rate of randomized methods for Lipschitz integrands is known to be $\Oh(n^{-1/2} n^{-1/d})$ \citep{BAKHVALOV2015502,novak2016some}. Furthermore, when the first $s$ derivatives of the integrand are bounded, the convergence rate becomes $\Oh(n^{-1/2} n^{-s/d})$. 
These complexity rates are informative as they show that for smooth integrands, the Monte Carlo estimate is suboptimal and leaves room for improvement by relying on the regularity of the integrand. Several approaches are already known for improving upon the Monte Carlo benchmark. They can be classified according to their convergence rates while keeping in mind the lower bound $\Oh(n^{-1/2} n^{-s/d})$.

The control variate method \citep{rubinstein1981simulation,newton1994variance,caflisch1998monte,evans2000,glynn2002some,glasserman2004monte} is a powerful technique that allows to reduce the variance of the Monte Carlo estimate by approximating the integrand with a function whose integral is known. The Monte Carlo rate can be improved by combining control variates to model $\func$ \citep{oates2017control,portier2019monte,leluc2021control,south2022regularized}. 
In \cite{portier2019monte}, when $\meas$ is the uniform distribution on $[-1, 1]^d$ and when using $m$ orthogonal polynomials as control variates, the convergence rate is $\Oh(n^{-1/2} m^{-s/d})$, with $s$ is the regularity of $\func$. However, $m$ is required to be of smaller order than $n^{1/3}$, so that the optimal rate cannot be achieved. By relying on control functions constructed in a reproducing kernel Hilbert space, \cite{oates2019convergence} obtained the rate $\Oh(n^{-1/2 - (a \wedge b)/d + \varepsilon})$ for a specific class of integrands where $a$ is related to the smoothness of the target density, $b$ is related to the smoothness of the integrand $\func$ and $\varepsilon>0$ hides logarithmic factors and can be arbitrary small. This is done at the expense of computational complexity in $\Oh(n^3)$ due to matrix inversion when solving the underlying kernel ridge regression problem.

Another reliable technique to improve the rate of convergence of standard Monte Carlo is stratification. The sample space is partitioned and particles are sampled from each piece of the partition separately. It has allowed to improve the convergence rate of Monte Carlo estimates \citep{haber1966modified,haber1967modified} and to derive a general framework called stochastic quadrature rules \citep{haber1969stochastic}. Recently, Haber's work has been extended to take advantage of higher smoothness in the integrand \citep{chopin2022higher}. To the best of our knowledge, the works of \citet{haber1966modified} and \citet{chopin2022higher} are the first ones to provide a feasible method achieving the best rate of convergence for smooth integrands; an algorithm with optimal rate of convergence but difficult to implement is proposed in \cite{krieg+novak:17}. Some related ideas for building optimal algorithms were already introduced in the seminal work of \cite{BAKHVALOV2015502}; see also \cite{novak2016some} for a short review and \cite{novak2006deterministic} for more details. A drawback of the methods proposed in \citet{haber1966modified} and \citet{chopin2022higher} is that they are only valid for integration over the unit cube and involve a geometric number ($n=\ell^d$) of evaluations of the integrand $\func$, which may be somewhat restrictive in case the integrand is hard to compute, as in complex Bayesian models.

Interestingly, the stratification method in \citet{chopin2022higher}
relies on a piecewise constant control function with a low bias compared to a traditional regression estimate. Such type of estimator is shown to attain the optimal rate of convergence in case of noiseless evaluations of the function to approximate \citep{kohler2013optimal,hinrichs2020power} and this explains the success of their approach to reach to optimal rate of convergence for the integration problem.
The same idea of using a low-bias estimate is the starting point of this paper too and is relevant because the integrand $\func$ is accessible without noise. Low-bias estimates have also been employed successfully in adaptive rejection sampling \citep{achddou2019minimax}, allowing to reach the optimal rate. Alternatively, stochastic quadrature rules based on determinantal point processes have been proposed in \cite{bardenet2020monte} and \cite{coeurjolly2021monte}. Such methods allow to reduce the root mean squared error to $\Oh(n^{-1/2} n^{-1/(2d)})$ when the integrand $\func$ has a continuous derivative. This interesting acceleration still remains slower than the optimal lower bound.

In this paper, a new Monte Carlo method called \textit{control neighbors} is introduced. By using $1$-nearest neighbor estimates as control variates, it produces an estimate $\hat{\meas}_n(\func)$ of the integral $\meas(\func)$ for a probability measure $\meas$ on a metric space $(\ms, \rho)$ such that the measure of a ball of radius $r > 0$ is of the order $r^d$ as $r \to 0$, uniformly over the space. This novel estimate is shown to achieve the convergence rate $\Oh(n^{-1/2} n^{-s/d})$ for Hölder integrands of regularity $s \in (0, 1]$. The main properties of the control neighbors estimate are as follows: 
\begin{enumerate}[label=(\alph*)]
\item 
	It can be obtained under the same framework as standard Monte Carlo, i.e., as soon as one can draw random particles from $\meas$ and evaluate the integrand $\func$. In contrast to the classical method of control variates, the existence of control variates with known integrals is not required.
\item
	The method takes the form of a linear integration rule $ \sum_{i=1}^n w_{i,n} \func(X_i)$ with weights $w_{i,n}$ not depending on the integrand $\func$ but only on the particles $X_1, \ldots, X_n$. This property is computationally beneficial when several integrals are to be computed with the same measure $\meas$.
\item 
	For Lipschitz integrands $\func$, when $s = 1$, the convergence rate is the optimal one $\Oh(n^{-1/2} n^{-1/d})$ \citep{novak2016some}. Other recent approaches for general $\meas$ \citep[e.g.][]{oates2017control,portier2019monte} do not achieve this rate. Such a rate is achieved by the methods in \cite{haber1966modified} and in  \cite{chopin2022higher} provided $\mu$ is the uniform distribution on a $d$-dimensional cube and the number of evaluations of the integrand is geometric, $n=\ell^d$ for some integer $\ell \ge 1$.
	\item 
	The approach is \textit{post-hoc} in the sense that it can be run after sampling the particles and independently from the sampling mechanism. Although the theory in our paper is restricted to independent random variables, the method
	can be implemented for other sampling designs including MCMC or adaptive importance sampling.
	\item
	The method not only applies to integration with respect to smooth measures on compact subsets of Euclidean space but also with respect to the normalized volume measure on compact Riemannian manifolds such as the unit sphere or the orthogonal group.
\end{enumerate} 

A similar approach to the one in this paper (and developed independently) has been proposed in \cite{blanchetcan} as an example of an algorithm that achieves the complexity bounds for the integration problem; see their Section~4 and in particular their Theorem~4.2. Their study is focused on the derivation of complexity rates and on the existence of methods to achieve these rates. They consider the case where the integrand is evaluated with noise, and this leads to different choices of the number of neighbors (as explained in the proof of their Theorem~4.2). When the noise vanishes, their approach is similar to our $1$-nearest neighbor control variate estimate. One major difference remains that, for the proof of their theorem, they need to split the sample into two independent parts, leading to a loss of efficiency in practice. This latter drawback is avoided in this study with the help of a careful analysis of the leave-one-out algorithm.

The outline of the paper is as follows. The mathematical foundations of nearest neighbor estimates are gathered in Section~\ref{sec:2_nn} with a formal introduction of two different nearest neighbors estimates. The theoretical properties of the control neighbor estimates are stated in Section~\ref{sec:results}. Finally, Section~\ref{sec:simu} reports on several numerical experiments along with some practical remarks on the implementation of the proposed estimates, and Section~\ref{sec:discussion} concludes.
The supplement  contains proofs, auxiliary results and additional experiments.

\section{From nearest neighors to control neighbors}
\label{sec:2_nn}

This section presents the mathematical framework of nearest neighbor estimates with reminders on Voronoi cells and central quantities for the analysis, namely the degree of a point and the average cell volume. Next, two \emph{control neighbor} estimates are introduced and some basic properties are stated.

\subsection{Nearest neighbors}
\label{subsec:NN}

Let $X_1,\ldots, X_n$ be independent and identically distributed random variables taking values in a metric space $(\ms, \rho)$, drawn from a distribution $\mu$ without atoms. Let $x \in \ms$ denote a generic point.

\begin{definition}[Nearest neighbors and distances] 
	\label{def:nn_dist}
	The nearest neighbor $\Nn(x)$ of $x$ among $X_1,\ldots, X_n$ and the associated distance $\hat \tau_n (x)$ are 
	\begin{equation*}
		\Nn(x) \in \argmin_{Y \in \{X_1,\dots,X_n\}} \rho(x,Y), \qquad \hat \tau_n (x) =\rho \rbr{\Nn(x), x}. 
	\end{equation*}
	When the $\argmin$ is not unique, $\Nn(x)$ is defined as the one point among the $\argmin$ having the smallest index. More generally, for $k \in \cbr{1,\ldots,n}$, let $\hat N_{n,k}(x)$ denote the $k$-nearest neighbor of $x$ and $\hat \tau_{n,k} (x) = \rho \rbr{ \hat N_{n,k}(x), x }$ the associated distance, breaking ties by the lexicographic order.
\end{definition}

The sample $\cX_n = \{X_1,\ldots,X_n\}$ defines a natural (random) partition of the integration domain when considering the associated Voronoi cells. Any such cell is associated to a given sample point, say $X_i$, and contains all the points $x$ of which the nearest neighbor is $X_i$.

\begin{definition}[Voronoi cells and volumes]
	\label{def:vor_cells}
	The Voronoi cells of $X_1,\ldots,X_n$ are
	\begin{equation*}
		\forall i = 1,\ldots,n, \qquad
		S_{n,i} = \left\{ x \in \ms :\, \Nn (x) = X_i \right\},
	\end{equation*}
	with Voronoi volumes $V_{n,i} = \meas(S_{n,i})$. 
\end{definition}

The $1$-NN estimate of $\func$ is defined as $\hat \func_n (x) = \func\rbr{\Nn(x)}$ for all $x \in \reals^d$ and is piece-wise constant on the Voronoi cells, i.e., $\hat \func_n (x) =  \sum_{i=1}^n \func(X_i) \ind _{ S_{n,i}} (x) $.

The leave-one-out rule is a general technique to introduce independence between the prediction and the evaluation points. It is used as a cross-validation strategy in order to tune hyper-parameters of statistical procedures \citep{stone1974cross,craven1978smoothing}. The leave-one-out version of $\hat \func_n$ without $X_i$ is denoted by $\hgni $ 
and is obtained in the exact same way as $\hat \func_n$ except that a slightly different sample---in which the $i$-th observation has been removed---is used.

\begin{definition}[Leave-one-out neighbors, Voronoi cells and volumes]
	\label{def:loo_vor_cells}
	Let $i \in \{1,\ldots,n\}$ and $\cX_n^{(i)} = \{X_1,\ldots,X_n\}\setminus \{X_i\}$. The leave-one-out neighbor of $x \in \ms$ is
	\[ 
	\Nni(x) \in \argmin_{Y \in \cX_n^{(i)}} \rho(x, Y).
	\]
	When the above $\argmin$ is not unique, $\Nni(x)$ is defined as the one point among the $\argmin$ having the smallest index.
	The leave-one-out Voronoi cell $S_{n,j} ^{(i)}$ denotes the $j$-th Voronoi cell in $\cX_n^{(i)}$, i.e.,
	\begin{equation*}
		\forall j \in \{1,\ldots,n\} \setminus \{i\}, \qquad
		S_{n,j} ^{(i)} = \cbr{ x \in \ms :\, \Nni(x) = X_j }.
	\end{equation*}
	The leave-one-out Voronoi volume is defined as
	$V_{n,j}^{(i)} = \meas\bigl(S_{n,j}^{(i)}\bigr)$. 
\end{definition} 

The leave-one-out $1$-NN predictor used in Section~\ref{subsec:CN} to define the proposed integral estimate is $\hgni(x) = \func\bigl(\Nni(x)\bigr)$. A key property is that $ \hgni$ and $ \hat \func_{n}$ coincide on $S_{n,j} $ for $j\neq i$. On the cell $S_{n,i}$, when the function $\func$ is Hölder with regularity $s \in (0, 1]$, the supremum distance between $ \hgni$ and $ \hat \func_{n}$ is of the same order as the nearest neighbor distance to the power $s$. In terms of the $L^1(\meas)$-norm, their difference is even smaller if the cell $S_{n,i}$ has a small $\mu$-volume. Relevant for our numerical integration problem is that the average of the integrals $\meas\bigl(\hgni\bigr)$ is close to $\meas(\hat{\func}_{n})$.

\begin{lemma}\label{prop:modification_integral}
	In terms of $\bar \func_n (x) =  \sum_{i=1}^n \hgni (x) \ind_{S_{n,i}} (x) $, we have	$	 \sum_{i=1}^n \bigl( \meas\bigl(\hgni\bigr) - \meas(\hat{\func}_{n}) \bigr)
	=  \meas(\bar{\func}_n - \hat{\func}_{n}).$
\end{lemma}

Enumerating how often a point $X_j$ is the nearest neighbor of points $X_i$ for $i \ne j$ reflects how much $X_j$ is surrounded within the sample. Another quantification of the isolation of a point $X_j$ is obtained by summing the Voronoi volumes $V_{n,j}^{(i)}$ over $i \ne j$. These two notions are formally stated in the next definition. 

\begin{definition}[Degree and cumulative volume] \label{def:degree_volume} For all $j=1,\ldots,n$, the degree $\hat{d}_{n,j}$ represents the number of times $X_j$ is a nearest neighbor of a point $X_i$ for $i \ne j$. The associated $j$-th cumulative Voronoi volume is denoted by $\hat{c}_{n,j}$, that is,
	\begin{equation*}
		\hat{d}_{n,j} = \sum_{i : i \neq j} \ind_{S_{n,j}^{(i)}} (X_i)
		\qquad \text{and} \qquad
		\hat{c}_{n,j} = \sum_{i : i \neq j} V_{n,j}^{(i)}. 
	\end{equation*}
\end{definition}

Interestingly, the degree of a point and its cumulative Voronoi volume have the same expectation: $\expec[\hat{d}_{n,j}] = \expec[\hat{c}_{n,j}] = 1$.

The two quantities $\hat d_{n,j}$ and $\hat c_{n,j}$ will be useful in Proposition~\ref{prop:quad} below to express the control neighbors estimate as a linear integration rule. For now, we note that  weighted sums of $\func(X_j)$ using $\hat d_{n,j}$ and $\hat c_{n,j}$ as weights are related to the leave-one-out estimate. 

\begin{lemma}
	\label{prop:di_ci} It holds that
	\begin{align*}
		\sum_{i=1}^n \func(X_i ) \, \hat d_{n,i} = 	\sum_{i=1}^n \hgni(X_i )
		\qquad  \text{ and } \qquad
		\sum_{i=1}^n \func(X_i ) \, \hat c_{n,i} = \sum_{i=1}^n  \meas\bigl(\hgni\bigr)  .
	\end{align*} 
\end{lemma}

\subsection{Control neighbors}
\label{subsec:CN}

With the help of previous notation, we now introduce the two \textit{control neighbor} estimates
\begin{align}
	\label{first_estimate2} 
	\munNN(\func)   
	&= \frac{1}{n} \sum_{i=1}^n \left[
	\func(X_i) - \left\{ \hgni (X_i) -   \meas ( \hat \func_n  ) \right\}
	\right]\\
	\label{first_estimate}  \munNNloo(\func)
	&= \frac{1}{n} \sum_{i=1}^n \left[ 
	\func(X_i) - \left\{\hgni (X_i) -   \meas ( \hgni ) \right\}
	\right].		
\end{align}
The first one, $\munNN(\func) $, is the output in the nearest neighbor algorithm in Section \ref{subsec:num_details}. The second one, $\munNNloo(\func)$, is a slight modification of  $\munNN(\func) $ that is an unbiased estimate of $\mu(\varphi)$.
Indeed, a simple conditioning argument implies
\[
	\expec \sbr{ \hgni (X_i) -   \meas\bigl(\hgni\bigr)} 
	= \expec \sbr{ \expec \sbr{
	\hgni (X_i) -   \meas\bigl(\hgni\bigr)
	\,\Big\vert\, \cX_n^{(i)} 
	} } 
	= 0,
\] 
which is sufficient to get the zero-bias property 
\[
\expec\left[ \munNNloo(\func) \right] = \meas(\func).
\] 
The estimate $\munNN$ is not unbiased, but under the conditions of Theorem~\ref{th:cv_rate} below, its bias is of smaller order than its root mean squared error; see Remark~\ref{rem:bias} below.
Moreover, since $\hgni$ is similar to the $1$-NN estimate $\hat \func_n$ of $\func$ based on the full sample $X_1,\ldots,X_n$, their integrals should be close. This intuition is confirmed in Proposition~\propdifflooNN{} in the Supplementary material.
Consequently, the estimate \eqref{first_estimate} will play an important role in the theory. 
However, computing the terms $\meas\bigl(\hgni\bigr)$ for $i = 1, \ldots, n$ requires the evaluation of $n$ additional integrals. In practice, the working estimate is \eqref{first_estimate2} as it involves fewer computations. Further, it is worthwhile to note that our two control neighbor estimates are \emph{not} based on a single approximation of the integrand $\func$ but rather on $n$ different approximations $\hgni$.

The two control neighbor estimates in~\eqref{first_estimate2} and~\eqref{first_estimate} can be expressed as linear integration rules of the form $ \sum_{i=1}^n w_{i,n} \, \func(X_i)$ with weights $w_{i,n}$ not depending on the integrand $\func$. The weights involve the degrees $\hat d_{n,i}$ and the (cumulative) volumes ${V}_{n,i}$ and $\hat c_{n,i}$ in Definitions~\ref{def:vor_cells} and~\ref{def:degree_volume}.

\begin{proposition}[Quadrature rules] 
	\label{prop:quad}
	The estimates $\munNN(\func) $ and $\munNNloo(\func)$ can be expressed as linear estimates of the form
	\begin{align*}
		\munNNloo(\func)
		= \sum_{i=1}^ n w_{i,n}^{\mathrm{(NN-loo)}} \func(X_i)
		\quad \text{ and } \quad 
		\munNN(\func) 
		= \sum_{i=1}^ n w_{i,n}^{\mathrm{(NN)}} \func(X_i)
	\end{align*}
	where $w_{i,n}^{\mathrm{(NN-loo)}} = (1 + \hat c_{n,i} - \hat d_{n,i})/n$ and $w_{i,n}^{\mathrm{(NN)}} = (1 + n V_{n,i} - \hat d_{n,i})/n$.
\end{proposition}

The weights of the two estimates satisfy $\sum_{i=1}^n w_{i,n}=1$, meaning that the integration rules are exact for constant functions. 
In the light of Proposition~\ref{prop:quad}, the proposed estimate $\munNN(\func)$ consists in a simple modification of $\munNNloo(\func)$ by replacing $\hat c_{n,i}$, which involves $n-1$ Voronoi volumes, by $n V_{n,i}$. For Hölder integrands of order $s \in (0, 1]$, the difference between both estimates is of the order $n^{-1/2-s/d}$, in view of Eq.~\eqref{eq:diffmunloomunn} below and Proposition~\propdifflooNN{} in the Supplementary material.

\section{Main results}
\label{sec:results}

This section gathers some preliminary considerations on distributions on metric spaces and a bound on moments of nearest-neighbor distances before stating the main theoretical properties of the control neighbor estimates \eqref{first_estimate2} and \eqref{first_estimate}.
Under suitable conditions on the measure $\mu$ and the integrand $\func$, their root mean squared errors are shown to reach the rate $\Exp[|\hat{\meas}_n(\func) - \meas(\func)|^2]^{1/2} \lesssim n^{-1/2}n^{-s/d}$, where two positive sequences $a_n$ and $b_n$ satisfy $a_n \lesssim b_n$ provided $a_n \le C b_n$ for some constant $C > 0$. 
Further, their finite-sample performance is studied through a concentration inequality providing a high-probability bound on the error $|\hat{\meas}_n(\func) - \meas(\func)|$.

\subsection{Distributions on metric spaces}
\label{sec:metric}

Let $\meas$ be a Borel probability measure on a bounded metric space $(\ms, \rho)$ with diameter $\diam(\ms)$. To control the nearest neighbor distances, we need a handle on the probabilities of small balls. Let $B(x, r) = \cbr{y \in \ms: \rho(x,y) \le r}$ denote the ball with center $x \in \ms$ and radius $r > 0$. Intuitively, if $\meas$ is an absolutely continuous measure on a compact subset of $\reals^d$, then $\meas(B(x, r))$ is of the order $r^d$. The same is true if $\meas$ is a smooth measure on a $d$-dimensional manifold. These considerations motivate the following assumption, which implicitly defines a kind of intrinsic dimension of $\ms$ as measured by $\meas$.

\begin{enumerate}[label=(A\arabic*)]  
\item 
\label{hyp:ballmeasure} There exist $d > 0$ and $0 < C_0 \le C_1 < \infty$ such that
\[
	\forall x \in \ms, \;
	\forall r \in (0, \diam(\ms)], \qquad 
	C_0 r^{d} \le \meas(B(x,r)) \le C_1 r^{d}.
\]	
\setcounter{assumpnum}{\value{enumi}}
\end{enumerate}

\begin{remark}
\label{rem:radius}
In \ref{hyp:ballmeasure}, it is sufficient to require the stated inequalities for all $r \in (0, r_0]$ for some $r_0 > 0$. Upon changing $C_0$ and $C_1$, the inequalities then hold for all $r \in (r_0, \diam(\ms)]$ too. Indeed, for such $r$, we then have the upper bound $\meas(B(x,r)) \le 1 \le r_0^{-d} \cdot r^d$ and the lower bound $\meas(B(x, r)) \ge \meas(B(x,r_0)) \ge C_0 r_0^d \ge C_0 \{r_0/\diam(M)\}^d \cdot r^d$.
\end{remark}

\begin{remark}
\label{rem:fractal}
In theory, the exponent $d$ in \ref{hyp:ballmeasure} could be non-integer. As an example, suppose $\meas$ is the normalized Hausdorff measure on a fractal-like set, like the Cantor set ($d = \log 2 / \log 3$) or the Sierpinski triangle ($d = \log 3 / \log 2$), see \cite{SteinShakarchi+2005}. For the integration problem we study, this does not seem of much practical interest.
\end{remark}

A sufficient hypothesis for \ref{hyp:ballmeasure} is that $\mu$ is supported by a bounded subset of $\reals^d$ with positive $d$-dimensional Lebesgue measure and is absolutely continuous with a well-behaved density:
\begin{enumerate}[label=(A\arabic*')]  
\item \label{hyp:density} The measure $\meas$ is supported by a compact subset $\ms$ of $\reals^d$ and admits a density $\dens$ with respect to the $d$-dimensional Lebesgue measure $\Leb_d$. There exist constants $b,U,c,r_0 \in (0, \infty)$ with $b \le U$ such that
\begin{itemize}
	\item $\forall x \in \ms, \quad b \leq \dens(x) \leq U$;
	\item $\forall r \in (0, r_0], \forall x \in \ms, \quad \Leb_d(\ms \cap B(x,r)) \geq c \Leb_d(B(x,r))$.
\end{itemize}
\end{enumerate}
Assumption~\ref{hyp:density} is related to the \textit{strong density assumption} of \citet{Audibert_Tsybakov07}. It ensures that the density $\dens$ of the measure $\meas$ has a sufficiently regular support and that it is bounded away from zero and infinity.
In view of Remark~\ref{rem:radius}, it is easy to see that \ref{hyp:density} implies \ref{hyp:ballmeasure}: since $\mu(B(x,r)) = \int_{M \cap B(x,r)} \dens(y) \, \diff y$, we have, for $x \in M$ and $r \in (0, r_0]$, the bounds
\begin{align*}
	\mu(B(x,r)) &\le U \Leb_d(B(x,r)) \le U V_d \cdot r^d, \qquad \text{and} \\
	\mu(B(x,r)) &\ge b \Leb_d(\ms \cap B(x,r)) \ge bc \Leb(B(x,r)) \ge bc V_d \cdot r^d,
\end{align*}
where $V_d = \Leb_d(B(0, 1))$ is the volume of the unit ball in $\reals^d$.

Another case of interest is when $\meas$ is supported by a subset of a closed $d$-dimensional Riemannian manifold $(\mani,g)$. Let $\rho_g$ be the geodesic distance on $\mani$ induced by the metric $g$ and let $B_g(x,r) = \cbr{y \in \mani: \rho_g(x,y) \le r}$ denote the geodesic ball with center $x \in \mani$ and radius $r > 0$.
\begin{enumerate}[label=(A\arabic*'')]  
	\item \label{hyp:mani_density} The measure $\meas$ admits a density $\dens$ having support $\ms = \cbr{x \in \mani : f(x) > 0 }$ with respect to the Riemannian volume measure $\vol_g$ on $\mani$. The sectional curvatures of $\mani$ are not greater than some $\delta > 0$ and there exist constants $b,U,c,r_0 \in (0, \infty)$ with $b \le U$ such that
	\begin{itemize}
		\item $\forall x \in \ms, \quad b \leq \dens(x) \leq U$;
		\item $\forall 0 < r \leq r_0, \forall x \in \ms, \quad \vol_g(\mani \cap B_g(x,r)) \geq c \vol_g(B_g(x,r))$.
	\end{itemize}
\end{enumerate}
Since $\meas(B_g(x,r)) = \int_{\ms \cap B_g(x,r)} f(y) \vol_g(\diff y)$, the argument that \ref{hyp:mani_density} implies \ref{hyp:ballmeasure} is the same as the one of the sufficiency of \ref{hyp:density} above, provided we can relate the volume of a geodesic ball to its radius. This is the content of the next result.

\begin{theorem*}[Corollary~2.2 in \cite{kokarev2021}] 
\itshape
Let $(\mani, g)$ be a closed $d$-dimensional Riemannian manifold whose sectional curvatures are not greater than $\delta$, where $\delta > 0$. Let $\inj(\mani)$ denote the injectivity radius of $\mani$. Then, for any point $x \in \mani$ and any positive radius $r\le \rad(g)$, the volume of the geodesic ball $B_g(x,r)$ satisfies the inequalities
\[
	2^{1-d} V_d r^d
	\le \vol_g(B_g(x,r)) 
	\le 2^{d-1} \frac{\vol_g(\mani)}{\rad(g)^d} r^d, 
\]
where $V_d$ is the volume of a unit ball in the d-dimensional Euclidean space, $\rad(g)$ stands for $\min \bigl\{ \inj(\mani), \pi/(2\sqrt{\delta}) \bigr\}$.
\end{theorem*}

In Section~\ref{sec:simu}, we perform numerical experiments featuring manifolds such as the unit sphere $\sphere^{m-1}$ in $\reals^m$ and the orthogonal group. The cited theorem is valid for $\sphere^{m-1}$ since the sectional curvature of a $(m-1)$-sphere of radius $r$ is $1/r^2$ \citep{lee2019introduction}. The orthogonal group is a special case of the Stiefel manifold, which has a positive bounded sectional curvature \citep{ZimmermannStoye}. Other examples are the complex projective space $\operatorname{CP}^d$ and the quaternion projective space $\operatorname{HP}^d$, both having sectional curvatures bounded between $1$ and $4$ \citep{Ichida2000}. 

\subsection{Moments of nearest-neighbor distance}

Consider drawing independent random samples from a distribution on a metric space.
\begin{enumerate}[label=(A\arabic*)]  
	\setcounter{enumi}{\value{assumpnum}}
	\item \label{hyp:iid}  $X, X_1, X_2, \ldots$ are independent and identically distributed random variables in a bounded metric space $(\ms,\rho)$ with common distribution $\meas$. 
	\setcounter{assumpnum}{\value{enumi}}
\end{enumerate}
For Hölder-continuous functions, the quality of the approximation of the integrand $\func$ by its nearest-neighbor estimate $\hat{\func}_n$ depends on the nearest-neighbor distance $\hat{\tau}_n$ (see Section~\ref{subsec:NN}). The next lemma provides control on the $k$-nearest neighbor distance $\hat{\tau}_{n,k}(x)$ and is inspired by the $k$-NN literature \citep{biau2015lectures}. Let $\Gamma$ denote the Euler gamma function.

\begin{lemma}[Bounding moments of nearest neighbor distances]
\label{lem:1NN-distance}
	Under \ref{hyp:ballmeasure} and \ref{hyp:iid}, we have, for any $q > 0$ and for integer $k \in \cbr{1,\ldots,n}$,
	\[
	\expec \sbr{\hat{\tau}_{n,k}(x)^q}
	\le \rbr{C_0 n}^{-q/d} \frac{\Gamma(k+q/d)}{\Gamma(k)}.
	\]
\end{lemma}

\begin{remark}[Lower bounds]
	In Lemma~\ref{lem:1NN-distance}, we only need the lower bound in \ref{hyp:ballmeasure}, not the upper bound. This uniform lower bound in Condition~\ref{hyp:ballmeasure} plays an important role in our analysis as it allows a uniform control on the Voronoi cell diameters. Some refinement might be possible given the recent progress in $k$-NN regression for covariates with unbounded support \citep{kohler2006rates} and $k$-NN classification using some tail assumption on the covariates \citep{gadat2016classification}. Extending the present analysis to such general measures is left for further research.
\end{remark}

\subsection{Root mean squared error bounds}

The main result of this section provides a finite-sample bound on the root mean squared error of the two control neighbors estimates in~\eqref{first_estimate2} and~\eqref{first_estimate}. The bound depends on the regularity of the integrand. 

\begin{enumerate}[label=(A\arabic*)]
	\setcounter{enumi}{\value{assumpnum}}
	\item  \label{hyp:lip} The function $\func:\ms \to \reals$ is Hölder continuous, i.e., there exist $s \in (0, 1]$ and $L>0$ such that
	\[
	\forall x,y \in \ms, \qquad 
	\abs{\func(x) - \func(y)} \leq L \rho(x,y)^s.
	\]
	\setcounter{assumpnum}{\value{enumi}}
\end{enumerate}

The leave-one-out version $\munNNloo$ involves $n$ additional integrals and is computationally cumbersome. The proposed estimate $\munNN$ requires only a single additional integral while remaining close to the leave-one-out estimate, since their difference is

\begin{equation}\label{eq:diffmunloomunn}
\munNN(\func) - \munNNloo(\func)
= \meas(\hat{\func}_n) - \frac{1}{n} \sum_{i=1}^n \meas(\hat{\func}_n^{(i)}).
\end{equation}
Using this property and Lemma~\ref{prop:modification_integral}, we obtain that the root mean squared distance between the leave-one-out version and the proposed estimate is of the order $\Oh(n^{-1/2 - s/d})$ as $n \to \infty$; see Proposition~\propdifflooNN{} in the Supplementary material for a precise statement. Therefore, the two estimates share the same convergence rate. 

In the next theorem, we assume \ref{hyp:ballmeasure}, for which each of \ref{hyp:density} and \ref{hyp:mani_density} is a sufficient condition.

\bigskip

\begin{theorem}[Root mean squared error bounds]
\label{th:cv_rate} 
	Under \ref{hyp:ballmeasure}, \ref{hyp:iid} and \ref{hyp:lip}, if $n\geq 4$, then

	\begin{align*}
		\Exp \left[ \bigl| \munNNloo   (\func)   - \meas(\func) \bigr|^2 \right]^{1/2}    
		&\leq  C_{\mathrm{NN-loo}} n^{-1/2} n^{-s/d}, \\
		\expec \left[ \bigl| \munNN (\func)  -  \meas(\func) \bigr|^2 \right]^{1/2} 
		&\leq C_{\mathrm{NN}} n^{-1/2} n^{-s/d},
	\end{align*}
	where 
	\begin{align*}
		C_{\mathrm{NN-loo}} &= L C_0^{-s/d} (C_1/C_0)^{1/2} 2^{s/d} \sqrt{10 \Gamma(2s/d+2)} , \\
		C_{\mathrm{NN}} &=L C_0^{-s/d} \sqrt{\Gamma(2s/d+2) }  \rbr{2 \sqrt{s/d + 1} + 2^{s/d} (C_1/C_0)^{1/2} \sqrt{10}}.
	\end{align*}
\end{theorem}

For Lipschitz functions ($s = 1$), the rates obtained in Theorem~\ref{th:cv_rate} match the complexity rate stated in \citet{novak2016some}, see Section~\ref{sec:1_intro}. The results in the aforementioned paper are concerned about a slightly more precise context as they assert that no random integration rule (see the article for more details) can reach a higher accuracy---measured in terms of mean-square error---than $\Oh(n^{-1-2/d})$ when the integration measure is the uniform distribution on the unit cube and the integrand is Lipschitz. Theorem~\ref{th:cv_rate} states that the optimal rate is in fact achieved for more general spaces and measures.

The control neighbors estimates involve an approximation of $\func$ by piece-wise constant functions $\hgni$ and therefore cannot achieve an acceleration of the convergence rate of more than $n^{-1/d}$. It is an interesting question whether a higher-order nearest neighbor approximation can achieve a greater accuracy for integrands that are smoother than just Lipschitz. The simulation experiments in Section~\ref{sec:simu} will give some feeling on the actual improvement of the control neighbors estimate over the basic Monte Carlo method at finite sample sizes.

\begin{remark}
\label{rem:bias}
As mentioned in Section~\ref{subsec:CN}, the control neighbor estimate $\munNN$ is in general not unbiased, in contrast to the leave-one-out version $\munNNloo$. Still, by the first statement of Proposition~\propdifflooNN{} in the Supplement, the bias of $\munNN$ is of the order $n^{-1-s/d}$ and thus constitutes an asymptotically negligible part of the estimator's root mean squared error, which decreases with rate $n^{-1/2-s/d}$.
\end{remark}

\subsection{Concentration inequalities}

In order to obtain a finite-sample performance guarantee of the proposed estimates, we apply an extension of McDiarmid's concentration inequality for functions with bounded differences on a high probability set $A$. The inequality is stated in Theorem~3 in Appendix~A of the supplementary material and is itself a minor extension of an inequality due to \cite{extensionofMcDiarmid2015}.

Because of our reliance on theory for nearest neighbor processes in \cite{portier2021nonasymptoticbound}, our bound is established only for measures supported by a compact subset $\ms$ of Euclidean space, which is why we assume the more specific setting of \ref{hyp:density} rather than \ref{hyp:ballmeasure}.

Assumption~\ref{hyp:lip} implies that $\func$ is uniformly bounded on $\ms$, since $\ms$ has a finite diameter. Write $C_\func = \frac{1}{2} \cbr{ \sup_{x \in \ms} \func(x) - \inf_{x \in \ms} \func(x) }$. The two control neighbor estimates satisfy the following concentration inequalities.

\begin{theorem}[Concentration inequalities] \label{th:NNlooconcentration}
	Under \ref{hyp:density}, \ref{hyp:iid} and \ref{hyp:lip}, we have, for any $\varepsilon \in (0,1)$, with probability at least $1-\varepsilon$,

	\begin{align*}
		\abs{\munNNloo(\func) - \meas(\func)} &\le K_1 \Delta_{n,\varepsilon} + \frac{6 C_{\func}  }{n^{1/2+s/d}} \qquad \text{and} \\
		\abs{\munNN(\func) - \meas(\func)} 
		&\le K_2 \Delta_{n,\varepsilon} + \frac{6  C_{\func}}{n^{1/2+s/d}} +  \frac{(s/d+2)  L  (V_d b c)^{-s/d}} {n^{1+s/d}},
	\end{align*}
	where 
	\begin{align*}
		\Delta_{n,\varepsilon}  = \frac{\sqrt{\log(2/\sqrt{\varepsilon})}  + \frac{1}{n^{s/d}}}{n^{1/2+s/d}}
		\cdot
		\begin{dcases}
			\log(24 n/ \varepsilon)^{1+s/d} &\text{ for } \varepsilon \le \frac{2}{n^{1/2+s/d}}, \\
			(3 \log(3 n))^{1+s/d} &\text{ for }  \varepsilon > \frac{2}{n^{1/2+s/d}}.
		\end{dcases}
	\end{align*}
	The values of $K_1, K_2 > 0$ depend on $(b, c, d, L, s, U)$ and are provided explicitly in the proof.
\end{theorem}

For small but fixed $\varepsilon \in (0, 1)$, the concentration inequalities state that the probability that the error $\abs{\hat{\mu}_n(\func) - \meas(\func)}$ is larger than a constant multiple of $\sqrt{\log(1/\varepsilon)} (\log n)^{1 + s/d} n^{-1/2}n^{-s/d}$ is less than $\varepsilon$. The additional logarithmic factor in Theorem~\ref{th:NNlooconcentration} in comparison to Theorem~\ref{th:cv_rate} comes from the need to exert uniform control on the largest nearest neighbor distance via Lemma~5 in Appendix~A in the supplement, which in turn originates from a concentration inequality in \cite{portier2021nonasymptoticbound}.

\vspace{-0.1cm} 
\subsection{Optimal $L_p$-approximation with nearest neighbors} 
\label{subsec:Lpapprox}

The upper bounds on nearest neighbor distances stated in Lemma~\ref{lem:1NN-distance} have been key to examine the control neighbor estimators  $\munNNloo(\func) $ and $\munNN(\func) $ as claimed in Theorems~\ref{th:cv_rate} and~\ref{th:NNlooconcentration}. Another application of these upper bounds is the $L_p$-approximation of the nearest neighbor function estimator $\hat \func_{n,k}(x) = k^{-1} \sum_{i=1}^k \func\bigl(\hat{N}_{n,i}(x)\bigr)$ of $ \func(x)$. Note that this type of estimator is often used in the setting of noiseless regression where one has access to $\func(X_i)$ without noise. A precise statement is given in the next proposition, the proof of which is given in the supplement.
   
\begin{proposition}[$L_p$-approximation with nearest neighbors]
	\label{lem:Lpconv}
	Under \ref{hyp:ballmeasure}, \ref{hyp:iid}, and \ref{hyp:lip}, we have, for any
	\vspace{-0.1cm} $k \in \cbr{1,\ldots,n}$ and $p\geq 1$,
	\[
		\forall x\in \ms, \qquad 
		\expec \sbr{
			\abs{\hat \func_{n,k} (x) - \func(x)}^p
		}^{1/p} 
		\leq  L (C_0 n)^{-s/d} \rbr{\frac{\Gamma(k + sp/d)}{\Gamma(k)}}^{1/p}.
	\]
\end{proposition}

Since $\Gamma(x+\alpha) \sim x^\alpha \Gamma(x)$ as $x \to \infty$ for fixed $\alpha \in \reals$, the upper bound is of the order $(k/n)^{s/d}$ as $k = k_n \to \infty$. We emphasize that similar results as in Proposition~\ref{lem:Lpconv} can be found among the noiseless regression literature. However, to the best of our knowledge, most of the results seem to have been obtained for Euclidean covariates, i.e., when $M = \mathbb R^d$. In Chapter 15 in \cite{biau2015lectures}, the weak convergence of $n^{1/d} (\hat \func_{n,k} (x) - \func(x)) $ is established. In Theorem 1 from \cite{kohler2013optimal}, the same rate of convergence for the integrated error is obtained.
Some other results are given in \cite{hinrichs2020power,krieg+novak:2022} where the authors investigate whether the distribution of the points can deteriorate the $L_p$-error or not. They show that the optimal rate of convergence can be obtained using either uniformly distributed points or deterministic points from the uniform grid. Our Proposition~\ref{lem:Lpconv} confirms, not only for the Euclidean space but for general metric spaces, the following principle, that when there are enough points everywhere, as implied by \ref{hyp:ballmeasure}, the optimal rates of convergence for the $L_p$-error can be achieved. 

Even though the analysis of the control neighbor estimate cannot be deduced from Proposition \ref{lem:Lpconv} because of the probabilistic dependence that is introduced in the algorithm (see the proofs of Theorems~\ref{th:cv_rate} and~\ref{th:NNlooconcentration}), the previous statement can be used for a better understanding of  the control neighbor estimator as it shows that the rate obtained in Theorem~\ref{th:cv_rate} can be decomposed as 
the usual Monte Carlo error, $n^{-1/2}$, times the $L_p$-approximation error of the control variates estimating the integrand, $n^{-s/d}$. This is in line with the results from \cite{portier+s:2018} where the same type of decomposition can be found for another control variates based estimator.

\vspace{-0.1cm}
\section{Numerical experiments}\label{sec:simu}

In Section \ref{subsec:num_details}, we start by presenting the \textit{control neighbors} algorithm\footnote{For reproducibility, the code is available at \href{https://github.com/RemiLELUC/ControlNeighbors}{Control Neighbors}.} along with some practical remarks related to its implementation. We present in Section~\ref{subsec:toy_data} some examples involving integration problems on different spaces. Finally, Section~\ref{subsec:OT} deals with the application of Monte Carlo estimates for optimal transport when computing the Sliced-Wasserstein distance. Additional experiments on a finance application with Monte Carlo exotic option pricing are available in Appendix E. 

\vspace{-0.25cm}

\subsection{Practical implementation}
\label{subsec:num_details}
\vspace{-0.25cm}

The procedure for computing the control neighbor estimate $\munNN(\func)$ is presented in Algorithm~\ref{algo:cvnn}. Required is a collection of points $X_{1},\dots,X_{n}$ generated from a distribution $\mu$. The estimate is based on the evaluations $\func(X_i)$ of the integrand and the evaluations $\hgni (X_i)$ of the leave-one-out $1$-nearest neighbors ($1$-NN) estimates of $\func(X_i)$. Here, $\hgni$ denotes the $1$-NN estimate of $\func$ constructed from the sample $\cX_n^{(i)} = \{X_1,\ldots, X_n\} \setminus \{X_i\}$ without the $i$-th particle, for $i \in \{1,\ldots,n\}$; see Section~\ref{subsec:NN} for precise definitions. Also required is the integral $\meas(\hat \func_n)$ of the $1$-NN estimate $\hat \func_n$ based on the whole sample. Practical remarks regarding the computation of all these quantities are given next.
 
\begin{algorithm}[h]
	\linespread{1.25}
	\caption{Control Neighbors for Monte Carlo integration}\label{alg:cvnn}
	\begin{algorithmic}[1]
		\Require integrand $\func$, probability measure $\meas$, number of particles $n$.
		\State Generate an independent random sample $X_{1},\dots,X_{n}$ from $\meas$ 
		\State Compute evaluations $\func(X_1), \ldots, \func(X_n)$
		\State Compute nearest neighbor evaluations $\hat{\func}_{n}^{(1)}(X_1), \ldots, \hat{\func}_{n}^{(n)}(X_n)$
		\State Compute the integral of the nearest neighbor estimate  $\meas(\hat \func_n)$ 
		\State Return $\frac{1}{n} \sum_{i=1}^n [\func(X_i) - \{\hgni (X_i) -   \meas ( \hat \func_n  ) \} ] $
	\end{algorithmic}
	\label{algo:cvnn}
		\vspace{-0.15cm}
\end{algorithm}

	\vspace{-0.2cm}
\mypar{Nearest neighbors and tree search.}
	The naive neighbor search implementation involves the computation of distances between all pairs of points in the training samples and may be computationally prohibitive. To address such practical inefficiencies, a variety of tree-based data structures have been invented to reduce the cost of a nearest neighbors search. The KD-Tree \citep{bentley1975multidimensional} is a binary tree structure which recursively partitions the space along the coordinate axes, dividing it into nested rectangular regions into which data points are filed. The construction of such a tree requires $\Oh(dn \log n) $ operations \citep{friedman1977algorithm}. Once constructed, the query of a nearest neighbor in a KD-Tree can be done in $\Oh(\log n) $ operations. However, in high dimension, the query cost increases and the structure of Ball-Tree \citep{omohundro1989five} is favoured. Where KD trees partition data points along the Cartesian axes, Ball trees do so in a series of nested hyper-spheres, making tree construction more costly than for a KD tree, but resulting in an efficient data structure even in high dimensions. Many software libraries contain KD-tree and Ball-Tree implementations with efficient compression and parallelization \citep{pedregosa2011scikit,johnson2019billion}.

\mypar{Evaluation of $\hgni (X_i)$.}
When the computing time is measured through the evaluation of the integrand, the $n$ additional evaluations $\hgni (X_i)$  are not computationally difficult as no additional calls to $\func$ are necessary.
The leave-one-out nearest neighbor evaluation can be easily obtained as follows: first fit a KD-Tree  on the particles $X_1,\ldots, X_n$, then query the $2$-nearest neighbor of each $X_i$ to produce the vector of leave-one-out nearest neighbors.

\mypar{Evaluation of $\meas(\hat \func_n)$ and Voronoi volumes.}
The quantity $\meas(\hat \func_n)$ is the sum of the evaluations $\func(X_i)$ weighted by the values of the Voronoi volumes associated to the sample points $X_i$ (see Definition~\ref{def:vor_cells} in the next section). 	The Voronoi volumes may be hard to compute but can always be approximated. The integral $\meas(\hat \func_n)$ of the nearest neighbor estimate may be replaced by a Monte Carlo approximation based on $N$ particles, such as $\tilde{\meas}_N(\hat \func_n) = N^{-1} \sum_{ i = 1} ^N \hat \func_n(\tilde X_i)$, where the variables $\tilde X_i$ are drawn independently from $\meas$. No additional evaluations of $\func$ are required.
Conditionally on the first sample $X_1,\ldots,X_n$, the error of this additional Monte Carlo approximation is  $\Oh(N^{-1/2})$, meaning that large values of the form $N=n^{1+2/d}$ permit to preserve the $\Oh(n^{-1/2} n^{-1/d})$ rate of the control neighbors estimate for Lipschitz functions.  

\begin{figure}[H]
	\vspace{-0.5cm}
	\centering
	\begin{subfigure}[h]{0.44\linewidth}
		\centering
		\includegraphics[width=\linewidth]{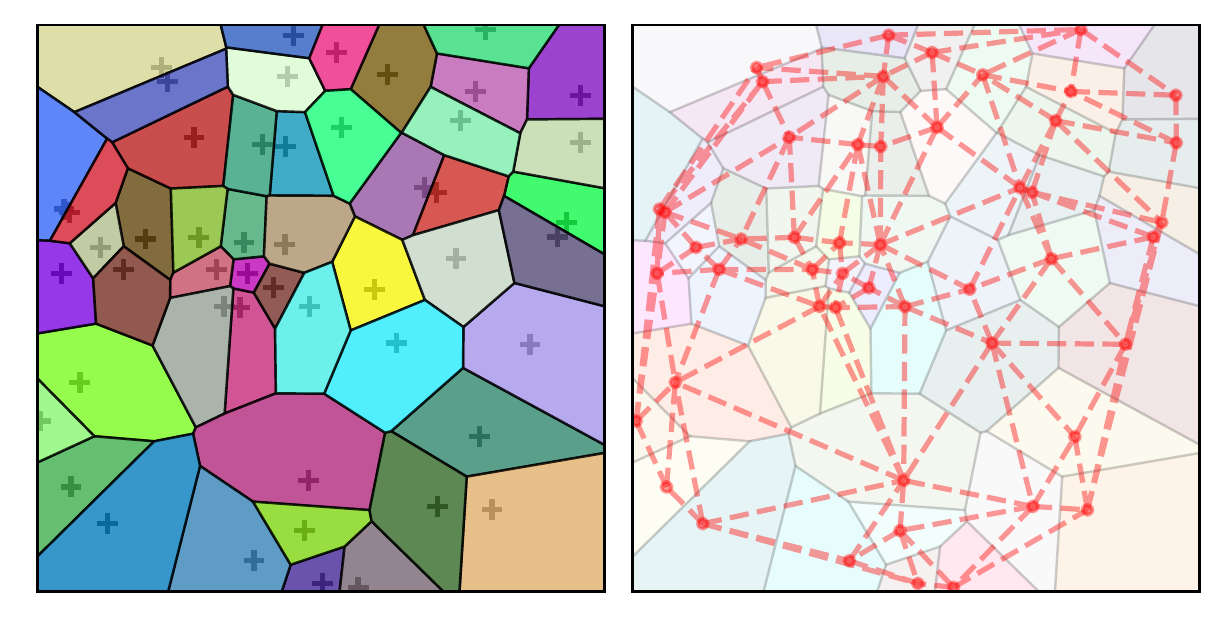}
		\label{fig:voro2d}
	\end{subfigure}
	\begin{subfigure}[h]{0.55\linewidth}
		\centering
		\includegraphics[width=\linewidth]{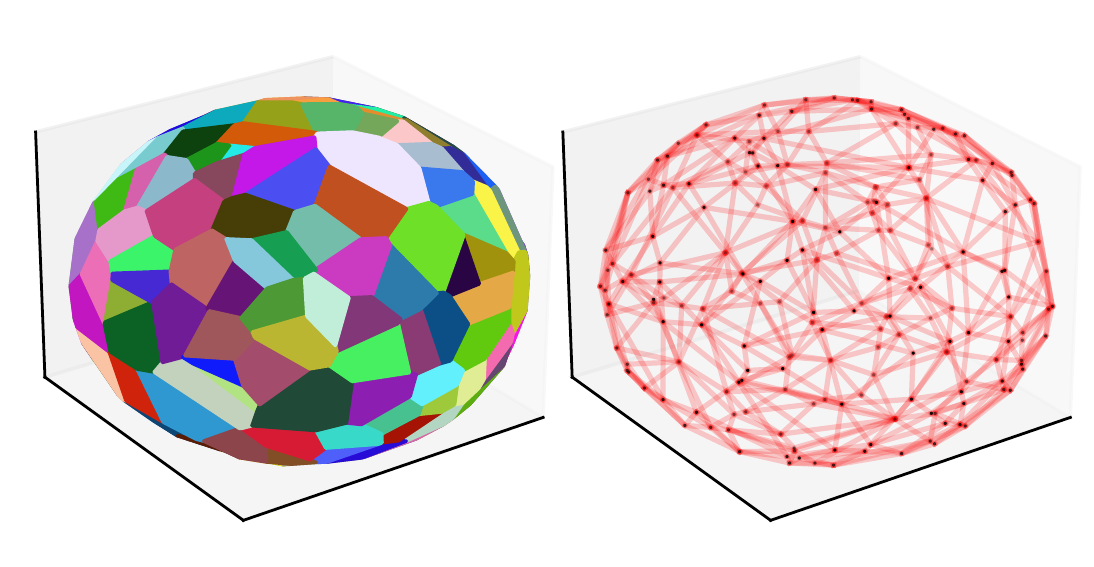}
		\label{fig:voro3d}
	\end{subfigure}
	\vspace{-1.cm}
	\caption{Voronoï cells and Delaunay triangulations on $[0,1]^2$ and $\sphere^2$. \\ (Made with Python package \textsf{Matplotlib} \citep{hunter2007matplotlib})}
	\label{fig:voro} 
\end{figure}

In case the measure $\meas$ is the uniform measure on $[0,1]^d$, one may be able to explicitly compute the Voronoi volumes. Starting from the pioneering work of \citet{richards1974interpretation} in the context of protein structures, there have been advances to perform efficient Voronoi volume computations using Delaunay triangulations and taking advantage of graphic hardware \citep{hoff1999fast}. Computations for Voronoi tessellations in dimensions $d = 2$ and $d = 3$ are implemented in the Voro++ software \citep{rycroft2009voro++}. However, this type of algorithm might become inefficient when $d$ is large.

\mypar{Computing time.}
When using the standard KD-tree approach, the computation time of our algorithm can be estimated in the light of the previous remarks. The different computation times are: $\Oh(nd \log n)$ for building the KD-tree, $\Oh(n \log n)$ for the evaluations $\{\hgni (X_i)\}_{i=1}^n$ and finally, $\Oh(N \log n)$ operations for the estimation of $\meas(\hat \func_n)$. Choosing $N$ as recommended before, the overall complexity is $\Oh\bigl( (nd + n^{1+2/d}) \log n \bigr)$ operations. Note that when the Voronoi volumes are available, the computation time reduces to $\Oh(nd  \log n)$.

\mypar{Extensions.}
A natural variant of the proposed method is obtained by replacing the $1$-NN estimate $\hat \func_n$ in Eq.~\eqref{first_estimate2} by a $k$-NN estimate $\hat \func_n^{(k)}$ which averages the evaluations of the $k$ nearest neighbors of a given point. The estimate is then defined by $\hat \func_n^{(k)}(x) = k^{-1}\sum_{j=1}^k \func(\hat N_{n,j}(x))$ where $\hat N_{n,j}(x)$ is the $j$-nearest neighbor of $x$. This involves both the tuning of the hyperparameter $k \geq 1$ and some extra computation due to the associated nearest neighbors search. In regression or classification, high values of $k$ can reduce the variance of the estimate by averaging the model noise at the cost of added computations. In contrast, the control neighbors estimate $(k=1)$ is free of these additional costs and takes advantage of the noiseless evaluations \citep[Chapter~15]{biau2015lectures} of the integrand.

\subsection{Integration on various spaces: $[0,1]^d, \Rd, O_m(\rset)$ and $\sphere^{q-1}$}
\label{subsec:toy_data}

The aim of this section is to empirically validate the $\Oh(n^{-1/2} n^{-1/d})$ convergence rate of the control neighbors estimate in a wide variety of integration problems ranging from the unit cube $[0,1]^d$, the group of orthogonal matrices $O_m(\rset)$, which is an $m(m-1)/2$ dimensional manifold embedded in $\reals^{m \times m}$, and the unit sphere $\sphere^{q-1}$ in $\reals^q$. In the different settings, the sample size evolves from $n=10^1$ to $n=10^4$. The figures report the evolution of the root mean squared error $n \mapsto \expec [ | \munNN  (\func) - \meas (\func) |^2]^{1/2}$, where the expectation is computed over $100$ independent replications. In all experiments, MC represents the naive Monte Carlo estimate and CVNN returns the value of $\munNN (\func)$ for which the integral $\meas(\func_n)$ is replaced by a Monte Carlo estimate that uses $n^2$ particles, with $n$ the number of evaluations of $\func$. 

\begin{figure}
	\centering
	\begin{subfigure}[h]{0.32\linewidth}
		\centering
		\includegraphics[width=\linewidth]{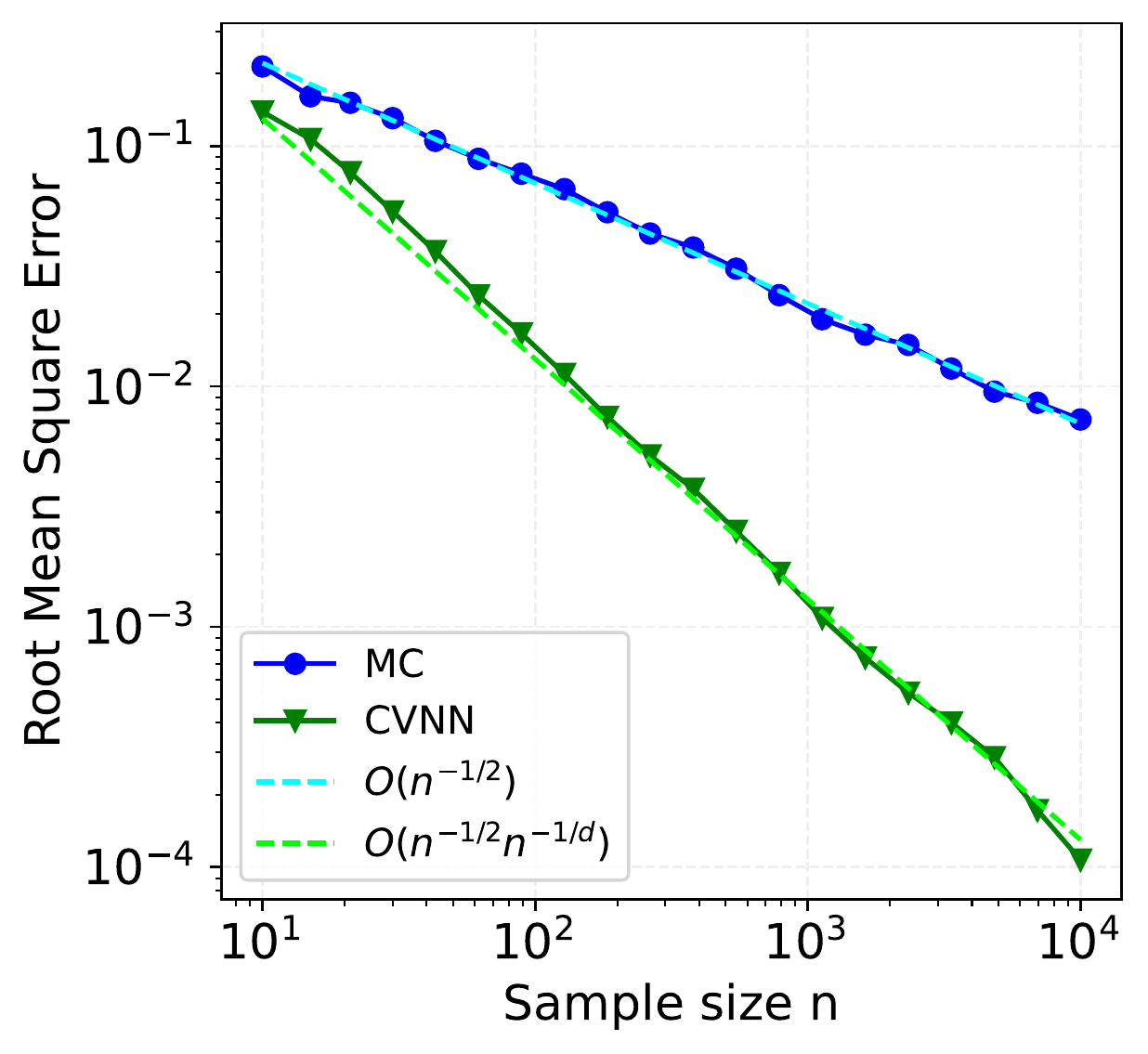}
		\caption{$\func_1, d=2$}
		\label{fig:func1_d2}
	\end{subfigure}
	\hfill
	\begin{subfigure}[h]{0.32\linewidth}
		\centering
		\includegraphics[width=\linewidth]{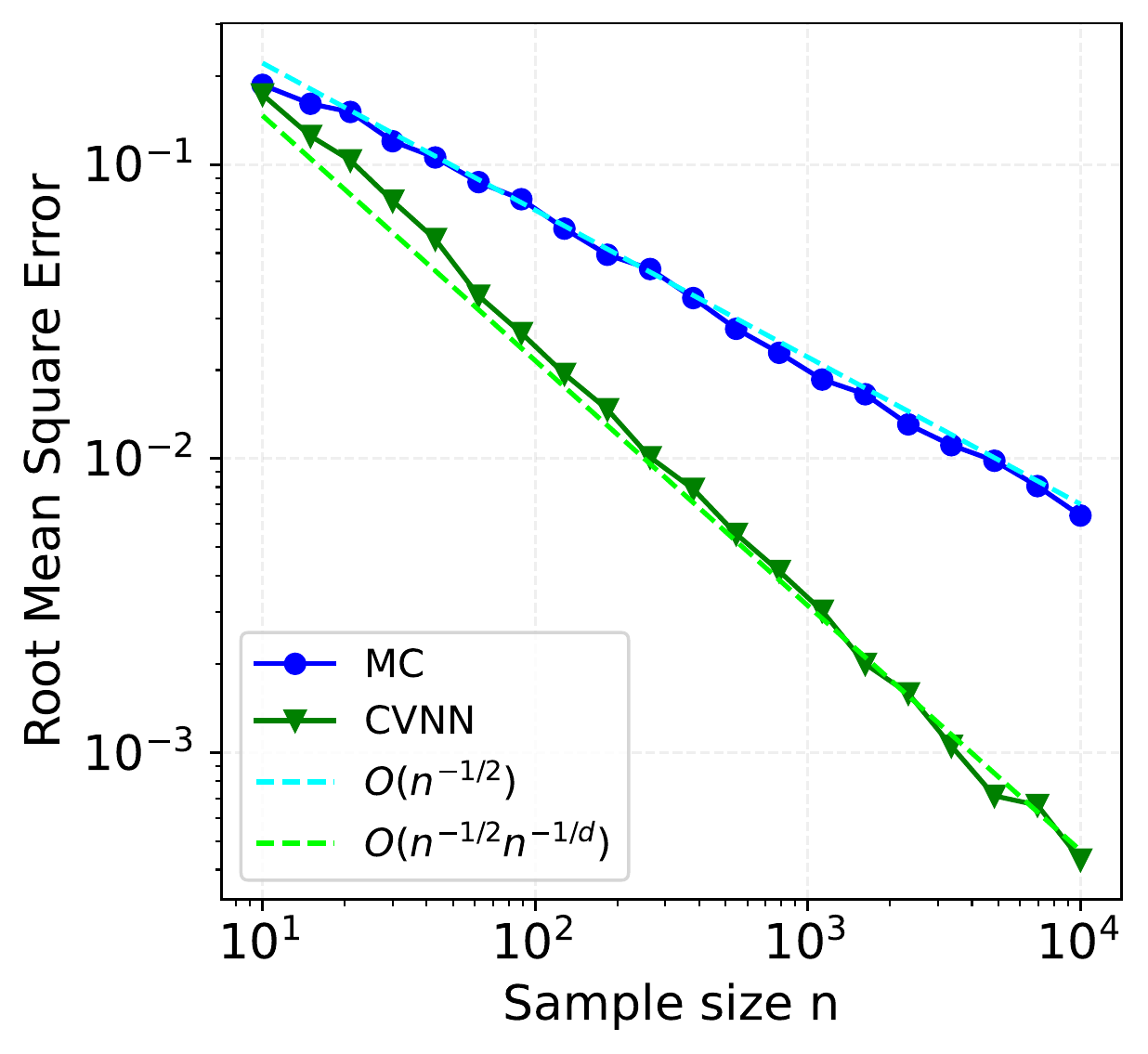}
		\caption{$\func_1, d=3$}
		\label{fig:func1_d3}
	\end{subfigure}
	\hfill
	\begin{subfigure}[h]{0.32\linewidth}
		\centering
		\includegraphics[width=\linewidth]{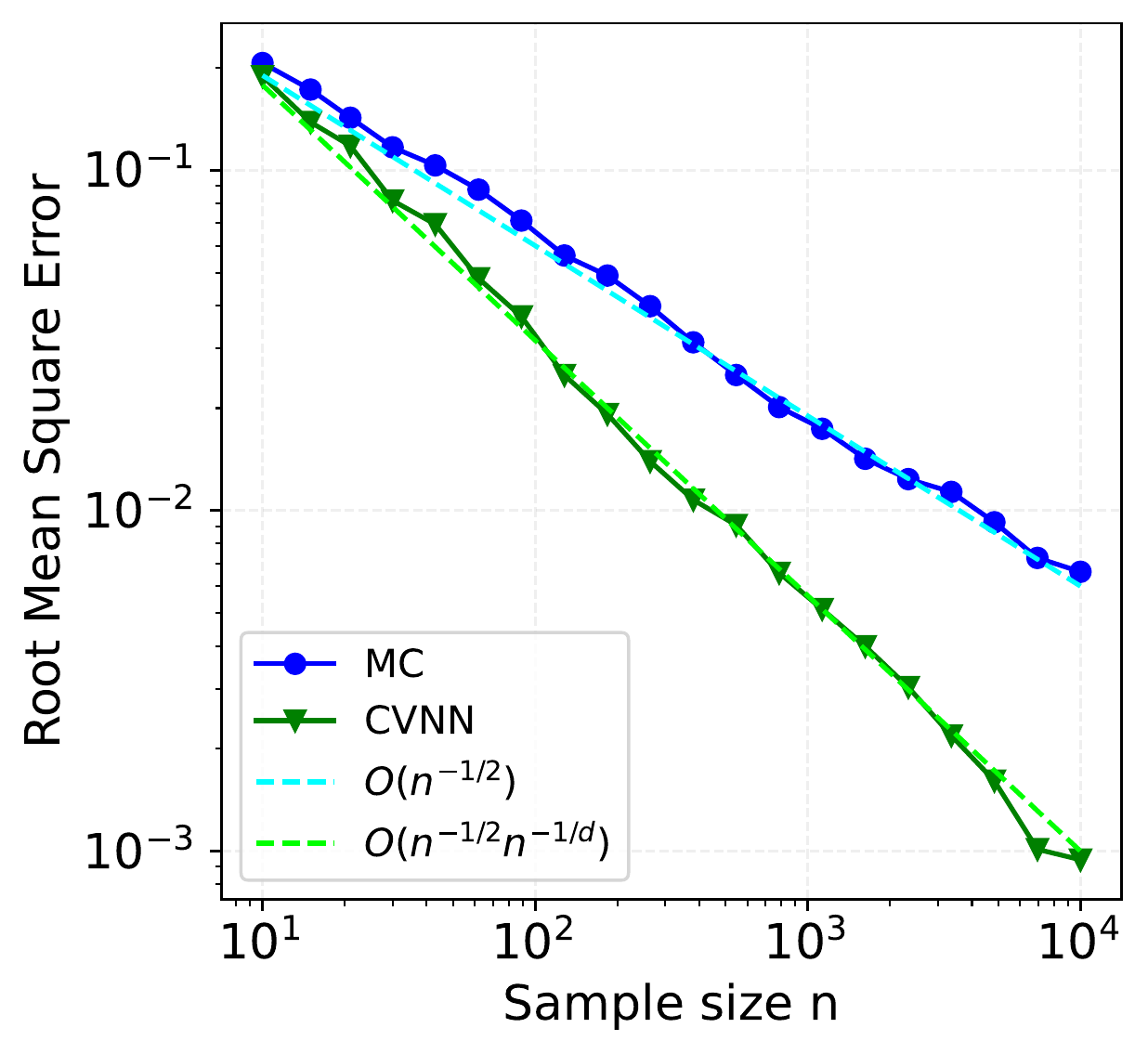}
		\caption{$\func_1, d=4$}
		\label{fig:func1_d4}
	\end{subfigure}
	\hfill
	\begin{subfigure}[h]{0.32\linewidth}
		\centering
		\includegraphics[width=\linewidth]{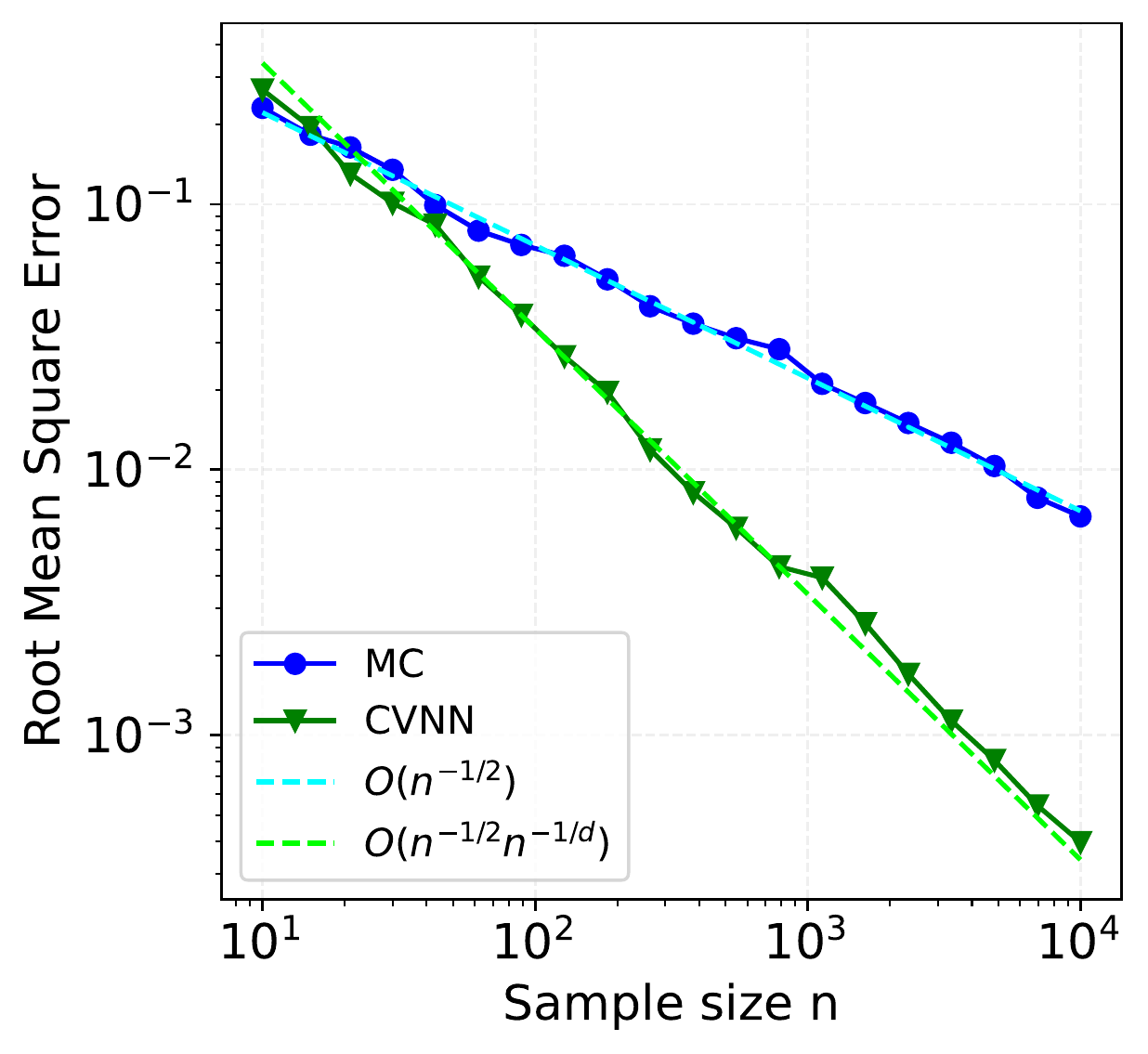}
		\caption{$\func_2, d=2$}
		\label{fig:func2_d2}
	\end{subfigure}
	\hfill
	\begin{subfigure}[h]{0.32\linewidth}
		\centering
		\includegraphics[width=\linewidth]{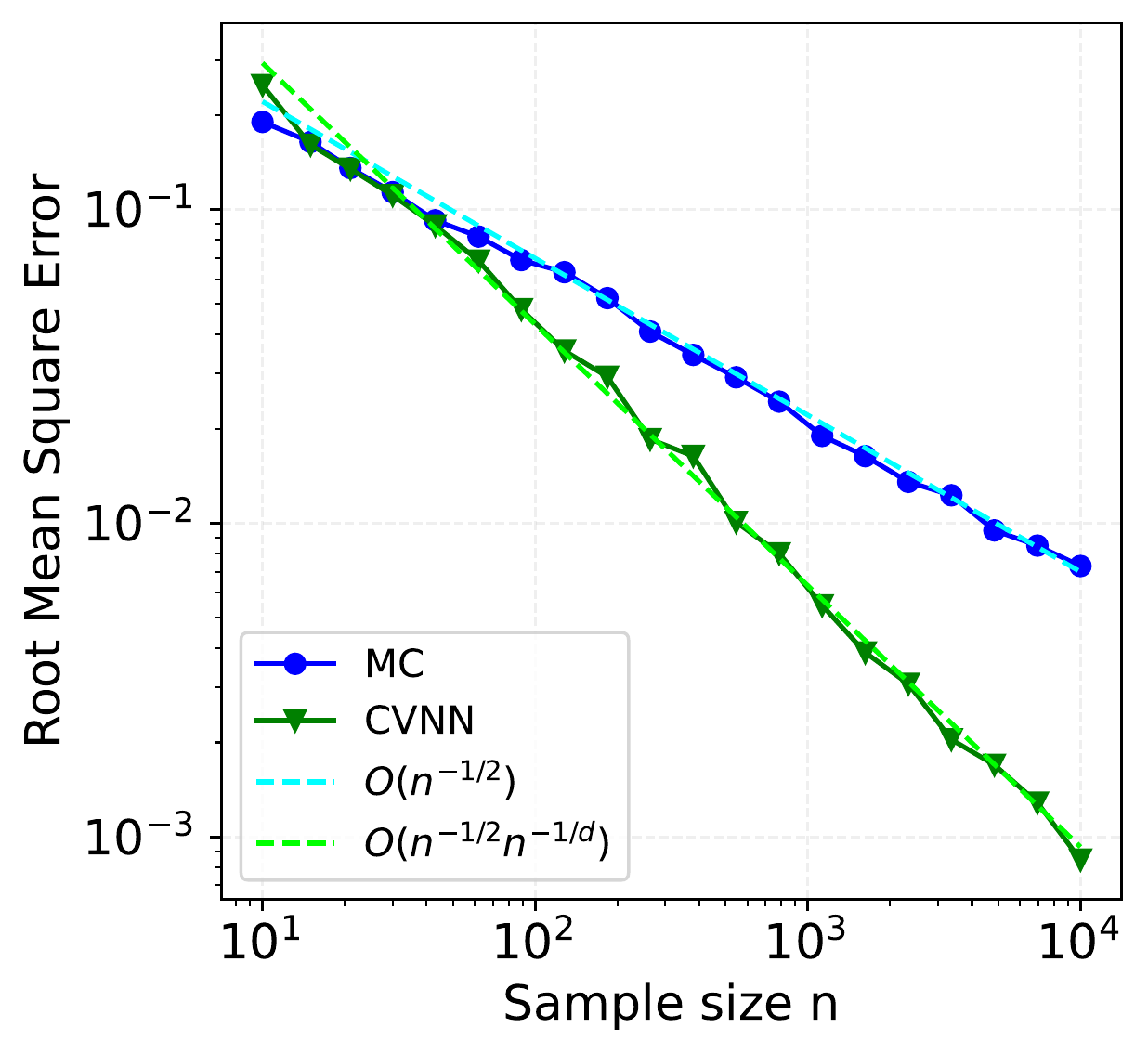}
		\caption{$\func_2, d=3$}
		\label{fig:func2_d3}
	\end{subfigure}
	\hfill
	\begin{subfigure}[h]{0.32\linewidth}
		\centering
		\includegraphics[width=\linewidth]{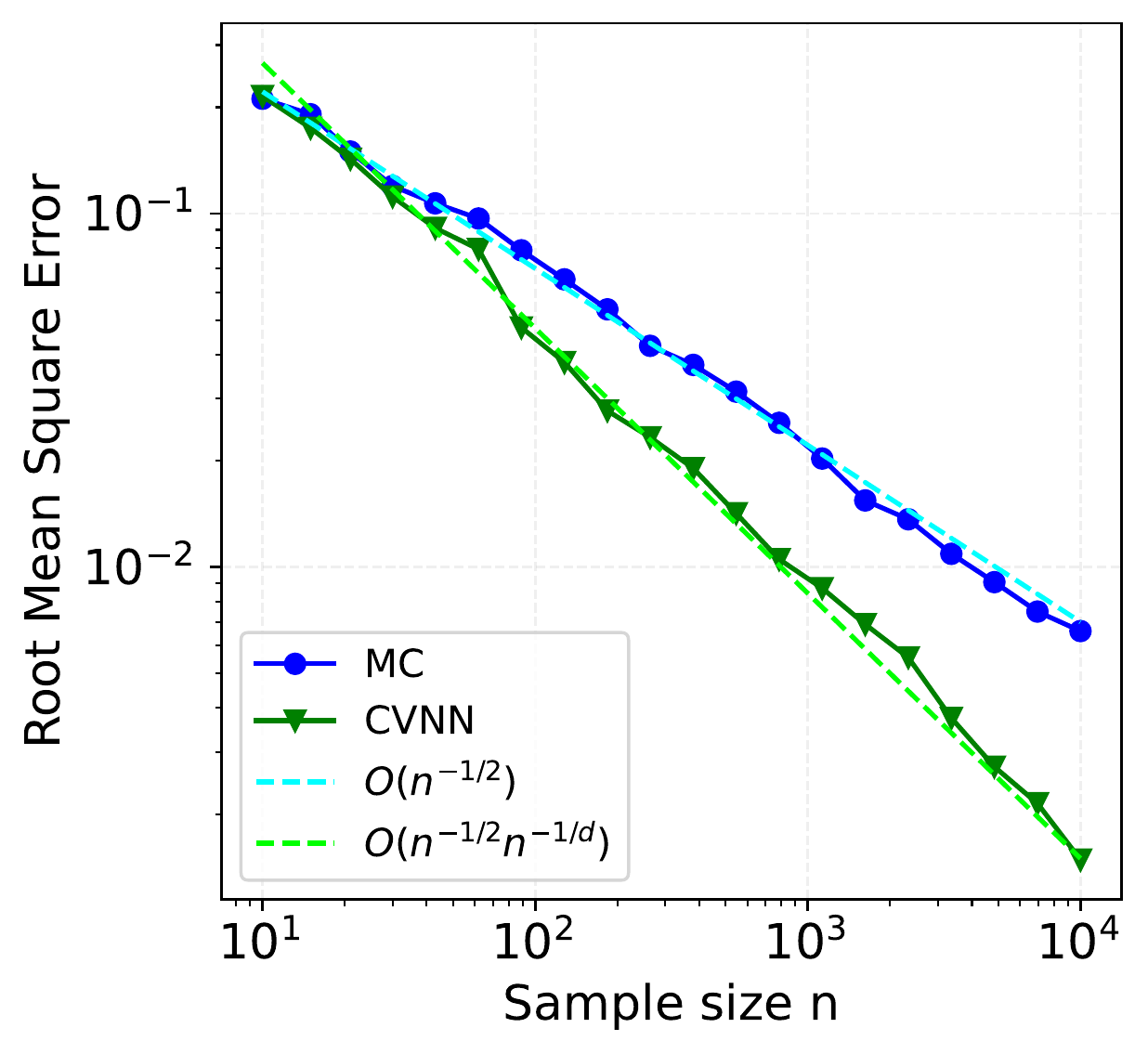}
		\caption{$\func_2, d=4$}
		\label{fig:func2_d4}
	\end{subfigure}
	\caption{Root mean squared errors obtained over $100$ replications for functions $\func_1$ (top) and $\func_2$ (bottom) in Eq.~\eqref{eq:phi1phi2} in dimension $d \in \{2,3,4\}$ (left to right).}
	\label{fig:sin_func} 
\end{figure}

\mypar{Integration on $[0,1]^d$ and $\Rd$.} Consider first the integration problem $\meas(\func) = \int \func  \,  \diff \meas$ where the measure $\meas$ is either the uniform distribution over the unit cube $[0,1]^d$ or the multivariate Gaussian measure on $\Rd$. The latter case is not covered by our theory, since the support is unbounded.
The goal is to compute $\int \func_1(x) \1_{[0,1]^d}(x) \, \diff x$ and $\int \func_2(x) \, \psi(x) \, \diff x$ with
\begin{equation}
	\label{eq:phi1phi2}
	\textstyle
	\func_1(x_1,\ldots,x_d) 
	=  \sin\left(\pi \left(\frac{2}{d}\sum_{i=1}^d x_{i}-1\right)\right)
	\quad \text{ and } \quad
	\func_2(x_1,\ldots,x_d) 
	=  \sin\left(\frac{\pi}{d} \sum_{i=1}^d x_{i}\right)
\end{equation}
and with $\psi(\cdot)$ the probability density function of the multivariate Gaussian distribution $\mathcal{N}(0,I_d)$. Figure~\ref{fig:sin_func} displays the evolution of the root mean squared error for the two integrals in dimensions $d \in \{2, 3, 4\}$.
The different error curves confirm the optimal convergence rate $\Oh(n^{-1/2} n^{-1/d})$ for the control neighbors estimate and illustrate the error reduction with respect to the Monte Carlo method.

\begin{figure}
	\centering
	\begin{subfigure}[h]{0.24\linewidth}
		\centering
		\includegraphics[width=\linewidth]{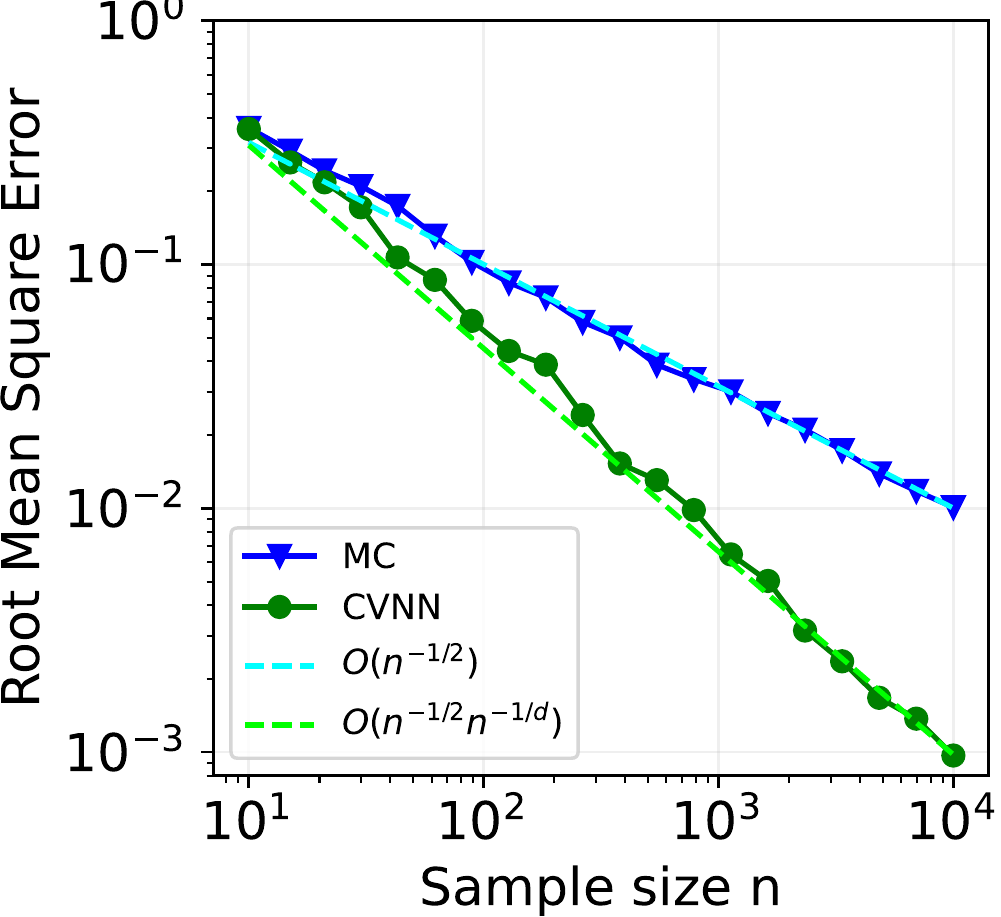}
		\caption{RSME for $a_1$}
		\label{fig:func1_d2}
	\end{subfigure}
	\begin{subfigure}[h]{0.24\linewidth}
		\centering
		\includegraphics[width=\linewidth]{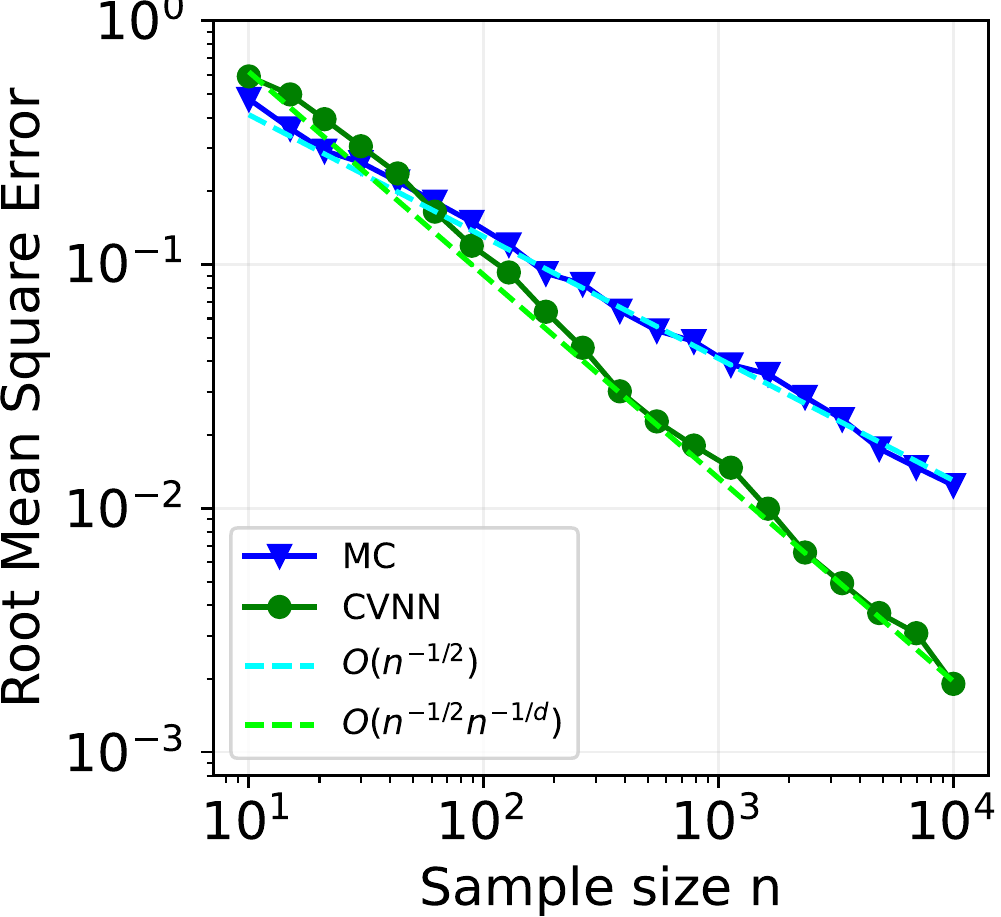}
		\caption{RMSE for $a_2$}
		\label{fig:func1_d3}
	\end{subfigure} 
	\begin{subfigure}[h]{0.25\linewidth}
		\centering
		\includegraphics[width=\linewidth]{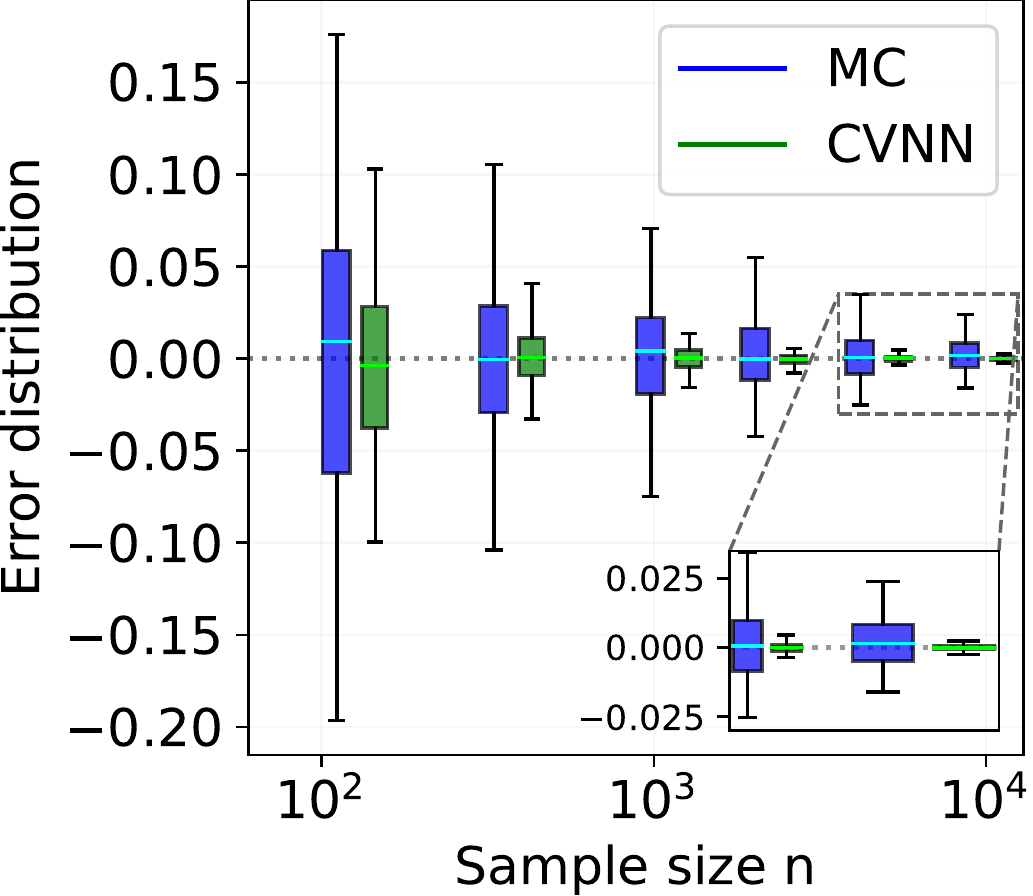}
		\caption{Boxplot for $a_1$}
		\label{fig:func1_d4}
	\end{subfigure}
	\begin{subfigure}[h]{0.25\linewidth}
		\centering
		\includegraphics[width=\linewidth]{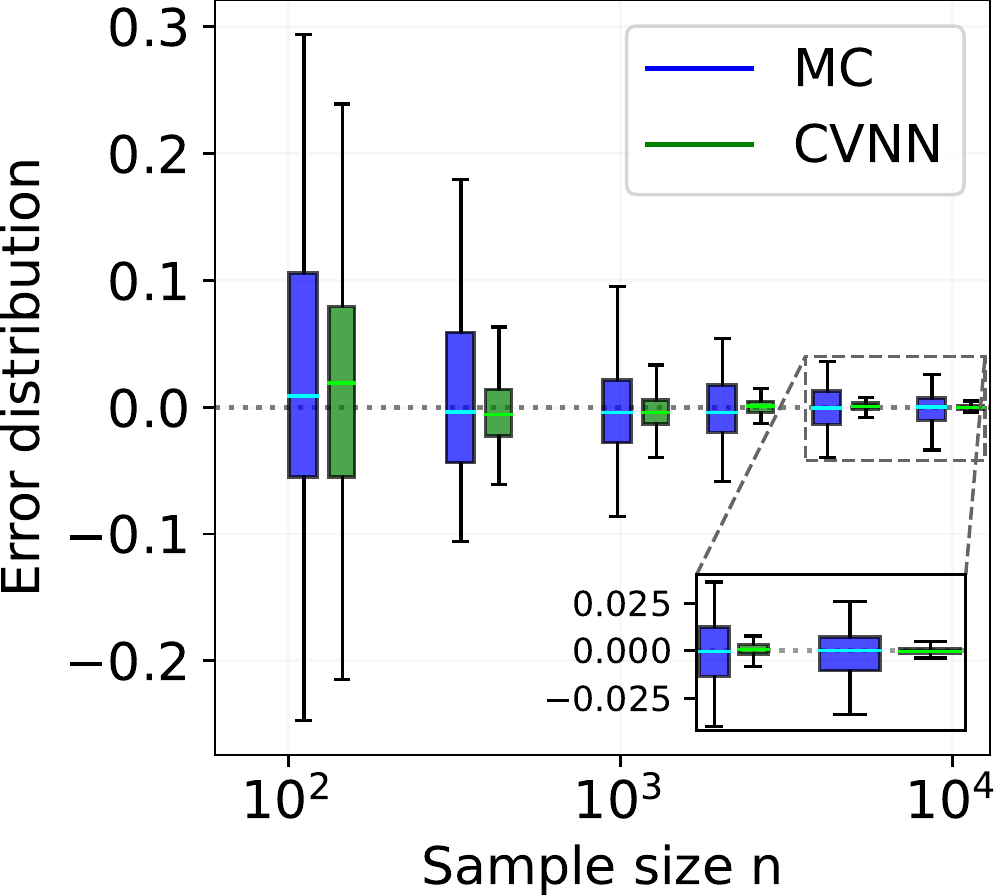}
		\caption{Boxplot for $a_2$}
		\label{fig:func2_d2}
	\end{subfigure}
	\caption{Root mean squared errors (left) and boxplots of the errors (right) obtained over $100$ replications for integrands $X \mapsto \operatorname{tr}(X)^k$ for $k \in \{1,2\}$ in Eq.~\eqref{eq:Odphik} with respect to the unit Haar measure on the orthogonal group $O_3(\rset)$.}
	\label{fig:trace_moments} 
\end{figure}

\mypar{Integration on the orthogonal group $O_m(\reals)$.}
Consider the group of real orthogonal matrices $O_m(\reals) = \cbr{X \in \mathrm{GL}_m(\reals): X^\top X = X X^\top = \mathrm{I}_m}$. As a toy example, we compute moments of the trace of a random orthogonal matrix, i.e., integrands of the form $O_m(\rset) \to \rset : X \mapsto \operatorname{tr}(X)^k$, when $X$ is drawn uniformly on $O_m(\reals)$. In fact, exact expressions are known; see for instance \cite{diaconis2001linear} and \cite{pastur2004moments} for background and results. The goal is thus to compute the integrals 
\begin{equation}
	\label{eq:Odphik}
	a_k 
	= \int_{O_m(\reals)} \operatorname{tr}(X)^k \, \diff\meas(X),
\end{equation}
with $\diff\meas(X)$ the unit Haar measure on $O_m(\reals)$.
Assumption~\ref{hyp:ballmeasure} is satisfied with $d = m(m-1)/2$. In practice, random orthogonal matrices $X_1,\ldots,X_n \in O_m(\rset)$ are generated using the function \textsf{ortho\_group} of the Python package \textsf{scipy} \citep{virtanen2020scipy}. This function returns random orthogonal matrices drawn from the Haar distribution using a careful QR decomposition\footnote{QR decomposition refers to the factorization of a matrix $X$ into a product $X=QR$ of an orthonormal matrix $Q$ and an upper triangular matrix $R$.} 
as in \citet{mezzadri2007generate}.  The nearest neighbors are computed using the norm associated to the Frobenius inner product $\langle X, Y \rangle = \operatorname{tr}(X^\top Y)$ for $X,Y \in O_m(\rset)$.
Figure~\ref{fig:trace_moments} reports the evolution of the root mean squared error and boxplots of the integration error over the group $O_3(\rset)$ which corresponds to dimension $d=3$. 
Both integrands are continuously differentiable and thus Lipschitz, since $O_m(\reals)$ is compact.
Once again, the experiments empirically validate the convergence rates of the Monte Carlo methods and reveal the variance reduction obtained with the control neighbors estimate.

\begin{figure}
	\centering
	\begin{subfigure}[h]{0.32\linewidth}
		\centering
		\includegraphics[width=\linewidth]{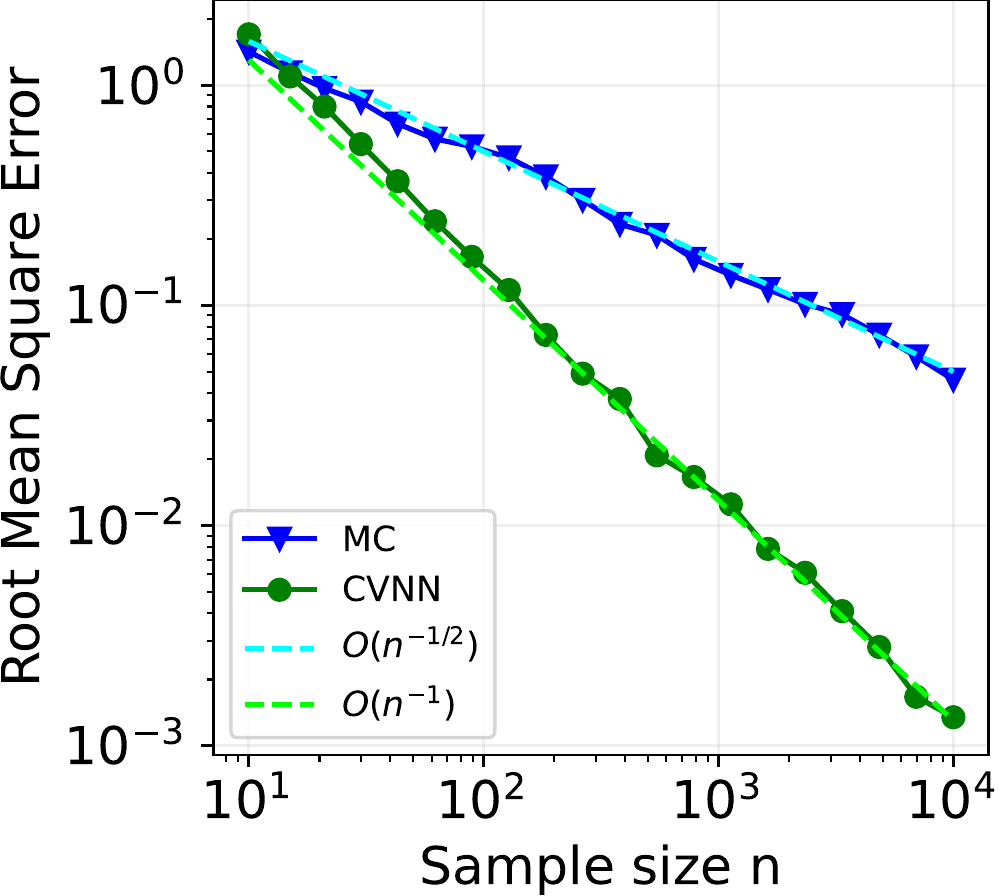}
		\caption{RMSE $\func_3$}
		\label{fig:func3}
	\end{subfigure}
	\hfill
	\begin{subfigure}[h]{0.32\linewidth}
		\centering
		\includegraphics[width=\linewidth]{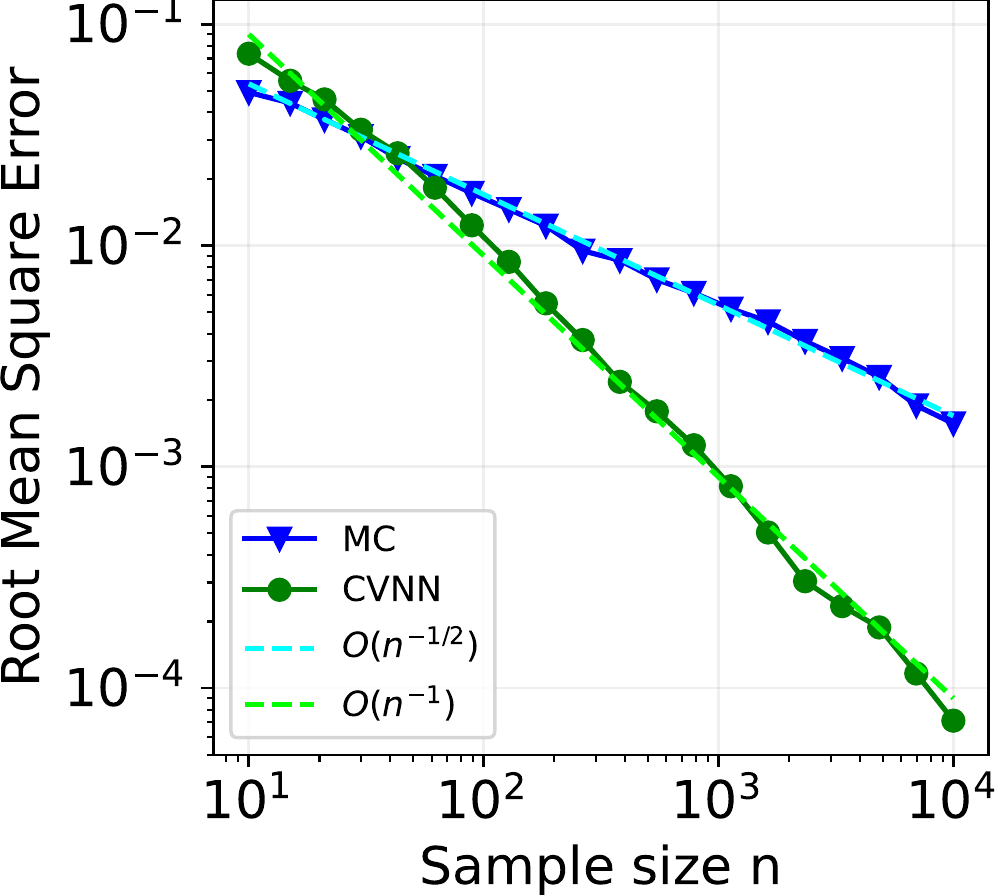}
		\caption{RMSE $\func_4$}
		\label{fig:func4}
	\end{subfigure}
	\hfill
	\begin{subfigure}[h]{0.32\linewidth}
		\centering
		\includegraphics[width=\linewidth]{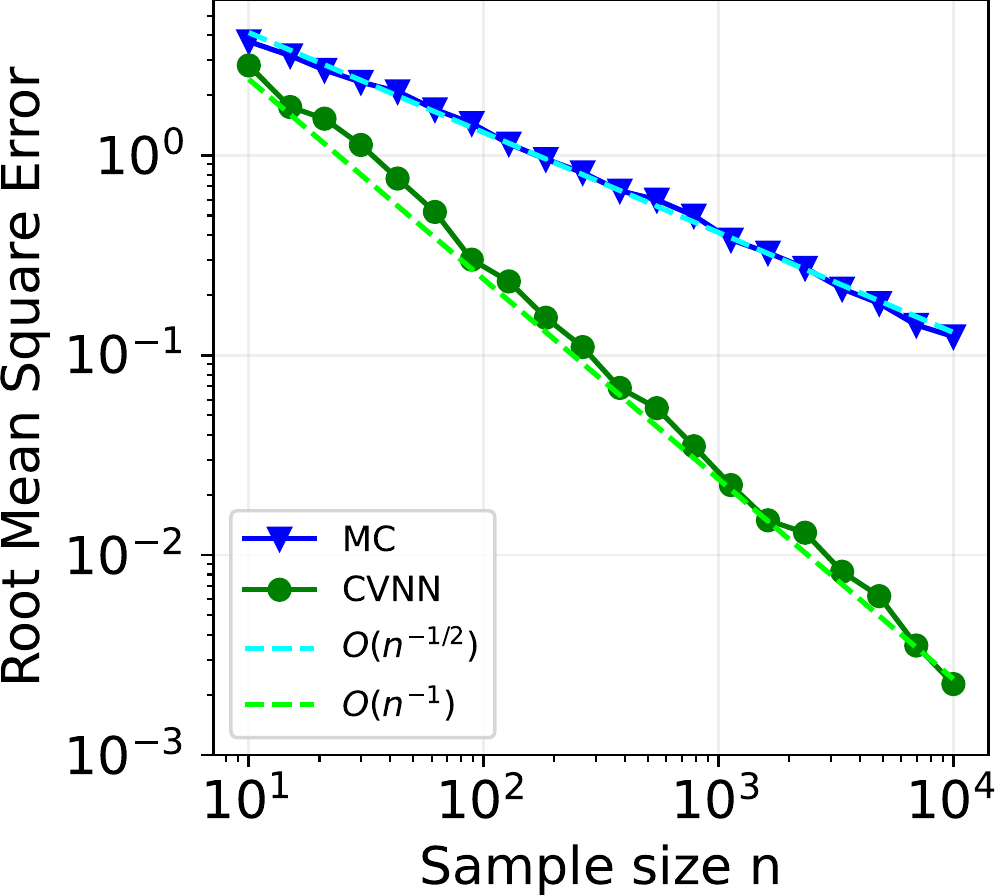}
		\caption{RMSE $\func_5$}
		\label{fig:func5}
	\end{subfigure}
	\hfill
	\begin{subfigure}[h]{0.32\linewidth}
		\centering
		\includegraphics[width=\linewidth]{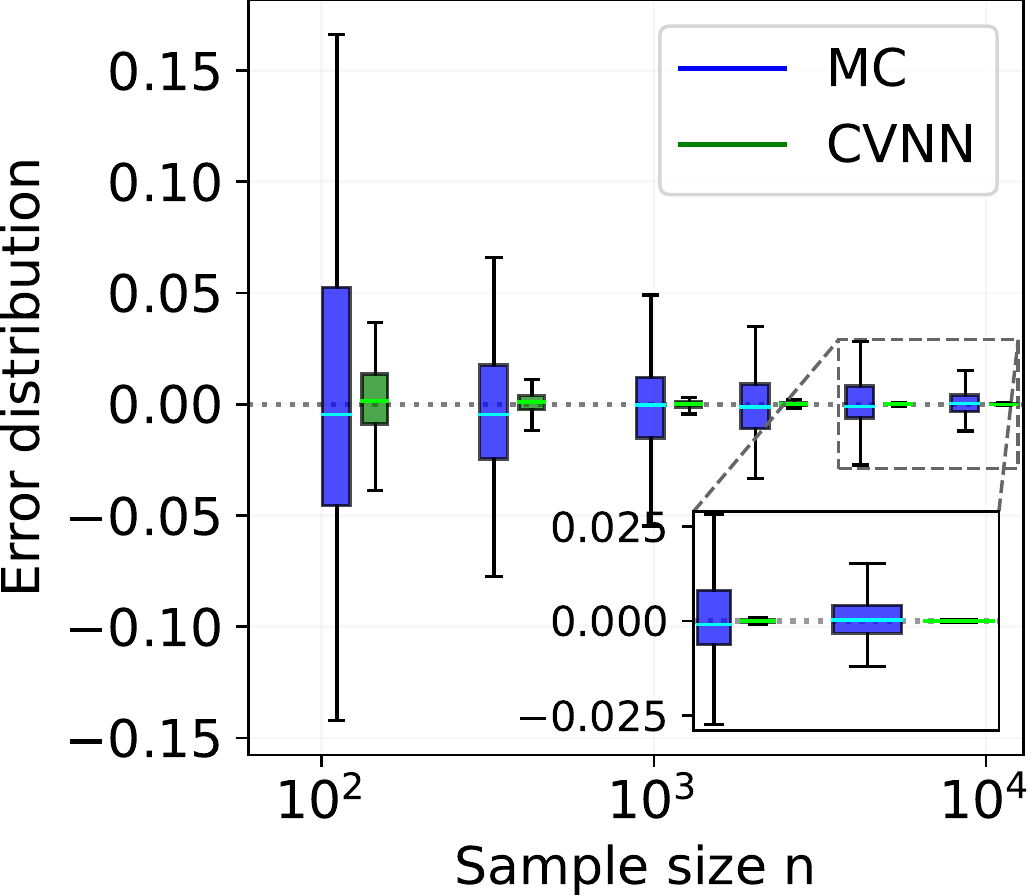}
		\caption{Boxplot $\func_3$}
		\label{fig:func3box}
	\end{subfigure}
	\hfill
	\begin{subfigure}[h]{0.32\linewidth}
		\centering
		\includegraphics[width=\linewidth]{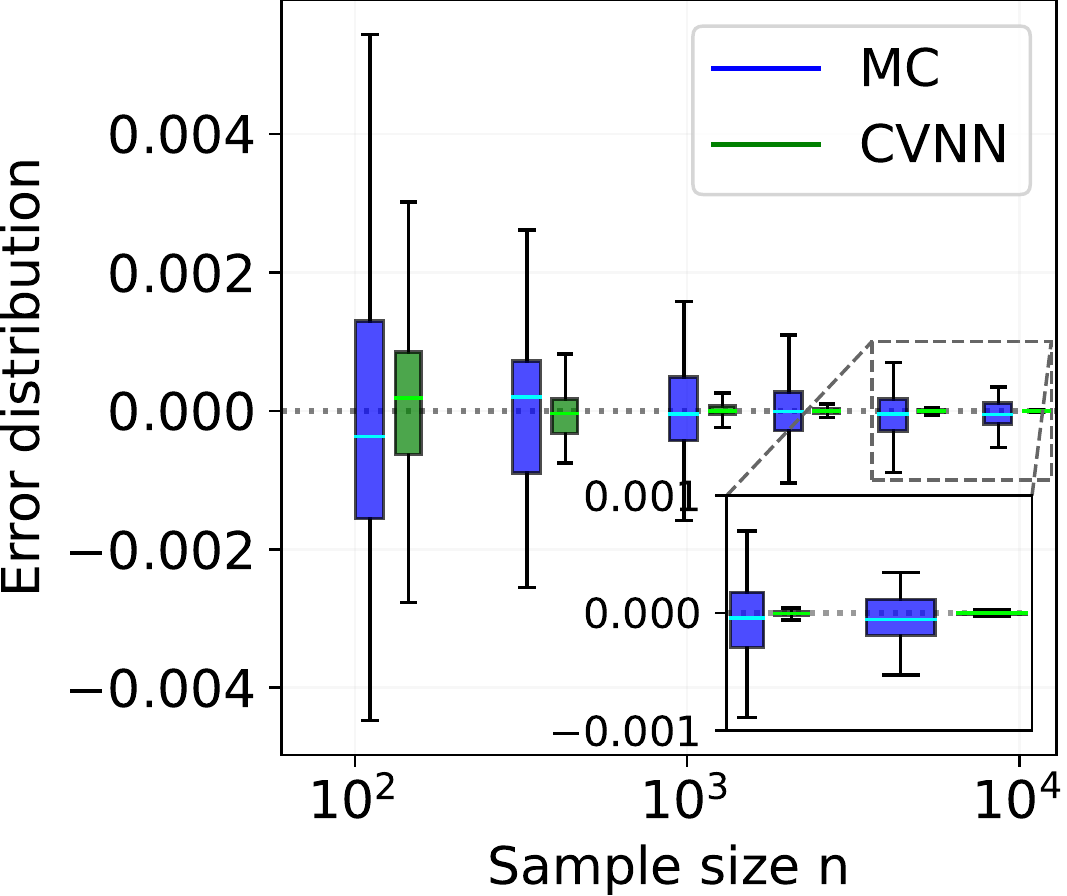}
		\caption{Boxplot $\func_4$}
		\label{fig:func4box}
	\end{subfigure}
	\hfill
	\begin{subfigure}[h]{0.32\linewidth}
		\centering
		\includegraphics[width=\linewidth]{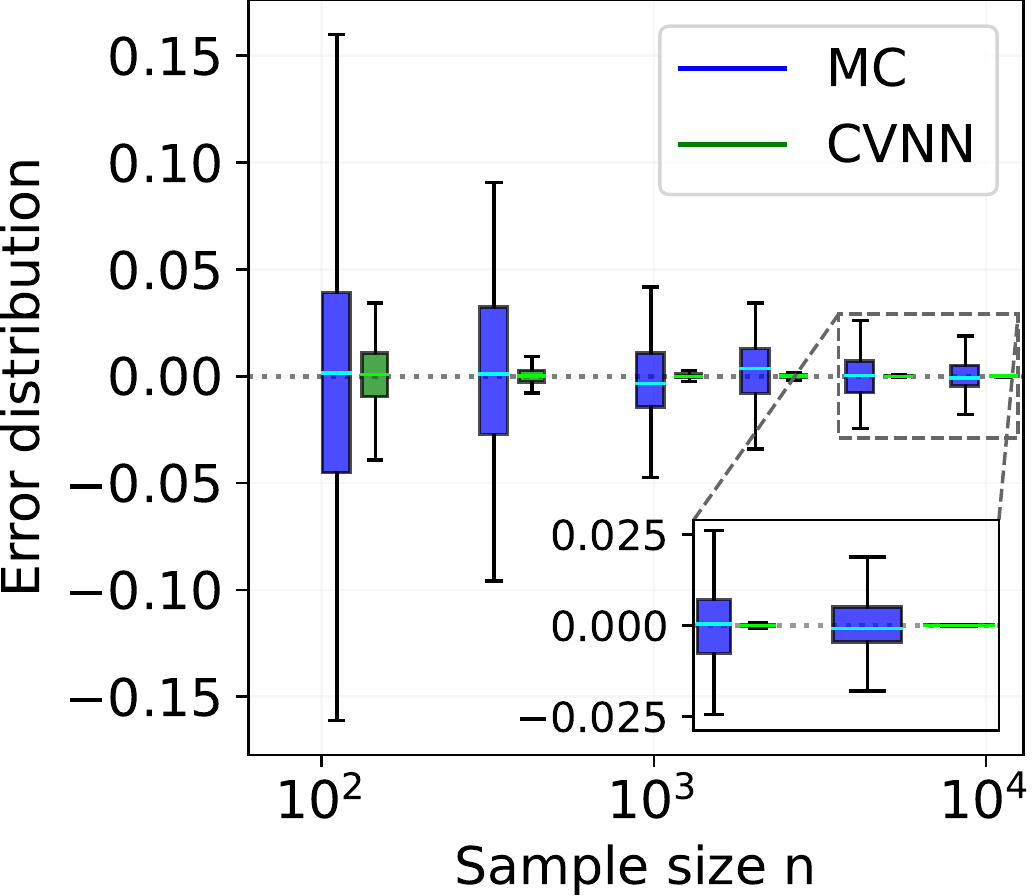}
		\caption{Boxplot $\func_5$}
		\label{fig:func5box}
	\end{subfigure}
	\caption{Root mean squared errors (top) and boxplots of the errors (bottom) obtained over $100$ replications for functions $\func_3$ (left), $\func_4$ (middle) and $\func_5$ (right) in Eq.~\eqref{eq:sphereintegrands} when integrating over $\sphere^2$.}
	\label{fig:sphere} 
\end{figure}

\mypar{Integration on the sphere $\sphere^2$.} 
The unit sphere $\sphere^2$ in $\reals^3$ is a compact Riemannian manifold of dimension $d=2$ with positive curvature. The distances on the sphere are computed using the equations of the great circles. Consider the integral $\int_{\sphere^2} \func \, \diff\Omega$ with $\Omega$ the uniform distribution on $\sphere^2$ and integrands
\begin{equation}
	\label{eq:sphereintegrands}
	\func_3(x,y,z) = \cos(x + y + z), \quad 
	\func_4(x,y,z) = \cos(x)\cos(y)\cos(z), \quad 
	\func_5(x,y,z) = \exp(x-y),
\end{equation}
so that $\meas(\func_3)=\meas(\func_4)= (4\pi/\sqrt{3})\sin(\sqrt{3})$ and $\meas(\func_5)= \pi \sqrt{8} \sinh(\sqrt{2})$. Figure~\ref{fig:sphere} reports the evolution of the root mean squared error and the boxplots of the errors over the sphere $\sphere^2$. The error curves empirically validate the $\Oh(n^{-1})$ convergence rate for the control neighbors estimate while the boxplots highlight the variance reduction of the proposed method.

\subsection{Monte Carlo Integration for Optimal Transport} \label{subsec:OT}

\mypar{Optimal transport.} Optimal transport (OT) is a mathematical framework for measuring the distance between probability distributions. It has proven to be particularly useful in machine learning applications, where it can be used for tasks such domain adaptation \citep{courty2017joint} and image generation \citep{gulrajani2017improved,genevay2018learning}. The recent development of efficient algorithms for solving optimal transport problems, e.g., the Sliced-Wasserstein (SW) distance  \citep{rabin2012wasserstein} and Sinkhorn  distance \citep{cuturi2013sinkhorn}, has made it a practical tool for large-scale machine learning applications.

\mypar{OT distances.} For $\mathsf{X} \subseteq \rset^q$, let $\Prob(\mathsf{X})$ denote the set of probability measures supported on $\mathsf{X}$ and let $p \in [1, \infty)$. The Wasserstein distance of order $p$ between $P, Q \in \Prob(\mathsf{X})$ is 
\begin{align*}
	W_p^p(P,Q) = \inf_{\pi \in \Pi(P,Q)} \int_{\mathsf{X} \times \mathsf{X}} \|x-y\|^p \, \diff \pi(x,y),
\end{align*}
where $\Pi(P,Q) \subset \Prob(\mathsf{X} \times \mathsf{X})$ denotes the set of couplings for $(P,Q)$, i.e., probability measures whose marginals with respect to the first and second variables are $P$ and $Q$ respectively. While the Wasserstein distance enjoys attractive theoretical properties \citep[Chapter 6]{villani2009optimal}, it suffers from a high computational cost. When computing $W_p(P_m,Q_m)$ for discrete distributions $P_m$ and $Q_m$ supported on $m$ points, the worst-case computational complexity scales as $\Oh(m^3 \log m)$ \citep{peyre2019computational}. To overcome this issue, the Sliced-Wasserstein distance takes advantage of the fast computation of the Wasserstein distance between univariate distributions $P,Q \in \Prob(\rset)$. Indeed, for $P_m = (1/m)\sum_{i=1}^m \delta_{x_i}$ and $Q_m = (1/m)\sum_{i=1}^m \delta_{y_i}$ with $x_i, y_i \in \rset$, the $W_p$-distance involves sorting the atoms $x_{(1)} \leq \ldots \leq x_{(m)}$ and $y_{(1)} \leq \ldots \leq y_{(m)}$, yielding
\begin{align*}
	W_p^p(P_m,Q_m) = \frac{1}{m} \sum_{i=1}^m |x_{(i)} - y_{(i)}|^p,
\end{align*} 
leading to a complexity of $\Oh(m \log m)$ operations induced by the sorting step.

\medskip
Recall that $\sphere^{q-1} = \{ \theta \in \rset^q: \|\theta\|=1\}$ is the unit sphere in $\rset^q$ and for $\theta \in \sphere^{q-1}$ let $\theta^\star: \rset^q \to \rset$ denote the linear map $\theta^\star(x) = \langle \theta, x \rangle$ for $x \in \rset^q$. Let $\meas$ denote the uniform distribution on $\sphere^{q-1}$, that is, the normalized volume measure.
The Sliced-Wasserstein (SW) distance \citep{rabin2012wasserstein,bonneel2015sliced,kolouri2019generalized} of order $p$ based on $\meas$ is defined for $P, Q \in \Prob(\mathsf{X})$ as 
\begin{equation}\label{eq:def_sw}
	SW_p^p(P,Q,\rho) 
	= \E_\theta[W_p^p(\theta_{\#}^\star P,\theta_{\#}^\star Q)]
	= \int_{\sphere^{q-1}} W_p^p(\theta_{\#}^\star P,\theta_{\#}^\star Q) \, \diff \meas(\theta),
\end{equation}
where $f_{\#}\xi$ is the \textit{push-forward} measure of $\xi \in \Prob(\rset^q)$ by a measurable function $f$ on $\reals^q$. For $\theta \in \sphere^{q-1}$, the measures $\theta_{\#}^\star P$ and $\theta_{\#}^\star Q$ are the distributions of the projections $\langle \theta, X \rangle$ and $\langle \theta, Y \rangle$ of random vectors $X$ and $Y$ with distributions $P$ and $Q$, respectively. The integrand is Lipschitz in view of Theorem~2.4 in \citet{han2023sliced}.

In practice, the random directions $\theta_1,\ldots,\theta_n$ are sampled independently from the uniform distribution on the unit sphere $\sphere^{q-1}$ and the SW-distance of Eq.~\eqref{eq:def_sw} is approximated using a standard Monte Carlo estimate with $n$ random projections as
\begin{align*}
	\widehat{SW}_p^p(P_m,Q_m) 
	= \frac{1}{n} \sum_{i=1}^n W_p^p \left(
		(\theta_i)_{\#}^\star P_m, 
		(\theta_i)_{\#}^\star Q_m 
	\right).
\end{align*}
The use of more accurate numerical integration methods to compute the SW-distance is an active area of current research; see, e.g., \cite{nguyen2024sliced} and \cite{leluc2024sliced}. Below, we illustrate the benefits of using control neighbors, while at the same time we anticipate that more specialized methods may do even better for this particular integration problem.

\begin{figure}
	\centering
	\begin{subfigure}[h]{0.38\linewidth}
		\centering
		\includegraphics[width=\linewidth]{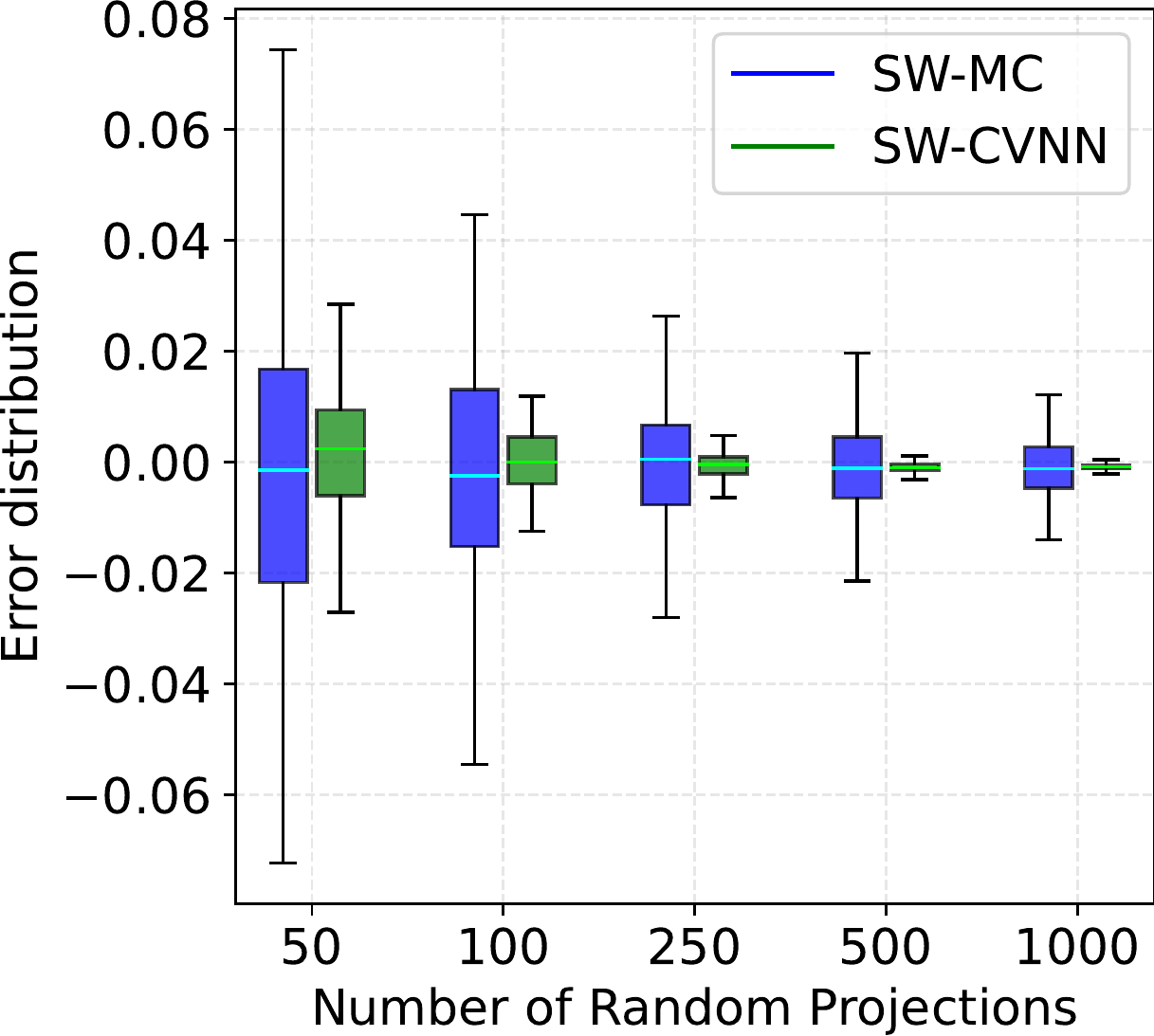}
		\caption{$SW_2, \ q=3$}
		\label{fig:ot_d3}
	\end{subfigure}
	\begin{subfigure}[h]{0.38\linewidth}
		\centering
		\includegraphics[width=\linewidth]{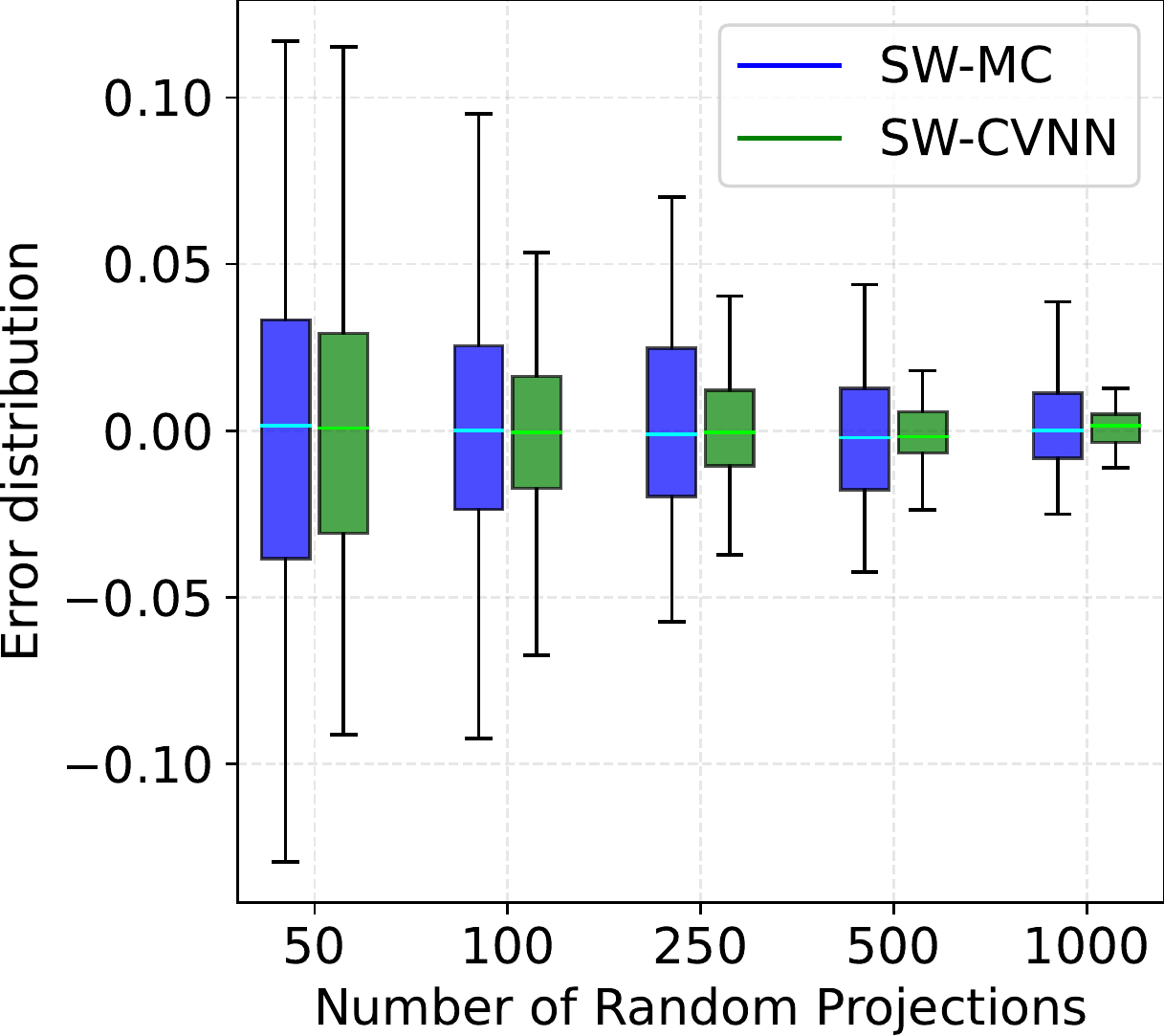}
		\caption{$SW_2, \ q=6$}
		\label{fig:ot_d6}
	\end{subfigure}
	\caption{Boxplots of Sliced-Wasserstein estimates SW-MC and SW-CVNN for Gaussian distributions on $\rset^q$ with $q \in \{3;6\}$. The boxplots are obtained over $100$ replications.}
	\label{fig:OT_box} 
\end{figure}

\mypar{Multivariate Gaussians.} The goal is to compare the variance of the standard Monte Carlo estimate (SW-MC) with the proposed control neighbors estimate (SW-CVNN) when computing the Sliced-Wasserstein distance between two Gaussian distributions. More precisely, we want to compute $SW_2(P,Q)$, where $P = \mathcal{N}_q(m_X,\sigma_X^2 \mathrm{I}_q)$ and $Q = \mathcal{N}_q(m_Y,\sigma_Y^2 \mathrm{I}_q)$ with $m_X, m_Y \sim \mathcal{N}(0,\mathrm{I}_q)$ and $\sigma_X = 2$ and $\sigma_Y=5$. We consider the corresponding empirical distributions $P_m$ and $Q_m$ based on $m=2\,000$ samples and compute the Monte Carlo estimates of the $SW_2$ distance using a number of projections $n \in \{50;100;250;500;1000\}$ in dimension $q \in \{3;6\}$. In the notation of the general theory, we have $d=q-1 \in \{2;5\}$ since we are integrating over the $(q-1)$-dimensional sphere in $\rset^q$. Figure~\ref{fig:OT_box} shows the error distribution of the different Monte Carlo estimates (SW-MC and SW-CVNN) of the $SW_2$ distance where the exact value \citep[see][Appendix S3.1]{nadjahi2021fast} is given by 
\begin{align*}
	SW_2^2 \left( 
		\mathcal{N}_q \big(
			m_X,\sigma_X^2 \mathrm{I}_q
		\big),
		\mathcal{N}_q \big(
			m_Y,\sigma_Y^2 \mathrm{I}_q
		\big) 
	\right) 
	= q^{-1} \|m_X - m_Y \|_2^2 + \rbr{\sigma_X - \sigma_Y}^2.
\end{align*}
The different boxplots highlight the good performance of the control neighbors estimate in terms of variance reduction.

\section{Discussion} \label{sec:discussion}

We have explored the use of nearest neighbors in the construction of control variates for variance reduction in Monte Carlo integration. We have shown that for Hölder integrands of regularity $s \in (0, 1]$ on bounded metric spaces of dimension $d$ as measured by a sufficiently regular probability distribution, a faster rate of convergence, namely $\Oh(n^{-1/2} n^{-s/d})$ as $n \to \infty$, is possible through the construction of a control variate via leave-one-out neighbors. Theoretical guarantees are given both in terms of bounds on the root mean squared error and as concentration inequalities (requiring an additional logarithmic factor). In numerical experiments, the method enjoyed a notable error reduction with respect to Monte Carlo integration.

A drawback is that our method is not able to leverage higher-order smoothness of the integrand besides the Lipschitz property: this comes from the fact that the control variate is piecewise constant on the Voronoi cells induced by the Monte Carlo sample. As a consequence, the accuracy gain is limited on high-dimensional domains. Further, integration of the control function requires an additional Monte Carlo sample, of larger size than the original one, although without needing additional evaluations of the integrand.

\section*{Acknowledgments}
Aigerim Zhuman gratefully acknowledges a research grant from the \textit{National Bank of Belgium} and of the Research Council of the UCLouvain.

\bibliographystyle{plainnat}
\bibliography{references}

\newpage
\newpage

\begin{center}
{\large {\bf\textsc{Appendix: \\ Speeding up Monte Carlo Integration: \\ Control Neighbors for Optimal Convergence}}}
\end{center}

Appendix \ref{app:auxiliary} contains some auxiliary results. Appendix \ref{app:lemmas} gathers the proofs of all the lemmas, while Appendix \ref{app:propositions} is concerned with the proofs of the different propositions. Appendix \ref{app:theorems} comprises the proofs of the theorems and Appendix \ref{app:option_pricing} presents additional numerical experiments for option pricing.

\appendix
\startcontents[sections]
\printcontents[sections]{l}{1}{\setcounter{tocdepth}{2}}

\section{Auxiliary results} \label{app:auxiliary}

The following theorem is a corollary to Theorem~2.1 in \cite{extensionofMcDiarmid2015}.

\begin{theorem}[Extension of McDiarmid’s inequality] \label{th:mcdiarmid_extension}
	Let $X_1,\dots,X_n$ independent random variables with $X_\ell$ taking values in some measurable space $\Omega_\ell$. Define $X=(X_1,\dots,X_n) \in  \prod_{\ell=1}^{n} \Omega_\ell =: \Omega$. Let $\phi: \Omega \to \reals$ be a measurable function and assume there exists an event $A \subseteq \Omega$ and constants $c_1,\ldots,c_n \ge 0$ such that, for any $\ell \in \{1,\ldots,n\}$ and for any $x \in A$ and $x' \in A$  that differ only in the $\ell$-th coordinate, we have
	\begin{equation} \label{eq:bounded_diff_propeprty}
		\abs{\phi(x)-\phi(x')} \le c_\ell.
	\end{equation} 
	Then, for any $t \ge 0$, we have, writing $m = \E[\phi \mid A]$, $p = 1 - \prob(A)$ and $\bar{c} = \sum_{\ell=1}^n c_\ell$,
	\begin{equation*}
		\prob \left(|\phi(X) - m|\ge t \right) 
		\le p + 2 \exp \left(- \frac{2 \max(0,t-p \bar{c})^2}{\sum_{\ell=1}^n c_\ell^2} \right).
	\end{equation*}
\end{theorem}

\begin{lemma}[\cite{kabatyanskii1974bounds, zeger1994number}] \label{lem:degree_nn}
	Given a point set $P \subset \Rd$ and $x \in P$, the maximum number of points in $P$ that can have $x$ as nearest neighbor is bounded by the $d$-dimensional kissing number $\psi_d$, where $\psi_d \le 2^{0.401 d (1+ \oh(1))}$ as $d \to \infty$.
\end{lemma}

\begin{proposition}\label{prop:difference_loonn_nn}
	Under Assumptions~(A1), (A2) and (A3) in the paper, if $n \ge 4$, we have
	\[
		\Exp \sbr{ \abs{ \frac{1}{n}  \sum_{i=1}^n \meas ( \hgni )   - \meas(\hat \func_n ) } }
		\le\rbr{s/d +2} L   C_0^{-s/d} \, \Gamma( s/d +1) \, n^{-1-s/d} 
	\]
	and
	\[
		\Exp \sbr{ \abs{ \frac{1}{n}  \sum_{i=1}^n \meas ( \hgni )   - \meas(\hat \func_n ) }^2 } 
		\le 4 \rbr{s/d +1}  L^2  C_0^{-2s/d} \, \Gamma ( 2s/d +1) \, n^{-1-2s/d}.
	\]
\end{proposition}

\begin{lemma}[\cite{portier2021nonasymptoticbound}] \label{lm:nonasymptoticbound}
	Under Assumptions~(A1') and~(A2) in the paper, for all $n \ge 1$, all $\delta \in (0,1)$ and all $k \in \{1,\ldots,n\}$ such that 
$8 d \log(12 n/\delta) \le k \le r_0^d n b c V_d / 2$, it holds, with probability at least $1-\delta$, that 
	\begin{equation*}
		\sup_{x \in \cX} \hat{\tau}_{n,k}(x) 
		\le \left(\frac{2 k }{n b c V_d}\right)^{1/d} =: \bar{\tau}_{n,k}.
	\end{equation*}
\end{lemma}

\begin{lemma}[\cite{extensionofMcDiarmid2015}]  	
	\label{rk:expect_diff}
	Let $\phi: \Omega \to \reals$ be a bounded, measurable function and let $A \subseteq \Omega$ be an event with $\prob(A) > 0$. Writing $m = \E[\phi \mid A]$, $p = 1 - \prob(A)$ and $\ninf{\phi} = \sup_{x \in \Omega} |\phi(x)|$, we have
	\begin{equation*}
		\abs{\E[\phi] - m} \le 2 p \ninf{\phi}.
	\end{equation*}
\end{lemma}

\begin{lemma}
	\label{lem:Gammaineq}
	For $x, a > 0$, we have $\Gamma(x + 1 + a) > x^a \Gamma(x+1)$.
\end{lemma}

\section{Proofs of Lemmas} \label{app:lemmas}

\subsection{Proof of Lemma  \ref{prop:modification_integral}}

Given any collection $X_1,\ldots, X_n$ of $n$ distinct points, if $j\neq i$, then $ \hgni$ and $ \hat \func_{n}$ are the same on $ S_{n,j} $. It holds that
\begin{align*}
 \hgni (x) - \hat \func_{n}   (x) 
  =  \bigl\{ \hgni(x) - \hat{\func}_{n}(x) \bigr\}  \ind _{ S_{n,i} } (x) .
\end{align*}
Now using that $\bar \func_n $ and $\hgni$ are the same on $ S_{n,i} $, it follows that
\[
 \hgni (x) - \hat \func_{n}(x)
  =  \bigl\{ \bar{\func}_{n}(x) - \hat{\func}_{n}(x) \bigr\} \ind_{S_{n,i}}(x) .
\]
Taking the sum and using $\sum_{i=1}^n \ind _{S_{n,i}}(x) = 1$ gives 
\[
  \sum_{i=1}^n  \bigl\{ \hgni (x)  -  \hat \func_{n}(x) \bigr\} 
  =  \bar \func_n  (x) -   \hat \func_{n}(x),
\]
and the result follows by integrating with respect to $\meas$.

\subsection{Proof of Lemma \ref{prop:di_ci}}

Because the  Voronoi cells define a partition of $\ms$, we have for any $x \in \ms$,
	\[
		\hgni(x) = \sum_{j:j \ne i} \func(X_j) \ind_{S_{n,j}^{(i)}}(x)
	\]
and in particular
	\[
		\hgni(X_i) = \sum_{j:j \ne i} \func(X_j) \ind_{S_{n,j}^{(i)}}(X_i),
	\]
from which we deduce
	\[
		\sum_{i=1}^n \hgni(X_i)
		= \sum_{j=1}^n \func(X_j) \sum_{i:i \ne j} \ind_{S_{n,j}^{(i)}}(X_i) 
		= \sum_{j=1}^n \func(X_j) \, \hat{d}_{n,j}.
	\]
	Further, we have
	\[
		\meas(\hat{\func}_n) 
		= \sum_{i=1}^n \meas\left( \func(X_i) \ind_{S_{n,i}} \right) 
		= \sum_{i=1}^n \func(X_i) \, \meas(S_{n,i}) 
		= \sum_{i=1}^n \func(X_i) \, {V}_{n,i}
	\]
	and
	\[
		\sum_{i=1}^n \meas\bigl(\hgni\bigr) 
		= \sum_{i=1}^n \sum_{j:j \ne i} \meas\left( \func(X_j) \ind_{S_{n,j}^{(i)}} \right) 
		= \sum_{j=1}^n \func(X_j) \sum_{i:i \ne j} V_{n,j}^{(i)} 
		= \sum_{j=1}^n \func(X_j) \, \hat{c}_{n,j}.
	\]

\subsection{Proof of Lemma~\ref{lem:1NN-distance}}

For $x \in \ms$, $r \in (0, \diam(\ms)]$ and integer $k \in \{1,\ldots,n\}$, we have
\begin{align*}
	\prob \sbr{\hat{\tau}_{n,k}(x) > r}
	&= \prob \sbr{\sum_{i=1}^n \ind_{\{\rho(X_i,x) \le r\}} \le k-1} \\
	&= \prob \sbr{\Bin(n, \mu(B(x,r)) \le k-1} \\
	&\le \prob \sbr{\Bin(n, C_0r^d) \le k-1} \\
	&= \sum_{j=0}^{k-1} \binom{n}{j} (C_0 r^d)^j (1 - C_0 r^d)^{n-j}.
\end{align*}
The equality remains valid for $r > 0$ such that $C_0 r^d \le 1$, since if $r > \diam(M)$, the probability is zero. We find
\begin{align*}
	\expec \sbr{\hat{\tau}_{n,k}(x)^q}
	&= \int_0^\infty \prob \sbr{\hat{\tau}_{n,k}(x) > t^{1/q} } \, \diff t \\
	&\le \sum_{j=0}^{k-1} \binom{n}{j} \int_0^{C_0^{-q/d} }\rbr{C_0 t^{d/q}}^j \rbr{1 - C_0 t^{d/q}}^{n-j} \, \diff t \\
	&= \sum_{j=0}^{k-1} \binom{n}{j} \int_0^1 u^j (1-u)^{n-j} \cdot \frac{q}{d} C_0^{-q/d} u^{q/d-1} \, \diff u \\
	&= \frac{q}{d} C_0^{-q/d} \sum_{j=0}^{k-1} \binom{n}{j} \frac{\Gamma(j+q/d) \Gamma(n-j+1)}{\Gamma(n+q/d+1)},
\end{align*}
since the integral can be expressed in terms of the Beta function, which can in turn be developed in terms of the Gamma function. Writing the binomial coefficient in terms of the Gamma function too and simplifying yields
\[
	\expec \sbr{\hat{\tau}_{n,k}(x)^q}
	\le \frac{q}{d} C_0^{-q/d} \frac{\Gamma(n+1)}{\Gamma(n+q/d+1)} \sum_{j=0}^{k-1} 
	\frac{\Gamma(j+q/d)}{\Gamma(j+1)}
	\le \frac{q}{d} \rbr{C_0 n}^{-q/d} \sum_{j=0}^{k-1} \frac{\Gamma(j+q/d)}{\Gamma(j+1)},
\]
where we used Lemma~\ref{lem:Gammaineq} in the second inequality. 
By induction, one can obtain 
\begin{equation} \label{eq:gammaratioind}
\frac{q}{d} \sum_{j=0}^{k-1} \frac{\Gamma(j+q/d)}{\Gamma(j+1)} = \frac{\Gamma(k+q/d)}{\Gamma(k)}.
\end{equation}
For $k=1$, \eqref{eq:gammaratioind} reduces to 
\[
\frac{(q/d) \Gamma(q/d)}{\Gamma(1)} = \frac{\Gamma(q/d+1)}{\Gamma(1)}.
\]
Assume that \eqref{eq:gammaratioind} holds for $k-1$, i.e.,
\begin{equation} \label{eq:inductionstep}
	\frac{q}{d} \sum_{j=0}^{k-2} \frac{\Gamma(j+q/d)}{\Gamma(j+1)} = \frac{\Gamma(k-1+q/d)}{\Gamma(k-1)}.
\end{equation}
Then, by \eqref{eq:inductionstep} and $\Gamma(k)=(k-1) \Gamma(k-1)$, we have
\begin{align*}
\frac{q}{d} \sum_{j=0}^{k-1} \frac{\Gamma(j+q/d)}{\Gamma(j+1)} &= \frac{q}{d} \sum_{j=0}^{k-2} \frac{\Gamma(j+q/d)}{\Gamma(j+1)} + \frac{q}{d} \frac{\Gamma(k-1+q/d)}{\Gamma(k)}\\
&=\frac{\Gamma(k-1+q/d)}{\Gamma(k-1)}+\frac{q}{d}\frac{\Gamma(k-1+q/d)}{\Gamma(k)}\\
&=\frac{(k-1+q/d) \Gamma(k-1+q/d)}{\Gamma(k)} \\
&= \frac{\Gamma(k+q/d)}{\Gamma(k)}.
\end{align*}
Therefore,
\[
\expec \sbr{\hat{\tau}_{n,k}(x)^q}
\le \rbr{C_0 n}^{-q/d} \frac{\Gamma(k+q/d)}{\Gamma(k)}.
\]

\subsection{Proof of Lemma~\ref{lem:Gammaineq}}
	The Gamma function is strictly logarithmically convex. Writing
	\[
	x+1 = (1-\lambda) x + \lambda(x+1+a)
	\text{ with } \lambda = \frac{1}{1+a}
	\]
	yields
	\[
	\Gamma(x+1) 
	< \Gamma(x)^{1-\lambda} \, \Gamma(x+1+a)^\lambda.
	\]
	Since $\Gamma(x) = x^{-1}\Gamma(x+1)$, it follows that
	\begin{align*}
		\Gamma(x+1+a)
		> \cbr{\Gamma(x+1) \Gamma(x)^{\lambda-1}}^{1/\lambda} 
		= \cbr{\Gamma(x+1)^{1+\lambda-1} x^{1-\lambda}}^{1/\lambda} 
		= \Gamma(x+1) x^{1/\lambda-1} = \Gamma(x+1) x^a,
	\end{align*}
	as required.

\section{Proofs of Propositions} \label{app:propositions}

\subsection{Proof of Proposition \ref{prop:quad}}
Using Lemma \ref{prop:di_ci}, we find 
	\begin{align*}
		\munNN(\func)
		&= \frac{1}{n} \sum_{i=1}^n [\func(X_i) - \{\hgni(X_i) - \meas(\hat{\func}_n)\}] \\
		&= \frac{1}{n} \sum_{i=1}^n \func(X_i) 
		- \frac{1}{n} \sum_{j=1}^n \func(X_j) \, \hat{d}_{n,j} 
		+ \sum_{i=1}^n \func(X_i) \, {V}_{n,i} \\
		&= \frac{1}{n} \sum_{i=1}^n \left( 1 - \hat{d}_{n,i} + n {V}_{n,i} \right) \func(X_i)
	\end{align*}
	and, similarly,
	\begin{align*}
		\munNNloo
		&= \frac{1}{n} \sum_{i=1}^n [\func(X_i) - \{\hgni(X_i) - \meas(\hgni)\}] \\
		&= \frac{1}{n} \sum_{i=1}^n \func(X_i) - \frac{1}{n} \sum_{j=1}^n \func(X_j) \hat{d}_{n,j} + \frac{1}{n} \sum_{j=1}^n \func(X_j) \hat{c}_{n,j} \\
		&= \frac{1}{n} \sum_{i=1}^n \left( 1 - \hat{d}_{n,i} + \hat{c}_{n,i} \right) \func(X_i),
	\end{align*}
	as required.


\subsection{Proof of Proposition \ref{lem:Lpconv}}

	Let $N_{n,k}(x)$ denote the set of indices $i \in \{1,\ldots,n\}$ such that $X_i$ is among the $k$ nearest neighbors of $x \in \ms$. 
	Using the Hölder property and the definition of $\hat \tau_{n,k}$ we have, for each $i\in N_{n,k} (x)$,
	\[
	\abs{ \func(X_i) - \func(x) }
	\leq L \rho(X_i, x)^s
	\leq L \hat \tau_{n,k} (x)^s.
	\]
	It follows that
	\begin{align*}
		\abs{ \hat \func_{n,k} (x) - \func(x) }
		\leq  \frac{1}{k} \sum_{i\in N_{n,k} (x) } 
		\abs{ \func(X_i) - \func(x) }
		\leq L \hat \tau_{n,k} (x)^s.
	\end{align*}
	Next we apply  Lemma~\ref{lem:1NN-distance} with $q = sp$ to obtain
	\begin{align*}
		\expec [ |\hat \func_{n,k} (x) - \func(x) |^p ] 
		&\leq L^p \expec [ \hat \tau_{n,k} (x)^{sp} ] \\
		&\leq L^p (C_0 n)^{-sp/d} \frac{\Gamma(k + sp/d)}{\Gamma(k)}.
	\end{align*}
	Raising both sides to the power $1/p$ yields the stated inequality.

\subsection{Proof of Proposition  \ref{prop:difference_loonn_nn}}

Using the fact that $\hgni$ and $ \hat \func_n$ coincide outside $ S_{n,i} $  and that $\hat \func_n (x) = \func(X_i)$ for $x \in S_{n,i}$, we have
	\begin{align*}
		\meas(\hgni) - \meas(\hat \func_n) 
		& = \meas(\hgni - \hat \func_n) \\ 
		& = \meas \bigl((\hgni - \hat \func_n) \ind _{S_{n,i}} \bigr) \\
		& = \meas \bigl(\{\hgni -  \func(X_i)\} \ind _{S_{n,i}}\bigr).
	\end{align*}
		Write $R_n = \left| \frac{1}{n} \sum_{i=1}^{n} \bigl\{\meas(\hgni) - \meas(\hat \func_n) \bigr\} \right|$. Then 
	\begin{align*}
		R_n 
		& \le \frac{1}{n} \sum_{i=1}^{n} \abs{ \meas(\hgni) - \meas(\hat \func_n) } 	\\
		& =  \frac{1}{n} \sum_{i=1}^{n} \abs{ \meas\bigl(\{\hgni -  \func(X_i)\} \ind _{S_{n,i} }\bigr) }\\
		&\le \frac{1}{n} \sum_{i=1}^{n}  \int_{S_{n,i}} \abs{\hgni(x) - \func(X_i)} \, \diff \meas(x) \\
		& = \frac{1}{n} \sum_{i=1}^{n}  \int_{S_{n,i}} \abs{\func\bigl(\hat N_n^{(i)}(x)\bigr) - \func(X_i)} \, \diff \meas(x) .
\end{align*}		
Further, since $\func$ is Hölder continuous with exponent $s \in (0, 1]$, we have, by the triangle inequality and the sub-additivity of the function $x \mapsto x^s$,
\begin{align*}
	R_n 
	& \le \frac{L}{n} \sum_{i=1}^{n}  \int_{S_{n,i}} \rho\bigl(\hat N_n^{(i)}(x), X_i \bigr)^s \, \diff \meas(x) \\
		& \le \frac{L}{n} \sum_{i=1}^{n}  \int_{S_{n,i}} \cbr{ \rho\bigl(\hat N_n^{(i)}(x) , x\bigr) + \rho( x , X_i ) }^s \diff \meas(x) \\
		& \le \frac{L}{n} \sum_{i=1}^{n}  \int_{S_{n,i}} \cbr{ \rho\bigl(\hat N_n^{(i)}(x) , x \bigr)^s + \rho( x , X_i )^s } \, \diff \meas(x) \\
		& = \frac{L}{n} \sum_{i=1}^{n}  \int_{S_{n,i}} \cbr{ \hat \tau_n^{(i)}(x)^s + \hat \tau_n(x)^s } \, \diff \meas(x).
	\end{align*}
	Concerning $R_n^2$, one can use Jensen's inequality to obtain
	\begin{align*}
		R_n^2 & =  \abs{\frac{1}{n} \sum_{i=1}^{n} \cbr{\meas(\hgni) - \meas(\hat \func_n)}}^2 \\
		& \le  \frac{1}{n} \sum_{i=1}^{n} 
		\abs{\meas(\hgni) - \meas(\hat \func_n)}^2 \\
		& \le \frac{L^2}{n} \sum_{i=1}^{n} \int_{S_{n,i}} \cbr{ \rho\bigl(\hat N_n^{(i)}(x) , x \bigr)^s + \rho( x , X_i )^s }^2 \, \diff \meas(x) \\
		& \le \frac{2 L^2}{n} \sum_{i=1}^{n} \int_{S_{n,i}} \cbr{ \rho\bigl(\hat N_n^{(i)}(x) , x \bigr)^{2s} + \rho( x , X_i )^{2s} } \, \diff \meas(x). 
	\end{align*}
	For $x \in S_{n,i}$, the nearest neighbor in $\{X_1,\dots,X_n\}$ is $X_i$, and thus
	\begin{equation*}
		\hat \tau_n^{(i)}(x) = \hat \tau_{n,2}(x)
	\end{equation*}
	is the distance to the second nearest neighbor in $\{X_1,\dots,X_n\}$.
	We get
	\begin{align*}
		R_n  
		& \le \frac{L}{n} \sum_{i=1}^{n} \int_{S_{n,i}} \cbr{ \hat \tau_{n,2} (x) ^s + \hat \tau_n (x)^s } \, \diff \meas(x) \\
		& =  \frac{L}{n} \int_{\ms} \cbr{ \hat \tau_{n,2} (x)^s + \hat \tau_n(x)^s } \, \diff \meas(x)
	\end{align*}
	and, similarly,
	\begin{align*}
		R_n^2  
		& \le \frac{2 L^2}{n} \sum_{i=1}^{n} \int_{S_{n,i}} \cbr{ \hat \tau_{n,2} (x) ^{2s} + \hat \tau_n (x) ^{2s} } \, \diff \meas(x) \\
		& =  \frac{2 L^2}{n} \int_{\ms} \cbr{ \hat \tau_{n,2} (x)^{2s} + \hat \tau_n(x)^{2s} } \, \diff \meas(x). 
	\end{align*}
	Consequently, by Lemma \ref{lem:1NN-distance},
	\begin{align*}
		\Exp [R_n] 
		& \le \frac{L}{n} \int_{\ms}  \Exp \sbr{\hat \tau_{n,2} (x)^s + \hat \tau_n(x)^s} \diff \meas(x) \\
		& \le \frac{L}{n} \rbr{ \sup_{x\in \ms} \Exp \sbr{\hat \tau_{n,2} (x)^s} + \sup_{x\in \ms} \Exp \sbr{\hat \tau_n(x)^s } } \\
		& \le \frac{L}{n} \cbr{   \rbr{C_0 n}^{-s/d} \Gamma ( s/d +2) + \rbr{C_0 n}^{-s/d} \Gamma ( s/d +1)  } \\
		& = \rbr{s/d +2} L   C_0^{-s/d} \, \Gamma( s/d +1) \, n^{-1-s/d} 
	\end{align*}
	and
	\begin{align*}
		\Exp [R_n^2] 
		& \le \frac{2 L^2}{n} \int_{\ms}  \Exp \sbr{\hat \tau_{n,2} (x)^{2s} + \hat \tau_n(x)^{2s} } \diff \meas(x) \\
		& \le \frac{2 L^2}{n} \rbr{ \sup_{x\in \ms} \Exp \sbr{\hat \tau_{n,2} (x)^{2s}} + \sup_{x\in \ms} \Exp \sbr{ \hat \tau_n(x) ^{2s} } } \\
		& \le \frac{2 L^2}{n} \cbr{ \rbr{C_0 n}^{-2 s/d} \Gamma ( 2 s/d +2) + \rbr{C_0 n}^{-2 s/d} \Gamma ( 2 s/d +1) } \\
		& = 4 \rbr{s/d +1}  L^2  C_0^{-2s/d} \, \Gamma ( 2s/d +1) \, n^{-1-2s/d}. 
	\end{align*}

\section{Proofs of Theorems} \label{app:theorems}

\subsection{Proof of Theorem \ref{th:cv_rate} for $\munNNloo$}

Let $Y_{n,i} = \hgni (X_i) -   \meas ( \hgni )  $ and write 
\begin{align*}
\munNNloo   (\func)  - \meas   (\func)   &= \frac{1}{n} \sum_{i=1}^n \rbr{Y_i - Y_ {n,i} } 
\end{align*}
with $Y_i = \func(X_i) - \meas(\func) $. Then write
\begin{align}
\nonumber
	n^2 \expec \sbr{\cbr{ \munNNloo   (\func)  - \meas   (\func) }^2} 
	&=  \sum_{i=1 }^n   \expec \sbr{ (Y_i - Y_ {n,i})^2  }   +   \sum_{i\neq j }  \expec \sbr{ \rbr{Y_i - Y_ {n,i} }\rbr{Y_j - Y_ {n,j} } }  \\
	&= n  \expec \sbr{ (Y_1 - Y_ {n,1})^2  } +  n(n-1) \expec \sbr{ \rbr{Y_1 - Y_ {n,1} }\rbr{Y_2 - Y_ {n,2} } }.
\label{eq:munNNlooerrdecomp}
\end{align}
We decompose $Y_{n,1} $ into two terms, one of which does not depend on $X_2$. We also use the fact that the Voronoi partition made with $(n-1)$ points is more detailed than the one constructed with $(n-2)$ points, i.e., $S_{n-1,i}^{(1)} \subset  S_{n-2,i}^{(1,2)}$ for $i =3, \ldots, n$. Define the map $\mathcal N ^{(1,2)} :\ms \to \ms $ such that  $\mathcal N ^{(1,2)} (x) $ is the nearest neighbor to $ x$ among the sample $\{X_3,X_4,\ldots, X_n\}$. Using the identity $\mathcal N ^{(1,2)} ( x ) = X_i$ whenever $x\in  S_{n-1,i}^{(1)}$ for $i \ge 3$, we write
\begin{align*}
\hat \func_n^{(1,2)} (x) &=\func\bigl( \mathcal N ^{(1,2)} ( x )  \bigr) \\
&=  \func\bigl( \mathcal N ^{(1,2)} ( x )  \bigr)  \left(  \sum_{i=2}^{n} \ind_{ S_{n-1,i}^{(1)}  }(x)  \right)  \\
&=  \func\bigl( \mathcal N ^{(1,2)} ( x )  \bigr) \ind_{ S_{n-1,2}^{(1)}  }(x)    + \sum_{i=3 } ^n \func(X_i) \ind_{ S_{n-1,i}^{(1)}  }(x)   \\
& = \cbr{\func\bigl( \mathcal N ^{(1,2)} ( x ) \bigr) -  \func(X_2)  } \ind_{ S_{n-1,2}^{(1)}  }(x)    + \sum_{i=2} ^n \func(X_i) \ind_{ S_{n-1,i}^{(1)}  }(x)  .
\end{align*}
It follows that 
\[
	\hat \func_n^{(1)} (x) 
	=  \hat L^{(1)}  (x)  + \hat \func_n^{(1,2)} (x)
\]
with $ \hat L^{(1)}  (x) = \cbr{ \func(X_2)  - \func\bigl( \mathcal N ^{(1,2)} ( x )  \bigr)   } \ind_{ S_{n-1,2}^{(1)}  }(x) $. Therefore, since $Y_{n,1} = \hat \func_n^{(1)}(X_1) -   \meas ( \hat \func_n^{(1)} )$, we get
\[
	Y_1 - Y_ {n,1} = Y_1 - 
	\cbr{\hat L^{(1)}  (X_1)  - \meas(\hat L^{(1)})} -
	\cbr{\hat\func_n^{(1,2)}(X_1) - \meas(\hat\func_n^{(1,2)})}.
\]
Write
\begin{align*}
	A_1 &= Y_1 = \func(X_1) - \meas(\func), \\
	A_2 &= Y_2 = \func(X_2) - \meas(\func), \\
	B_1 &=  \hat L^{(1)}  (X_1)  - \meas(\hat L^{(1)}), \\
	B_2 &=  \hat L^{(2)}  (X_2)  - \meas(\hat L^{(2)}), \\
	C_1 &=  \hat \func_n^{(1,2)}   (X_1) - \meas(  \hat \func_n^{(1,2)} ), \\
	C_2 &=  \hat \func_n^{(1,2)}   (X_2) - \meas(  \hat \func_n^{(1,2)} ),
\end{align*} 
where $\hat L^{(2)}  (x) = \cbr{\func(X_1)  - \func\bigl(\mathcal N^{(1,2)}(x)\bigr)} \ind_{ S_{n-1,1}^{(2)}  }(x)$. Then
\begin{align*}
	\Exp \sbr{\rbr{Y_1 - Y_ {n,1} }\rbr{Y_2 - Y_ {n,2} }} 
	&= \Exp \sbr{\rbr{A_1 - B_1 - C_1} \rbr{A_2 - B_2 - C_2}} \\
	&= \Exp[A_1 A_2] - \Exp[A_1 B_2] - \Exp[A_1 C_2] \\
	&\quad\null - \Exp[B_1 A_2] + \Exp[B_1 B_2] + \Exp[B_1 C_2] \\
	&\quad\null - \Exp[C_1 A_2] + \Exp[C_1 B_2] + \Exp[C_1 C_2].
\end{align*}
Since $A_1$ and $A_2$ are independent, $\expec [A_1 A_2] = 0$. This also applies to $\Exp  [A_1 C_2]$ and $\Exp  [C_1 A_2]$.
Considering $\Exp[A_1 B_2]$ gives
\begin{align*}
	\Exp[A_1 B_2] 
	& = \Exp \sbr{  Y_1  \cbr{ \hat L^{(2)}  (X_2)  - \meas( \hat  L^{(2)})  } } \\
	& = \Exp \sbr{ \Exp \sbr{  Y_1  \cbr{ \hat L^{(2)}  (X_2)  - \meas( \hat  L^{(2)})  }  \mid X_1, X_3, \dots, X_n }} \\
	& = \Exp \sbr{ Y_1 \Exp \sbr{ \hat L^{(2)}  (X_2)  - \meas( \hat  L^{(2)})  \mid X_1, X_3, \dots, X_n } } = 0.
\end{align*} 
Due to similar reasoning, $\Exp[B_1 A_2] = 0$, $\Exp[B_1 C_2] = 0$ and $\Exp[C_1 B_2] = 0$. 
For $\Exp[C_1 C_2]$, we have
\begin{align*}
	\Exp[C_1 C_2] 
	& = \Exp \sbr{ \cbr{\hat \func_n^{(1,2)}   (X_1) - \meas\bigl(\hat \func_n^{(1,2)}\bigr)} \cbr{\hat \func_n^{(1,2)}   (X_2) - \meas\bigl(\hat \func_n^{(1,2)}\bigr)} } \\
	& = \Exp \sbr{ \Exp \sbr{
		\cbr{\hat \func_n^{(1,2)}   (X_1) - \meas\bigl( \hat \func_n^{(1,2)} \bigr)}
		\cbr{\hat \func_n^{(1,2)}   (X_2) - \meas\bigl( \hat \func_n^{(1,2)} \bigr)}
		\mid X_3, \dots,X_n 
	} }
	= 0.
\end{align*}
Therefore, we get
\[
	\expec \sbr{ \rbr{Y_1 - Y_ {n,1} }\rbr{Y_2 - Y_ {n,2} } } 
	= \expec [ B_1 B_2 ]
	= \expec \sbr{ 
		\cbr{  \hat L^{(1)}  (X_1)  - \meas( \hat L^{(1)}) }
		\cbr{  \hat L^{(2)}  (X_2)  - \meas( \hat L^{(2)}) }
	 } .
\]
Write the four terms on the right-hand side as $A,B,C,D$, respectively. The Cauchy--Schwarz inequality gives
$\| (A- B) (C-D) \|_1  \leq \|A- B\|_2 \, \| C-D\|_2$, while the fact that $B$ and $D$ are conditional expectations of $A$ and $C$, respectively, leads to $\| (A- B) (C-D) \|_1 \leq \|A\|_2 \, \|C\|_2 = \|A\|_2^2 $. As a result,
\[
	\expec \sbr{ \rbr{Y_1 - Y_ {n,1} }\rbr{Y_2 - Y_ {n,2} } } 
	\leq  \expec \sbr{ \hat L^{(1)}(X_1)^2}.
\]
Using the Hölder property, we obtain, by the triangle inequality and the sub-additivity of $x \mapsto x^s$,
\begin{align*}
	\abs{\hat L^{(1)}  (x)}
	&= \abs{\func(X_2)  - \func\bigl( \mathcal N ^{(1,2)} ( x ) \bigr) } \ind_{ S_{n-1,2}^{(1)}  }(x)  \\
	&= \abs{\func\bigl( \mathcal N ^{(1)} ( x ) \bigr) - \func\bigl( \mathcal N ^{(1,2)} ( x ) \bigr) } \ind_{ S_{n-1,2}^{(1)}  }(x) \\
	& \leq  L \, \rho\bigl(   \mathcal N ^{(1)} ( x )  ,  \mathcal N ^{(1,2)} ( x ) \bigr)^s \, \ind_{ S_{n-1,2}^{(1)}  }(x) \\
	&\leq 2 L \, \rho\bigl(   x , \mathcal N ^{(1,2)} ( x ) \bigr)^s \ind_{ S_{n-1,2}^{(1)}  }(x)\\
	&\leq 2 L \, \rho\bigl(   x , \mathcal N ^{(1,2)} ( x ) \bigr)^s \ind_{ B( x,\hat \tau^{(1)}  (x))  }(X_2)\\
	& \leq  2 L \, \rho\bigl(   x ,  \mathcal N ^{(1,2)} ( x ) \bigr)^s \ind_{ B( x,\hat \tau^{(1,2)}  (x))  }(X_2).
\end{align*}
Hence
\[
	\Exp \sbr{ \abs{\hat L^{(1)}(x)}^2 \mid X_3, \dots, X_n }
	\leq 4 L^2 \rho\bigl(x, \mathcal N^{(1,2)}(x) \bigr)^{2s} \, 
	\meas\rbr{B\bigl( x,\hat \tau^{(1,2)}(x) \bigr)} .
\]
Moreover,
\[
	\meas\rbr{B\bigl( x,\hat \tau^{(1,2)}  (x) \bigr)}
	\leq C_1 \hat \tau^{(1,2)}  (x) ^{d}.
\]
Using the equality $ \hat \tau^{(1,2)}  (x) = \rho\bigl(   x ,  \mathcal N ^ {(1,2)} ( x ) \bigr)  $, we obtain
\[
	\Exp \sbr{ \abs{\hat L^{(1)}(x)}^2 \mid X_3, \dots, X_n }
	\leq 4 C_1 L^2 \, \rho\bigl(x, \mathcal N^{(1,2)}(x)\bigr)^{2s+d}.  
\]
Applying to the term
\[
	 \rbr{Y_1 - Y_ {n,1} }^2 
	 =  \sbr{ \func(X_1) - \hat \func_n^{(1)} (X_1)  - \cbr{\meas(\func) - \meas ( \hat \func_n^{(1)} )}}^2
\]
the same reasoning as above with $A = C = \func(X_1) - \hat \func_n^{(1)} (X_1)$ and $B = D = \meas(\func) - \meas ( \hat \func_n^{(1)} )$,
we get
\begin{align*}
	\Exp \sbr{ \rbr{Y_1 - Y_ {n,1} }^2 }
	&\le \Exp \sbr{\cbr{\func(X_1) - \hat \func_n^{(1)}(X_1)}^2} \\
	&= \Exp \sbr{\cbr{\func(X_1) - \func\bigl(\mathcal N^{(1)}(X_1)\bigr)}^2} \\
	&\le L^2 \Exp \sbr{ \rho\bigl(X_1, \mathcal{N}^{(1)}(X_1)\bigr)^{2s}}.
\end{align*}
Combining all this with \eqref{eq:munNNlooerrdecomp} gives
\begin{align*}
	\lefteqn{\Exp \sbr{\abs{ \munNNloo   (\func)  - \meas   (\func) }^2 }} \\
	& \leq L^2 n^{-1} \Exp \sbr{ \rho\bigl( X_1 ,  \mathcal N ^ {(1)} ( X_1 ) \bigr)^{2s} }
	+ 4  C_1 L^2 \Exp \sbr{ \rho\bigl(   X_1 ,  \mathcal N ^ {(1,2)} ( X_1 ) \bigr)^{2s+d} } \\
	& = L^2 n^{-1} \Exp  [ \hat \tau_{n-1}(X_1)^{2s} ] + 4  C_1 L^2 \Exp [\hat \tau_{n-2}(X_1)^{2s+d}],
\end{align*}
where $X_1$ is understood to be independent of the nearest neighbor distance functions $\hat{\tau}_{n-1}$ (based on $X_2,\ldots,X_n$) and $\hat{\tau}_{n-2}$ (based on $X_3,\ldots,X_n$).
Applying Lemma \ref{lem:1NN-distance} to $\Exp  [\hat \tau_{n-1}(X_1)^{2s} \mid X_1]$ and to $\Exp [ \hat \tau_{n-2}(X_1)^{2s+d} \mid X_1]$, we get
\begin{align*}
	\Exp \sbr{ \hat \tau_{n-1}(X_1)^{2s} }
	&\leq \cbr{ (n-1)  C_0 }^{-2s/d} \Gamma ( 2s/d +1), \\
	\Exp \sbr{ \hat  \tau_{n-2}(X_1)^{2s+d} }
	&\leq \cbr{ (n-2)  C_0 }^{-2s/d-1} \Gamma ( 2s/d +2).
\end{align*}
Therefore,
\begin{align*}
	\Exp \sbr{\abs{\munNNloo(\func)  - \meas(\func)}^2 } 
	&\le L^2 n^{-1} \cbr{(n-1) C_0 }^{-2s/d} \Gamma ( 2s/d +1) \\
	&\quad \null + 4 C_1 L^2 \cbr{(n-2) C_0}^{-2s/d-1} \Gamma ( 2s/d +2)
\end{align*}
and thus
\begin{multline*}
	n^{1+2s/d} \Exp \sbr{\abs{\munNNloo(\func)  - \meas(\func)}^2 } 
	\le L^2 C_0^{-2s/d} \Gamma(2s/d+1) \cbr{n/(n-1)}^{2s/d} \\
	+ 4 C_1 L^2 C_0^{-2s/d-1} \Gamma(2s/d+2) \cbr{n/(n-2)}^{2s/d+1}. 
\end{multline*}
Since $n/(n-2) \le 2$ (from $n \ge 4$), the upper bound can be simplified to
\[
	5 L^2 C_1 C_0^{-2s/d-1} \Gamma(2s/d+2) 2^{2s/d+1}.
\]

\subsection{Proof of Theorem \ref{th:cv_rate} for $\munNN$}

The proof follows from combining the previous part of  Theorem \ref{th:cv_rate} and the following inequality from Proposition~\ref{prop:difference_loonn_nn},
\[
\Exp \sbr{ \abs{ \frac{1}{n}  \sum_{i=1}^n \meas ( \hgni )   - \meas(\hat \func_n ) }^2 } 
\le 4 \rbr{s/d +1}  L^2  C_0^{-2s/d} \, \Gamma ( 2s/d +1) \, n^{-1-2s/d}.
\]
By Minkowski's inequality, we have
\begin{multline*}
	\left(\Exp \Big[  \big| \munNN (\func)  -  \meas(\func) \big|^2\Big] \right)^{1/2} \\
	\le  \left(\Exp \Big[  \big| \munNN (\func) - \munNNloo (\func) \big|^2\Big] \right)^{1/2} + \left(\Exp \Big[  \big| \munNNloo (\func) -  \meas(\func) \big|^2\Big] \right)^{1/2}.
\end{multline*}
In view of \eqref{eq:diffmunloomunn}, $\Exp \Big[  \big| \munNN (\func) - \munNNloo (\func) \big|^2\Big] = \Exp \sbr{ \abs{ \frac{1}{n}  \sum_{i=1}^n \meas ( \hgni )   - \meas(\hat \func_n ) }^2 }$.
Therefore,
\begin{align*}
		\left(\Exp \Big[  \big| \munNN (\func)  -  \meas(\func) \big|^2\Big] \right)^{1/2} \le \sqrt{L^2 C_0^{-2s/d} \Gamma(2s/d+2) n^{-1-2s/d}}  \rbr{2 \sqrt{s/d + 1} + 2^{s/d} (C_1/C_0)^{1/2} \sqrt{10}}.
\end{align*}

\subsection{Proof of Theorem \ref{th:NNlooconcentration} for $\munNNloo$}
	We will apply Theorem~\ref{th:mcdiarmid_extension}, showing the bounded difference property in two parts (Step~2). First, in Step~1, we construct a large-probability event $A$ on which the bounded difference property will hold. In order to bound the gap between $\E[\munNNloo(\func) \mid A]$ and $\meas(\func)$, we rely in Step~3 on the identity $\E[\munNNloo(\func)] = \meas(\func)$ and on Lemma~\ref{rk:expect_diff}. 
	
	\mypar{Step 1}. Let $\lceil x \rceil$ denote the smallest integer upper bound to $x \in \reals$. By Lemma~\ref{lm:nonasymptoticbound}, there exists an event $A$ with probability $\prob(A)\ge 1-\delta_n$ such that on $A$, we have, for $k = \lceil 8 d \log(12n/\delta_n) \rceil $, 
	\begin{align*}
		\sup_{x \in \mathcal{X}} \hat{\tau}_{n,k}(x) 
		\le \left(\frac{2 k }{n b c V_d}\right)^{1/d} 
		&= \left(\frac{2 \lceil 8 d \log(12n/\delta_n) \rceil }{n b c V_d}\right)^{1/d} \\
		&\le \left(\frac{17 d \log(12n/\delta_n)}{n b c V_d}\right)^{1/d} 
		= \bar{\tau}_n.
	\end{align*}
	
	\mypar{Step 2}. On the event $A$, we consider separately two terms of $\munNNloo = n^{-1} \sum_{i=1}^n \{\func(X_i) - \hgni (X_i) + \meas(\hgni)\}$, namely $n^{-1} \sum_{i=1}^n \{\func(X_i) - \hgni (X_i)\}$ and $n^{-1} \sum_{i=1}^n \meas(\hgni)$, in order to make the bounded differences property \eqref{eq:bounded_diff_propeprty} for $\munNNloo$ satisfied. Further, we apply Theorem \ref{th:mcdiarmid_extension}.    
	
	\mypar{Step 2a}. Fix $\ell \in \{1,\ldots,n\}$. In the original set of points $X_1,\ldots,X_n$, replace $X_\ell$ by $\tX_\ell$. Let $\hat{\tilde N}_n^{(i)}(x)$ be the nearest neighbor of $x$ among $X_1,\dots,\tilde{X}_\ell,\dots,X_n$ without $X_i$ and define $\htgni (x) = \func(\hat{\tilde N}_n^{(i)}(x))$. Note that $\hat{\tilde N}_n^{(\ell)}(x)=\hat{N}_n^{(\ell)}(x)$ and hence $\hat{\tilde \func}_n^{(\ell)} (x) = \hat{\func}_n^{(\ell)} (x)$. Let 
	\begin{equation}
		\label{eq:phi1}
		\phi_1 \equiv \phi_1(X_1,\dots,X_n) = \frac{1}{n} \sum_{i=1}^n \{\func(X_i) - \hgni (X_i) \}.
	\end{equation}
	Since $\hat{\tilde \func}_n^{(\ell)} (\tX_\ell) = \hat{\func}_n^{(\ell)} (\tX_\ell)$, we have
	\begin{align*}
		D_{\ell,1} 
		&= \phi_1(X_1,\dots,X_\ell,\dots,X_n) - \phi_1(X_1,\dots,\tX_\ell,\dots,X_n) \\
		&= \frac{1}{n} \sum_{\substack{i=1 \\ i\neq \ell}}^{n} \{\func(X_i) - \hgni (X_i) \} + \frac{1}{n} \{\func(X_\ell) - \hat \func_n^{(\ell)} (X_\ell) \} \\
		&\quad \mbox{} - \frac{1}{n} \sum_{\substack{i=1 \\ i\neq \ell}}^{n} \{\func(X_i) - \hat{\tilde \func}_n^{(i)} (X_i) \} - \frac{1}{n} \{\func(\tX_\ell) - \hat{\tilde \func}_n^{(\ell)} (\tX_\ell) \} \\
		&= \frac{1}{n} \sum_{\substack{i=1 \\ i\neq \ell}}^{n} \{\hat{\tilde \func}_n^{(i)} (X_i) - \hgni (X_i) \} +  \frac{1}{n} \left[\{ \func(X_\ell) - \func(\tX_\ell)\} - \{\hat \func_n^{(\ell)} (X_\ell) - \hat{\func}_n^{(\ell)} (\tX_\ell) \}\right].
	\end{align*}
	
	Considering the first term of $D_{\ell,1}$, we have, in view of \ref{hyp:lip},
	\begin{align*}
		\frac{1}{n} \sum_{\substack{i=1 \\ i\neq \ell}}^{n} 
		\abs{\hat{\tilde \func}_n^{(i)} (X_i) - \hgni (X_i)}
		&=  \frac{1}{n} \sum_{\substack{i=1 \\ i\neq \ell}}^{n} 
		\abs{\func(\hat{\tilde N}^{(i)}_n(X_i)) - \func(\hat N^{(i)}_n(X_i))} \\
		&\le \frac{L}{n} \sum_{\substack{i=1 \\ i\neq \ell}}^{n} 
		\norm{ \hat{\tilde N}^{(i)}_n(X_i) - \hat N^{(i)}_n(X_i) }^s \\
		&= \frac{L}{n} \sum_{\substack{i=1 \\ i\neq \ell}}^{n} 
		\norm{\hat{\tilde N}^{(i)}_n(X_i) - \hat N^{(i)}_n(X_i)}^s
		\ind_{\{\hat{\tilde N}^{(i)}_n(X_i) = \tX_\ell \text{ or } \hat N^{(i)}_n(X_i)=X_\ell\}}.
	\end{align*}
	By the triangle inequality and the sub-additivity of $x \mapsto x^s$, we get, in hopefully obvious notation,
	\begin{align}
		\norm{ \hat{\tilde N}^{(i)}_n(X_i) - \hat N^{(i)}_n(X_i) }^s
		&\le \norm{ X_i - \hat{\tilde N}^{(i)}_n(X_i) }^s
		+ \norm{ X_i - \hat N^{(i)}_n(X_i) }^s \nonumber \\
		&\le \norm{ X_i - \hat{\tilde N}^{(i,\ell)}_n(X_i) }^s 
		+ \norm{ X_i - \hat N^{(i,\ell)}_n(X_i) }^s \nonumber \\
		&\le 2 \norm{ X_i - \hat N^{(i,\ell)}_n(X_i) }^s  
		\le 2 \sup_{x \in \cX} \hat \tau_{n,3}(x)^s, 
		\label{eq:diff_Ntildehat_Nhat}
	\end{align}
	where $\hat \tau_{n,3}(x)$ is the distance to the third nearest neighbor of $x$ among $X_1,\dots,X_n$. By Lemma~\ref{lem:degree_nn},
	\begin{equation}
		\sum_{\substack{i=1 \\ i\neq \ell}}^{n} \ind_{\{\hat{\tilde N}^{(i)}_n(X_i) = \tX_\ell \text{ or } \hat N^{(i)}_n(X_i)=X_\ell\}} \le 2 \psi_d.
		\label{eq:kissing_number_ineq}
	\end{equation}
	Therefore, by \eqref{eq:diff_Ntildehat_Nhat} and \eqref{eq:kissing_number_ineq},
	\begin{equation}
		\label{eq:first_term_Dl1}
		\frac{1}{n} \sum_{\substack{i=1 \\ i\neq \ell}}^{n} 
		\abs{\hat{\tilde \func}_n^{(i)} (X_i) - \hgni (X_i)}
		\le \frac{4 L \psi_d}{n} \sup_{x \in \cX} \hat \tau_{n,3}(x)^s.
	\end{equation}
	
	Considering the second term of $D_{\ell,1}$, we have
	\begin{align}
		\lefteqn{
			\frac{1}{n} \abs{ \{\func(X_\ell) - \func(\tX_\ell)\} - \{\hat \func_n^{(\ell)} (X_\ell) - \hat{\func}_n^{(\ell)} (\tX_\ell)\} } 
		}\nonumber \\
		&\le \frac{L}{n} \left[ \norm{X_\ell - \hat{N}^{(\ell)}_n(X_\ell)}^s + \norm{\tX_\ell - \hat{N}^{(\ell)}_n(\tX_\ell)}^s \right] \nonumber\\
		&\le \frac{2 L}{n} \sup_{x \in \cX} \hat \tau_{n,3}(x)^s.
		\label{eq:second_term_Dl1}
	\end{align}
	
	On the event $A$, by \eqref{eq:first_term_Dl1} and \eqref{eq:second_term_Dl1}, we therefore have, since $k \ge 3$ for our choice of $k$ in Step~1, the bound
	\begin{equation*}
		\abs{ D_{\ell,1} } \le \frac{2 L}{n} (2 \psi_d + 1) \bar{\tau}_n^s.
	\end{equation*}
	
	\mypar{Step 2b}. Let $\phi_2 \equiv \phi_2(X_1,\dots,X_n)= n^{-1} \sum_{i=1}^n \meas (\hgni)$. Therefore,
	\begin{align*}
		D_{l,2} 
		&= \phi_2(X_1,\dots,X_\ell,\dots,X_n) - \phi_2(X_1,\dots,\tX_\ell,\dots,X_n) \\
		&=  \frac{1}{n} \sum_{i=1}^n \{ \meas (\hgni) - \meas (\hat{\tilde \func}_n^{(i)}) \} = \frac{1}{n} \sum_{i=1}^n \meas (\hgni - \hat{\tilde \func}_n^{(i)}) .
	\end{align*}
	Let $S_{n+1,j}$, for $j = 1,\ldots,n+1$, be the Voronoi cells induced by  $\cX_{n+1} = \{ X_1,\ldots,X_n,\tX_\ell \}$: for $x \in S_{n+1,j}$ and $j = 1,\ldots,n+1$, the nearest neighbor of $x$ among $\cX_{n+1}$ is $X_j$, where $X_{n+1} := \tX_\ell$. Clearly,
	\[
	\forall j \in \{1,\ldots,n+1\} \setminus \{i,\ell,n+1\}, \;
	\forall x \in S_{n+1,j}, \qquad
	\htgni(x) = \func(X_j) = \hgni(x).
	\]
	It follows that
	\[
	D_{\ell,2} = \Pn \meas \left( (\hgni - \htgni) \ind_{S_{n+1,i} \cup S_{n+1,\ell} \cup S_{n+1,n+1}} \right).
	\]
	Moreover, by the triangle inequality, the sub-additivity of $x \mapsto x^s$ and \ref{hyp:lip}, we have, for $x \in S_{n+1,i} \cup S_{n+1,\ell} \cup S_{n+1,n+1}$,
	\begin{align*}
		\left| \htgni(x) - \hgni(x) \right|
		&\le \left| \htgni(x) - \func(x) \right| + \left| \hgni(x) - g(x) \right| \\
		&\le L \cdot \left\| \hat{\tilde{N}}_n^{(i)}(x) - x \right\|^s
		+ L \cdot \left\| \hat{N}_n^{(i)}(x) - x \right\|^s \\
		&\le 2L \sup_{x \in \cX} \hat{\tau}_{n,3}(x)^s.
	\end{align*}
	Clearly, $\sum_{i=1}^n \meas(S_{n+1,i}) = \meas\left(\bigcup_{i=1}^n S_{n+1,i}\right) \le 1$. On the event $A$, we thus obtain
	\begin{align*}
		|D_{\ell,2}|
		&\le \Pn \meas \left( 
		\left| \hgni - \htgni \right| 
		\ind_{S_{n+1,i} \cup S_{n+1,\ell} \cup S_{n+1,n+1}}
		\right) \\
		&\le 2L \bar{\tau}_n^s  \cdot
		\Pn \meas(S_{n+1,i} \cup S_{n+1,\ell} \cup S_{n+1,n+1}) \\
		&\le 2L \bar{\tau}_n^s  \cdot 
		\left\{ n^{-1} + \meas(S_{n+1,\ell}) + \meas(S_{n+1,n+1}) \right\}.
	\end{align*}
	
	Recall that $\lambda_d$ denotes the $d$-dimensional Lebesgue-measure and that the density, $\dens$, of $\meas$ satisfies $\sup_{x \in \ms} \dens(x) \le U$. For any $j = 1,\ldots,n+1$, the $\meas$-volume of a Voronoi cell satisfies
	\[
	\meas(S_{n+1,j})
	= \int_{S_{n+1,j}} \dens \, \diff \lambda_d
	\le U \cdot \lambda_d(S_{n+1,j} \cap \ms).
	\]
	Letting $\hat{\tau}_{n+1,1}(x)$ denote the nearest neighbor of $x \in \ms$ in $\cX_{n+1}$, we have
	\[
	\forall x \in S_{n+1,j} \cap \ms, \qquad 
	\norm{x - X_j} = \hat{\tau}_{n+1,1}(x) \le \sup_{x' \in \ms} \hat{\tau}_{n+1,1}(x').
	\]
	On the event $A$, the latter is bounded by $\bar{\tau}_n$, which implies that $S_{n+1,j} \cap \ms$ is contained in a ball of radius $\bar{\tau}_n$ centered at $X_{j}$. We find
	\begin{equation}
		\label{eq:maxmuSn1j}
		\max_{j=1,\ldots,n+1} \meas(S_{n+1,j}) 
		\le U \cdot \bar{\tau}_n^d V_d,
	\end{equation}
	with $V_d$ the volume of the unit ball in $\Rd$. We conclude that, on the event $A$,
	\[
	|D_{\ell,2}| \le 2 L \bar{\tau}_n^s \cdot \left( n^{-1} + 2 UV_d \bar{\tau}_n^d \right).
	\]

	\mypar{Step 2c}. Let $\phi_3=\phi_1+\phi_2 = n^{-1} \sum_{i=1}^n \{\func(X_i) - \hgni (X_i) + \meas(\hgni)\} = \munNNloo$. Then we get, with probability at least $1-\delta_n$, 
	\begin{align*}
		D_{\ell,3} &= \left| \phi_3(X_1,\dots,X_\ell,\dots,X_n) - \phi_3(X_1,\dots,\tX_\ell,\dots,X_n) \right|
		\le |D_{\ell,1}| + |D_{\ell,2}| \\
		&\le 2 L n^{-1} (2 \psi_d + 1) \bar{\tau}_n^s + 2 L \bar{\tau}_n^s \cdot \left( n^{-1} + 2 UV_d \bar{\tau}_n^d \right)\\
		&= 2 L n^{-1} (2 \psi_d + 2) \bar{\tau}_n^s + 4 L U V_d \bar{\tau}_n^{s+d} \\
		&= 4 L (\psi_d + 1) n^{-1} \left(\frac{17 d \log(12n/\delta_n) }{n b c V_d}\right)^{s/d} 
		+ 4 L U V_d \left(\frac{17 d \log(12n/\delta_n) }{n b c V_d}\right)^{1+s/d} \\
		&\le K_1 \left(\frac{\log(12 n/ \delta_n)}{n} \right)^{1+s/d} = c_{\ell,3}
		\quad \text{with} \quad
		K_1 = 4 L \left(\frac{17 d }{b c V_d}\right)^{s/d} \left( \psi_d +1 + \frac{17 d U }{b c}\right).
	\end{align*}
	We have 
	\begin{align*}
		\sum_{\ell=1}^n c_{\ell,3} 
		&= n c_{\ell,3} = K_1 n^{-s/d} (\log(12 n/ \delta_n))^{1+s/d}, \\
		\sum_{\ell=1}^n c_{\ell,3}^2 
		&= n c_{\ell,3}^2 = K_1^2 n^{-1-2s/d} (\log(12 n / \delta_n))^{2+2s/d}.
	\end{align*}
	We apply Theorem~\ref{th:mcdiarmid_extension} to $\phi_3=\phi_3(X_1,\dots,X_n)= n^{-1} \sum_{i=1}^n \{\func(X_i) - \hgni (X_i) + \meas(\hgni)\}$. For $t \ge \delta_n \sum_{\ell=1}^n c_{\ell,3} = \delta_n K_1 n^{-s/d} (\log(12 n/\delta_n))^{1+s/d}$, we get
	\begin{equation}
		\label{eq:phi3McD}
		\prob\left(\abs{\phi_3- \E(\phi_3\mid A)} \ge t\right) 
		\le 
		\delta_n + 
		2 \exp\left(
		- \frac{2 \left\{t-\delta_n K_1 n^{-s/d} (\log(12 n/\delta_n))^{1+s/d}\right\}^2}%
		{K_1^2 n^{-1-2s/d} (\log(12 n/\delta_n))^{2+2s/d}} 
		\right). 
	\end{equation}
	Let $0 < \delta < 1- \delta_n$. Then
	\begin{equation*}
		\prob(\abs{\phi_3- \E(\phi_3\mid A)} \ge t) \le \delta + \delta_n
	\end{equation*}
	provided $t$ is such that the exponential function in \eqref{eq:phi3McD} is bounded by $\delta/2$, which happens if
	\begin{align*}
		t 
		&\ge \sqrt{\frac{ K_1^2 \log(2/\delta) (\log(12 n/ \delta_n))^{2+2s/d}}{2 n^{1+2s/d}}} 
		+ \frac{K_1 \delta_n (\log(12 n/ \delta_n))^{1+s/d}}{n^{s/d}} \\
		&= K_1 \left(n^{-1/2} \sqrt{\log(2/\delta)/2}  + \delta_n \right) \frac{(\log(12 n/ \delta_n))^{1+s/d}}{n^{s/d}}.
	\end{align*}
	
	
	Hence, for $0 < \delta \le 1-\delta_n$, with probability at least $1-(\delta+\delta_n)$, we have
	\begin{equation*}
		\abs{\phi_3 - \E (\phi_3 \mid A)} 
		< K_1 \left(n^{-1/2} \sqrt{\log(2/\delta)/2}  + \delta_n \right) \frac{(\log(12 n/ \delta_n))^{1+s/d}}{n^{s/d}}
	\end{equation*}
	
	\mypar{Step 3}. We have $\ninf{\phi_3} \le 3 \ninf{\func}$. By Lemma~\ref{rk:expect_diff},
	\begin{equation*}
		\abs{\E(\phi_3) - \E(\phi_3 \mid A)}
		\le 2 \delta_n \ninf{\phi_3}
		\le 6 \delta_n \ninf{\func}.
	\end{equation*}
	
	Thus, we obtain, with probability at least $1-(\delta+ \delta_n)$, for $0 < \delta \le 1- \delta_n$,
	\begin{align*}
		\abs{\munNNloo(\func) - \meas(\func)} 
		&= \abs{\phi_3 - \E(\phi_3)} \\
		&\le \abs{\phi_3 - \E(\phi_3 \mid A)} + \abs{\E(\phi_3 \mid A) - \E(\phi_3)} \\
		&\le 
		K_1 \left(n^{-1/2} \sqrt{\log(2/\delta)/2}  + \delta_n \right) \frac{(\log(12 n/ \delta_n))^{1+s/d}}{n^{s/d}}
		+ 6 \delta_n \ninf{\func}. 
	\end{align*}
	
	\mypar{Step 4}. Fix $\varepsilon>0$ and choose $\delta_n = \min(\varepsilon/2, n^{-1/2-s/d})$, $\delta = \varepsilon - \delta_n$. Then $\delta = \varepsilon - \min(\varepsilon/2, n^{-1/2-s/d}) \\= \max(\varepsilon/2, \varepsilon-n^{-1/2-s/d})$. Considering the case when $\varepsilon \le \frac{2}{n^{1/2+s/d}}$, we get
	\begin{align*}
		&\abs{\munNNloo(\func) - \meas(\func)} \\
		&\le K_1 \left(\frac{\sqrt{\log(4/\varepsilon)/2}}{n^{1/2} }  + \frac{\varepsilon}{2} \right) \frac{(\log(24 n/ \varepsilon))^{1+s/d}}{n^{s/d}}	+ 3 \varepsilon \ninf{\func} \\
		&\le K_1 \left(\frac{\sqrt{\log(4/\varepsilon)/2}}{n^{1/2} }  + \frac{1}{n^{1/2+s/d}} \right) \frac{(\log(24 n/ \varepsilon))^{1+s/d}}{n^{s/d}}	+ \frac{6 \ninf{\func}  }{n^{1/2+s/d}} \\
		&\le K_1 \left(\sqrt{\log\left(\frac{4}{\varepsilon}\right) \bigg/2}  + \frac{1}{n^{s/d}} \right) \frac{(\log(24 n/ \varepsilon))^{1+s/d}}{n^{1/2+s/d}}	+ \frac{6 \ninf{\func}  }{n^{1/2+s/d}}.
	\end{align*} 
	Otherwise, when $\varepsilon > \frac{2}{n^{1/2+s/d}}$, using that $\log(12 n^{3/2+s/d}) \le \log((3 n)^3)$, we have
	\begin{align*}
		&\abs{\munNNloo(\func) - \meas(\func)} \\
		&\le K_1 \left( \sqrt{\log\left( \frac{2}{\varepsilon - n^{-1/2-s/d}} \right)\bigg/2}  + \frac{1}{n^{s/d}} \right) \frac{(\log(12 n^{3/2+1/d}))^{1+s/d}}{n^{1/2+s/d}} 	+ \frac{6  \ninf{\func}}{n^{1/2+s/d}} \\
		&\le K_1 \left( \sqrt{\log\left( \frac{4}{\varepsilon} \right)\bigg/2}  + \frac{1}{n^{s/d}} \right) \frac{(3 \log(3 n))^{1+s/d}}{n^{1/2+s/d}} + \frac{6  \ninf{\func}}{n^{1/2+s/d}}.
	\end{align*}
	Therefore, for any $\varepsilon \in (0,1)$, we have, with probability at least $1-\varepsilon$, 
	\begin{align*}
		&\abs{\munNNloo(\func) - \meas(\func)} 
		\le \frac{6 \ninf{\func}  }{n^{1/2+s/d}} \\
		+ & K_1 \left(\sqrt{\log\left(\frac{4}{\varepsilon}\right)  \bigg/2}  + \frac{1}{n^{s/d}} \right)
		\times
		\begin{cases}
			\begin{aligned}[c]
				\frac{(\log(24 n/ \varepsilon))^{1+s/d}}{n^{1/2+s/d}}, 
			\end{aligned}
			& \text{for $\varepsilon \le \frac{2}{n^{1/2+s/d}}$},	 \\		
			\begin{aligned}[c]
				\frac{(3 \log(3 n))^{1+s/d}}{n^{1/2+s/d}}, 
			\end{aligned}
			&\text{for $\varepsilon > \frac{2}{n^{1/2+s/d}}$}.
		\end{cases}	
	\end{align*}
	
Since $\sum_{i=1}^n w_{n,i} = 1$, we can replace $\func$ by $\func - C$ for any constant $C$ and thus, choosing $C$ optimally, replace $\ninf{\func}$ by half the diameter of the range of $\func$, that is, by $C_{\func} = \frac{1}{2} \left\{ \sup_{x \in \ms} \func(x) - \inf_{x \in \ms} \func(x) \right\}$.

\subsection{Proof of Theorem \ref{th:NNlooconcentration} for $\munNN$}
Consider the same event $A$ and the same upper bound $\bar{\tau}_n$ as in Step~1 in the proof of Theorem~\ref{th:NNlooconcentration} for $\munNNloo$. We write
\[
\munNN(\func) 
= \phi_5 \equiv \phi_5(X_1,\ldots,X_n) 
= \phi_1(X_1,\ldots,X_n) + \phi_4(X_1,\ldots,X_n)
\]
with $\phi_1$ as in Eq.~\eqref{eq:phi1} and with 
\[
\phi_4 \equiv \phi_4(X_1, \ldots, X_n) = \meas(\hat{\func}_n).
\]

\mypar{Step 1}. Let $\hat{\tilde N}_n(x)$ be the nearest neigbor of $x$ among $X_1,\dots,\tilde{X}_\ell,\dots,X_n$  and let $\hat{\tilde \func}_n (x) = g (\hat{\tilde N}_n(x))$. Then
\begin{align*}
	D_{\ell,4} 
	&= \phi_4(X_1,\dots,X_\ell,\dots,X_n) - \phi_4(X_1,\dots,\tX_\ell,\dots,X_n) \\
	&=  \meas (\hat \func_n) - \meas (\hat{\tilde \func}_n) = \meas (\hat \func_n - \hat{\tilde \func}_n).
\end{align*}

For $i = 1,\ldots,n+1$, define $S_{n+1,i}$ as in Step~2b in the proof of Theorem~\ref{th:NNlooconcentration} for $\munNNloo$. It follows that
\[
D_{\ell,4} = \meas \left( (\hat{\func}_n - \hat{\tilde \func}_n) \ind_{S_{n+1,\ell} \cup S_{n+1,n+1}} \right).
\]
Further, by the triangle inequality, the sub-additivity of the function $x \mapsto x^s$ and the Hölder-$s$ property, we have, for $x \in \ms \cap \left(S_{n+1,\ell} \cup S_{n+1,n+1}\right)$,
\begin{align*}
	\abs{ \hat{\func}_n(x) - \hat{\tilde \func}_n(x) }
	&\le \abs{ \hat{\func}_n(x) - g(x) } + \abs{ \hat{\tilde \func}_n - g(x) } \\
	&\le L \cdot \norm{ \hat{\tilde{N}}_n(x) - x }^s 
	+ L \cdot \norm{ \hat{N}_n(x) - x }^s \\
	&\le 2L \sup_{x' \in \ms} \hat{\tau}_{n,2}(x')^s.
\end{align*}

On the event $A$, we thus obtain
\begin{align*}
	\abs{D_{\ell,4}}
	\le 2L \bar{\tau}_n^s \cdot \meas (S_{n+1,\ell} \cup S_{n+1,n+1}) 
	\le 2L \bar{\tau}_n^s  \cdot 
	\left\{ \meas(S_{n+1,\ell}) + \meas(S_{n+1,n+1}) \right\}.
\end{align*}

By Eq.~\eqref{eq:maxmuSn1j} in Step~2b in the proof of Theorem~\ref{th:NNlooconcentration} for $\munNNloo$, we have, on the event $A$,
\[
\max_{j=1,\ldots,n+1} \meas(S_{n+1,j}) 
\le U \cdot \bar{\tau}_n^d V_d,
\]
with $V_d$ the volume of the unit ball in $\Rd$.
Therefore, we conclude that, on the event $A$, we have
\begin{equation}
	\label{eq:c4}
	\abs{D_{\ell,4}} \le 4 L UV_d \bar{\tau}_n^{s+d}.
\end{equation}

\mypar{Step 2}. By \eqref{eq:c4} and by the result of Step~2a in the proof of Theorem~\ref{th:NNlooconcentration} for $\munNNloo$, we have
\begin{align*}
	D_{\ell,5}
	&= \abs{ \phi_5(X_1,\ldots,X_\ell,\ldots,X_n) - \phi_5(X_1,\ldots,\tX_\ell,\ldots,X_5) } \\
	&\le \abs{D_{\ell,1}} + \abs{D_{\ell,4}} \\
	&\le 2L n^{-1} (2\psi_d + 1) \bar{\tau}_n^s + 4 L UV_d \bar{\tau}_n^{s+d} \\
	&=  2L n^{-1} (2\psi_d + 1) \left(\frac{17 d  \log(12 n/\delta_n) }{n b c V_d}\right)^{s/d} +  4 L UV_d \left(\frac{17 d \log(12 n/\delta_n) }{n b c V_d}\right)^{1+s/d} \\
	&\le K_2 \left(\frac{\log(12n/\delta_n)}{n} \right)^{1+s/d} = c_{\ell,5}
	\quad \text{with} \quad
	K_2 = 2 L \left( \frac{17 d }{b c V_d}\right)^{s/d} \left( 2 \psi_d +1 + \frac{34 d U}{b c}\right).
\end{align*}
We have
\begin{align*}
	\sum_{\ell=1}^n c_{\ell,5} 
	&= n c_{\ell,5} = K_2 n^{-s/d} (\log(12 n/\delta_n))^{1+s/d}, \\
	\sum_{\ell=1}^n c_{\ell,5}^2 
	&= n c_{\ell,5}^2 =  K_2^2 n^{-1-2s/d} (\log(12 n/ \delta_n))^{2+2s/d}.
\end{align*}
Thanks to Theorem~\ref{th:mcdiarmid_extension}, it follows that, for $t \ge \delta_n \sum_{\ell=1}^n c_{\ell,5} = \delta_n K_2 n^{-s/d} (\log(12 n/\delta_n))^{1+s/d}$,
\[
\prob \left(\abs{\phi_5 - \E(\phi_5 \mid A)} \ge t \right)
\le  
\delta_n + 
2 \exp \left( 
- \frac{2 \left\{t- \delta_n K_2 n^{-s/d} (\log(12 n/\delta_n))^{1+s/d}\right\}^2}%
{K_2^2 n^{-1-2s/d} (\log(12 n/\delta_n))^{2+2s/d}} 
\right).
\]
Let $0 < \delta \le 1 - \delta_n$. Then
\[
\prob \left(\abs{\phi_5 - \E(\phi_5 \mid A)} \ge t \right)
\le \delta + \delta_n
\]
provided
\[
t \ge  K_2 \left( n^{-1/2}\sqrt{\log(2/\delta)/2} + \delta_n \right) \frac{(\log(12 n/\delta_n))^{1+s/d}}{n^{s/d}}.
\]
Hence, for $0 < \delta < 1 - \delta_n$, with probability at least $1 - (\delta + \delta_n)$, we have
\begin{equation}
	\label{eq:phi5Combes}
	\abs{\phi_5 - \E(\phi_5 \mid A)} <  
	K_2 \left( n^{-1/2}\sqrt{\log(2/\delta)/2} + \delta_n \right) \frac{(\log(12 n/\delta_n))^{1+s/d}}{n^{s/d}}.
\end{equation}

\mypar{Step 3.} 
Since $\ninf{\phi_5} \le 3 \ninf{\func}$, we have, by Lemma~\ref{rk:expect_diff},
\[
\abs{\E(\phi_5 \mid A) - \E(\phi_5)}
\le 2 \delta_n \ninf{\phi_5}
\le 6 \delta_n \ninf{\func}.
\]
Recall $\expec[\munNNloo (\func)] = \meas (\func)$. On the event on which \eqref{eq:phi5Combes} holds, we have, in view of Proposition \ref{prop:difference_loonn_nn} with $C_0 = V_d b c$, for $n \ge 4$, 
\begin{align*}
	\lefteqn{
		\abs{\munNN(\func) - \meas(\func)}
	} \\
	\le& \abs{\phi_5 - \E(\phi_5 \mid A)}
	+ \abs{\E(\phi_5 \mid A) - \E(\phi_5)}
	+ \abs{\E\left[\munNN(\func)\right] - \E\left[\munNNloo(\func)\right]} \\
	\le& K_2 \left( n^{-1/2}\sqrt{\log(2/\delta)/2} + \delta_n \right) \frac{(\log(12 n/\delta_n))^{1+s/d}}{n^{s/d}} + 6 \delta_n \ninf{\func} \\
	&+ \frac{(s/d +2)  L  (V_d b c)^{-s/d} \Gamma ( s/d +1)} {n^{1+s/d}}	 \\
	\le&  K_2 \left( n^{-1/2}\sqrt{\log(2/\delta)/2} + \delta_n \right) \frac{(\log(12 n/\delta_n))^{1+s/d}}{n^{s/d}}  +  6 \delta_n \ninf{\func} \\
	&+ \frac{(s/d + 2)  L  (V_d b c)^{-s/d}} {n^{1+s/d}}.
\end{align*}

\mypar{Step 4.} Similarly as in Step 4 of the proof of Theorem \ref{th:NNlooconcentration} for $\munNNloo$, fix $\varepsilon>0$ and choose $\delta_n = \min(\varepsilon/2, n^{-1/2-1/d})$, $\delta=\varepsilon-\delta_n$. Then, for any $\varepsilon \in (0,1)$, we have, with probability at least $1-\varepsilon$, 
\begin{align*}
	\abs{\munNN(\func) - \meas(\func)} 
	\le  
	\frac{6  \ninf{\func}}{n^{1/2+s/d}} + & \frac{(s/d +2)  L  (V_d b c)^{-s/d}} {n^{1+s/d}}  
	+ K_2 \left(\sqrt{\log\left(\frac{4}{\varepsilon}\right) \bigg/2}  + \frac{1}{n^{s/d}} \right)\\ 
	&\times	
	\begin{cases}
		\begin{aligned}[c]
			\frac{(\log(24 n/ \varepsilon))^{1+s/d}}{n^{1/2+s/d}}, 
		\end{aligned}
		& \text{for $\varepsilon \le \frac{2}{n^{1/2+s/d}}$},	 \\		
		\begin{aligned}[c]
			\frac{(3 \log(3 n))^{1+s/d}}{n^{1/2+s/d}}, 
		\end{aligned}
		&\text{for $\varepsilon > \frac{2}{n^{1/2+s/d}}$}.
	\end{cases}	
\end{align*}

\subsection{Proof of Theorem \ref{th:mcdiarmid_extension}}
We have
\begin{align*}
	\prob \left(|\phi(X) - m|\ge t \right)
	&\le \prob \left(|\phi(X) - m|\ge t, X \in A \right) + \prob(X \not\in A) \\
	&\le \prob \left( \phi(X) - m \ge t, X \in A \right) + \prob \left( -\phi(X) + m \ge t, X \in A \right) + p.
\end{align*}
From the proof of Theorem 2.1 in \cite{extensionofMcDiarmid2015}, we get
\begin{equation*}
	\prob \left( \phi(X) - m \ge t, X \in A \right) \le \exp \left(- \frac{2 \max(0,t-p \bar{c})^2}{\sum_{\ell=1}^n c_\ell^2} \right).
\end{equation*}
By symmetry,
\begin{equation*}
	\prob \left( -\phi(X) + m \ge t, X \in A \right) \le \exp \left(- \frac{2 \max(0,t-p \bar{c})^2}{\sum_{\ell=1}^n c_\ell^2} \right).
\end{equation*}
Thus, we obtain the stated result.

\section{Additional Experiments}
\label{app:option_pricing}

\mypar{Finance background.} Options are financial derivatives based on the value of underlying securities. They give the buyer the right to buy (call option) or sell (put option) the underlying asset at a pre-determined price within a specific time frame. The price of an option may be expressed as the expectation, under the so-called risk-neutral measure, of the payoff discounted to the present value. Consider a contract of European type, which specifies a payoff $V(S_T)$, depending on the level of the underlying asset $S_t$ at maturity $t=T$. The value $V$ of the contract at time $t=0$ conditional on an underlying value $S_0$ is 
\begin{equation} \label{eq:price}
V(S_0) = \expec_Q[e^{-rT} V(S_T)],
\end{equation}
where $\expec_Q$ denotes the expectation under the risk-neutral measure and $r$ is the risk-free interest rate. Such a representation suggests a straightforward Monte Carlo based method for its calculation by simulating random paths of the underlying asset, calculating each time the resulting payoff and taking the average of the result. This approach is particularly useful when dealing with exotic options, for which the above expectation often does not permit a closed-form expression.

The payoff of a European call option with strike price $K$ is given by $V(S_T) = (S_T - K)_{+}$ and depends only on the level of the underlying asset $S_t$ at maturity time $t=T$. In contrast, the payoff a of barrier option \citep{merton1973theory} depends on the whole path $(S_t)_{t \in [0,T]}$. The option becomes worthless or may be activated upon the crossing of a price point barrier denoted $H$. More precisely, Knock-Out (KO) options expire worthlessly when the underlying's spot price crosses the pre-specified barrier level whereas Knock-In (KI) options only come into existence if the pre-specified barrier level is crossed by the underlying asset's price. The payoffs of \textit{up-in} (UI) and \textit{up-out} (UO) barrier options with barrier price $K$ are given by
\begin{equation}
	\label{eq:barrier_payoffs}
	V_{\mathrm{(UI)}}(S) = 
	(S_T - K)_{+} 
	\1 \bigl\{\max_{t \in [0,T]} S_t \geq H \bigr\}, \quad
	V_{\mathrm{(UO)}}(S) = 
	(S_T - K)_{+} 
	\1 \bigl\{\max_{t \in [0,T]} S_t < H \bigr\}.
\end{equation}

\mypar{Market Dynamics.} The Black--Scholes model \citep{black2019pricing} is a mathematical model for pricing option contracts. It is based on geometric Brownian motion with constant drift and volatility so that the underlying stock $S_t$ satisfies the following stochastic differential equation:
\begin{align*}
\diff S_{t}=\mu S_{t} \, \diff t+\sigma S_{t} \, \diff W_{t},
\end{align*}
where $\mu$ represents the drift rate of growth of the underlying stock, $\sigma$ is the volatility and $W$ denotes a Wiener process. Although simple and widely used in practice, the Black--Scholes model has some limitations. In particular, it assumes constant values for the risk-free rate of return and volatility over the option duration. Neither of those necessarily remains constant in the real world. The Heston model \citep{heston1993closed} is a type of stochastic volatility model that can be used for pricing options on various securities. For the Heston model, the previous constant volatility $\sigma$ is replaced by a stochastic volatility $v_t$ which follows an Ornstein--Uhlenbeck process. The underlying stock $S_t$ satisfies the following equations
\begin{align*}
\left\{
    \begin{array}{@{\,}l@{\,}l}
        \diff S_t &= \mu S_t \, \diff t + \sqrt{v_t} S_t \, \diff W_t^S, \\[1ex]
        \diff v_t &= \kappa(\theta-v_t) \, \diff t + \xi \sqrt{v_t}  \, \diff W_t^v, \qquad \diff W_t^S \, \diff W_t^v = \rho \, \diff t.
    \end{array}
\right.
\end{align*}
with stochastic volatility $v_t$, drift term $\mu$, long run average variance $\theta$, rate of mean reversion $\kappa$ and volatility of volatility $\xi$. Essentially the Heston model is a geometric Brownian motion with non-constant volatility, where the change in $S$ has relationship $\rho$ with the change in volatility. 

\mypar{Monte Carlo procedures.}
The application of standard Monte Carlo methods to option pricing takes the following form:
\begin{enumerate}[label=(\arabic*)]
\item Simulate a large number $n$ of price paths for the underlying asset: $(S_{(1)}, \ldots, S_{(n)})$. 
\item For each path, compute the associated payoff of the option, e.g., as in Eq.~\eqref{eq:barrier_payoffs}: $(V_1, \ldots, V_n).$ 
\item Average the payoffs and discount them to present value: $\hat V_n = (e^{-rT}/n) \sum_{i=1}^n V_i$.
\end{enumerate}
In practice, the price paths are simulated using an Euler scheme with a discretization of the time period $[0,T]$ comprised of $m$ times $t_1 = 0 < t_2 < \ldots < t_m = T$. Each price path $S_{(i)}$ for $i=1,\ldots,n$ is actually a vector $(S_{(i)}^{(1)}, \ldots, S_{(i)}^{(m)})$, so that the indicator function of the barrier options is computed on the discretized prices. Common values for $m$ are the number of trading days per year which is $m=252$ for $T=1$ year.

\mypar{Parameters.} Several numerical experiments are performed for the pricing of European barrier call options ``up-in'' and ``up-out''. The number of sampled paths is $n \in \{500, 1\,000, 2\,000, 3\,000, 5\,000\}$ and the granularity of the grid is equal to $m=240$. Two different mathematical models are considered when simulating the underlying asset price trajectories:
\begin{enumerate}[label=(\arabic*)]
\item the Black--Scholes model with constant volatility $\sigma = 0.30$;
\item the Heston model with initial volatility $v_0=0.1$, long-run average variance $\theta=0.02$, rate of mean reversion $\kappa=4$, instantaneous correlation $\rho=0.8$ and volatility of volatility $\xi=0.9$.
\end{enumerate}
In both cases the fixed parameters are: spot price $S_0=100$, interest rate $r=0.10$, maturity $T=2$ months,  strike price $K=S_0 = 100$ and barrier price $H=130$.

\begin{figure}[h]
		\centering
		\begin{subfigure}[h]{0.242\textwidth}
			\centering
			\includegraphics[width=\textwidth]{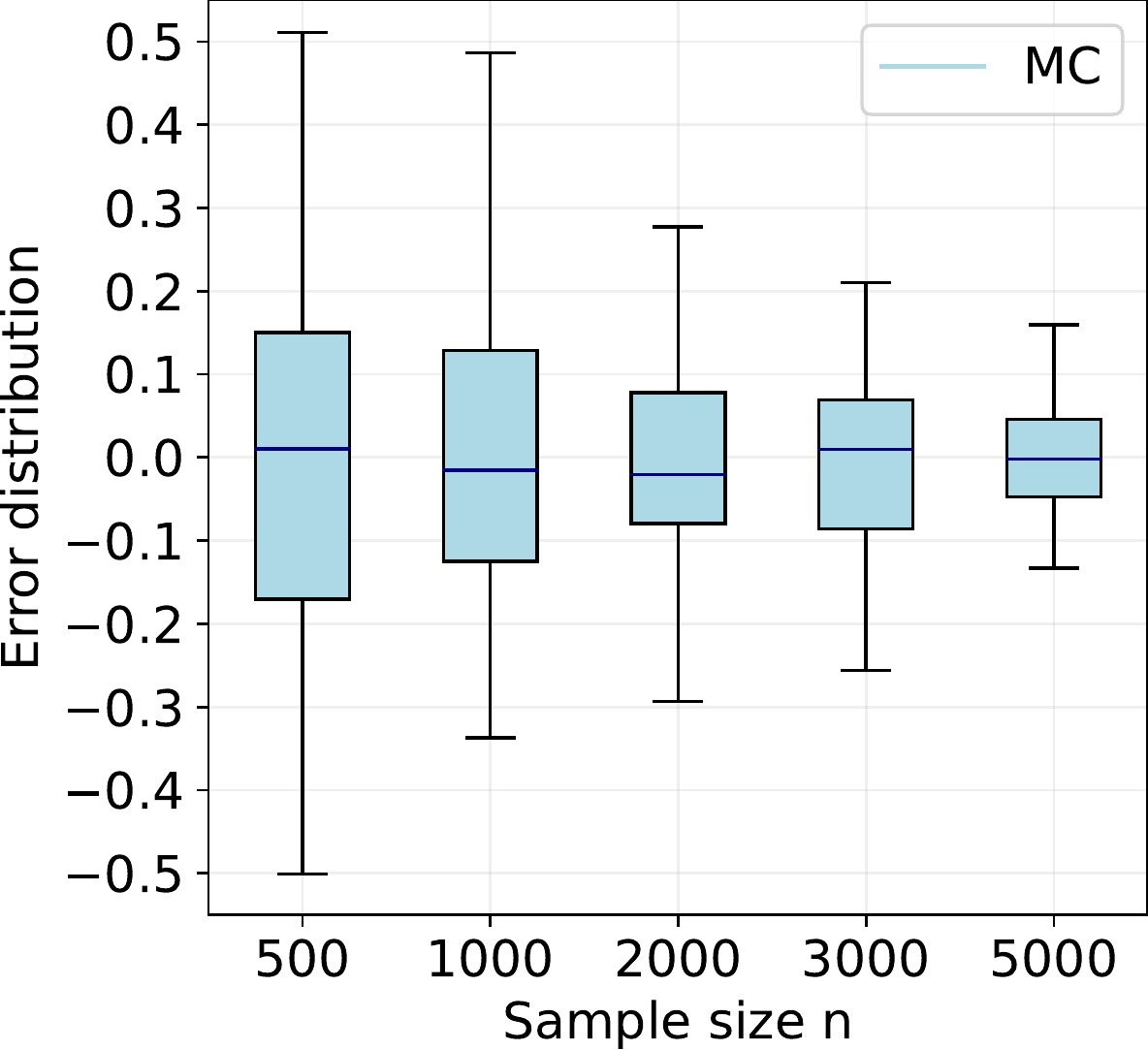}
			\caption{``Up-In'', MC}
			\label{fig:ui_mc}
		\end{subfigure}
		\begin{subfigure}[h]{0.242\textwidth}
			\centering
			\includegraphics[width=\textwidth]{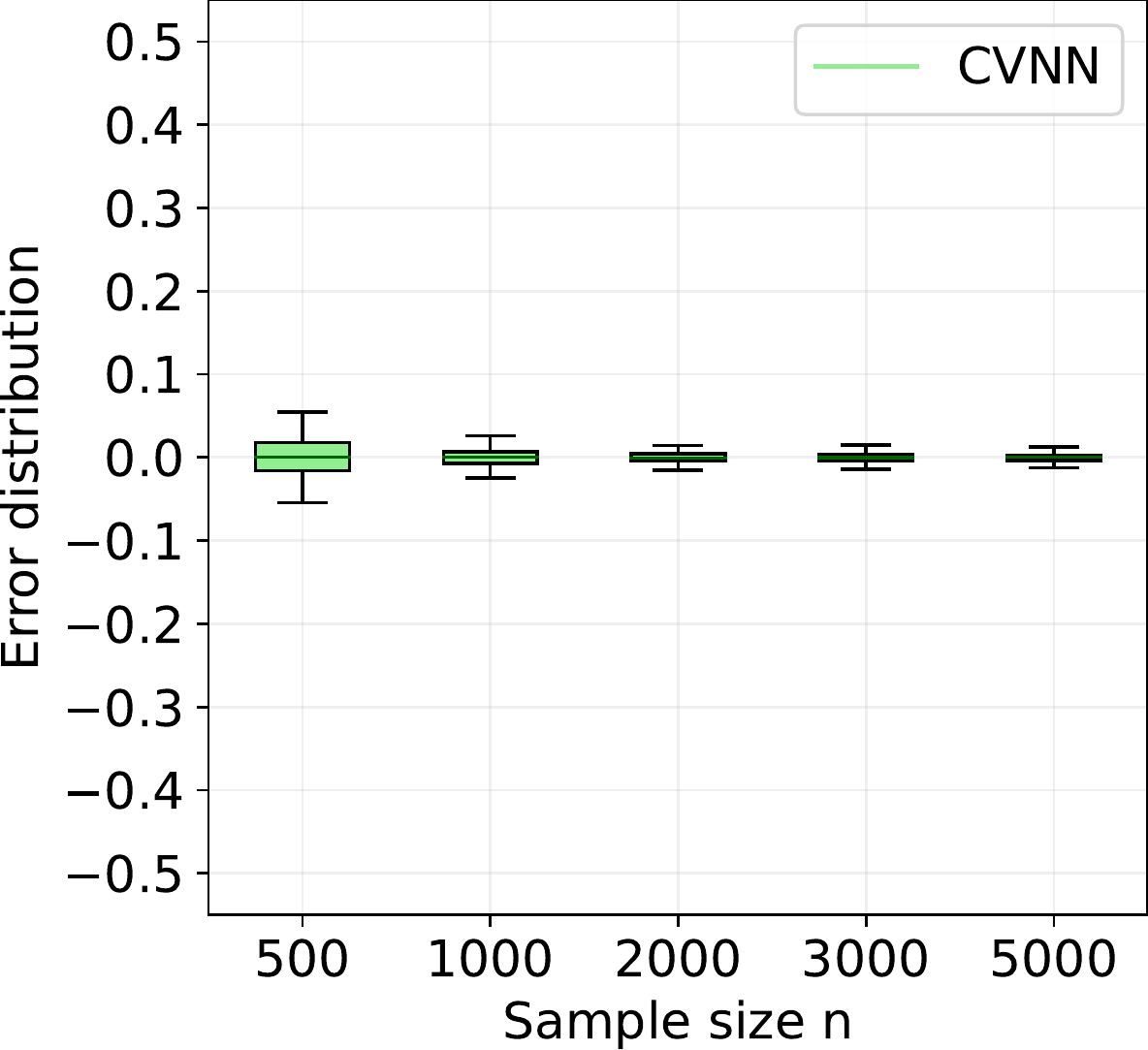}
			\caption{``Up-In'', CVNN}
			\label{fig:ui_cv}
		\end{subfigure}
		\begin{subfigure}[h]{0.242\textwidth}
			\centering
			\includegraphics[width=\textwidth]{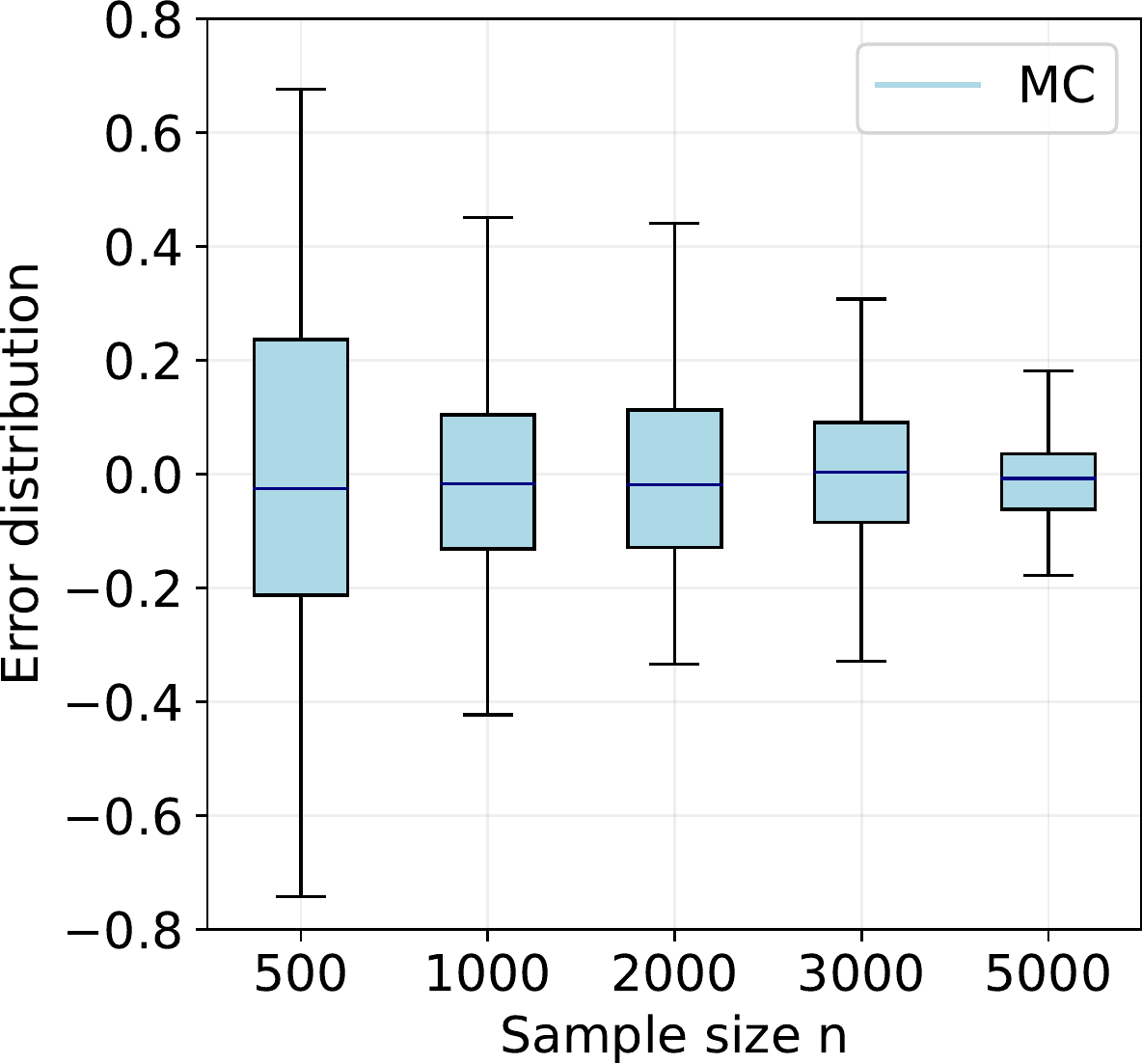}
			\caption{``Up-Out'', MC}
			\label{fig:uo_mc}
		\end{subfigure} 
		\begin{subfigure}[h]{0.242\textwidth}
			\centering
			\includegraphics[width=\textwidth]{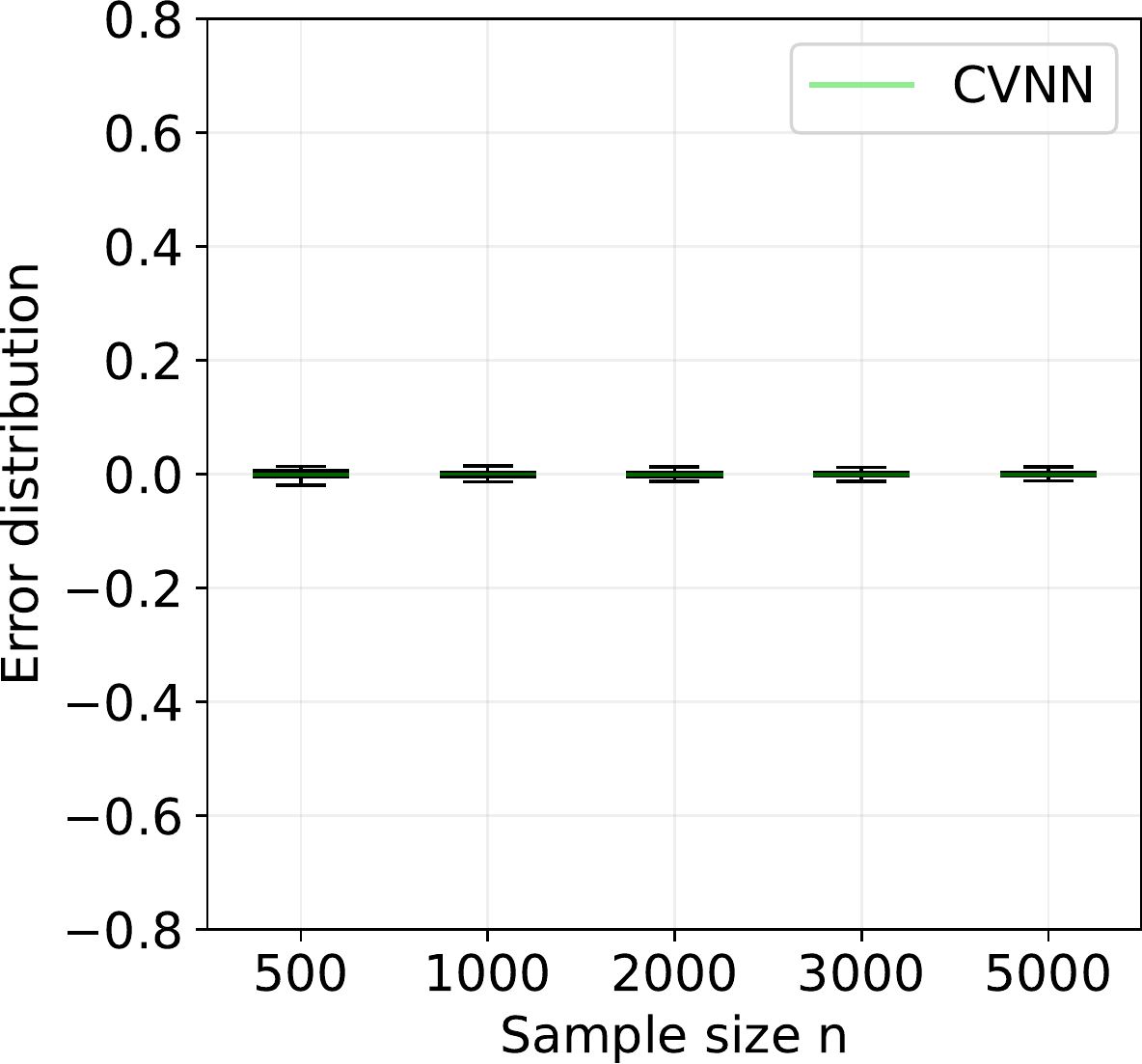}
			\caption{``Up-Out'', CVNN}
			\label{fig:uo_cv}
			\end{subfigure}
		\caption{Barrier option pricing under Black--Scholes model with spot price $S_0=100$, strike $K=S_0$, maturity $T=2$ months, risk-free rate $r=0.1$, constant volatility $\sigma=0.3$, barrier price $H=130$. The boxplots are obtained over $100$ replications.}
		\label{fig:barrier_options} 
	\end{figure}

\begin{figure}[h]
		\centering
		\begin{subfigure}[h]{0.242\textwidth}
			\centering
			\includegraphics[width=\textwidth]{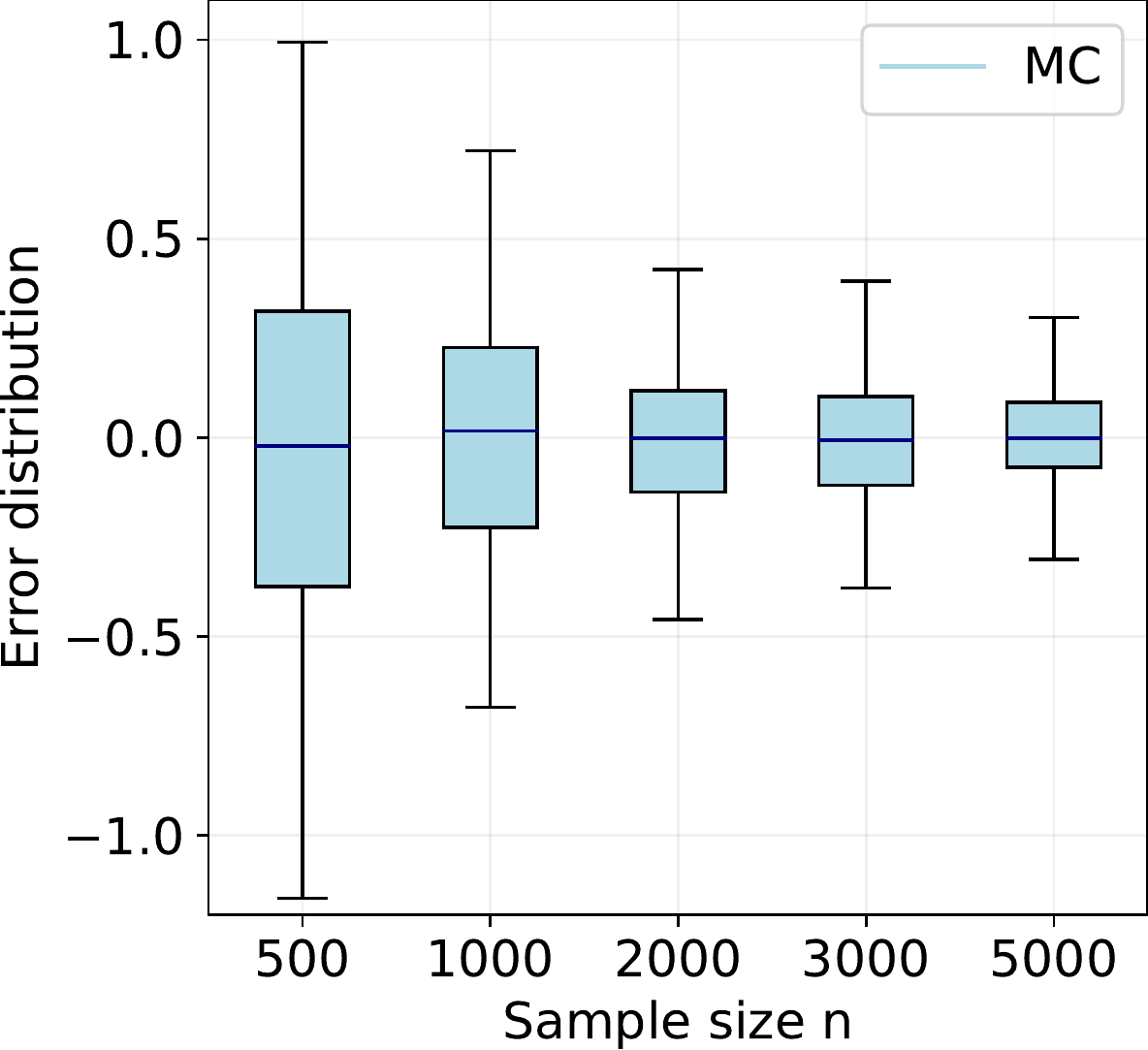}
			\caption{``Up-In'', MC}
			\label{fig:ui_mc_heston}
		\end{subfigure}
		\hfill
		\begin{subfigure}[h]{0.242\textwidth}
			\centering
			\includegraphics[width=\textwidth]{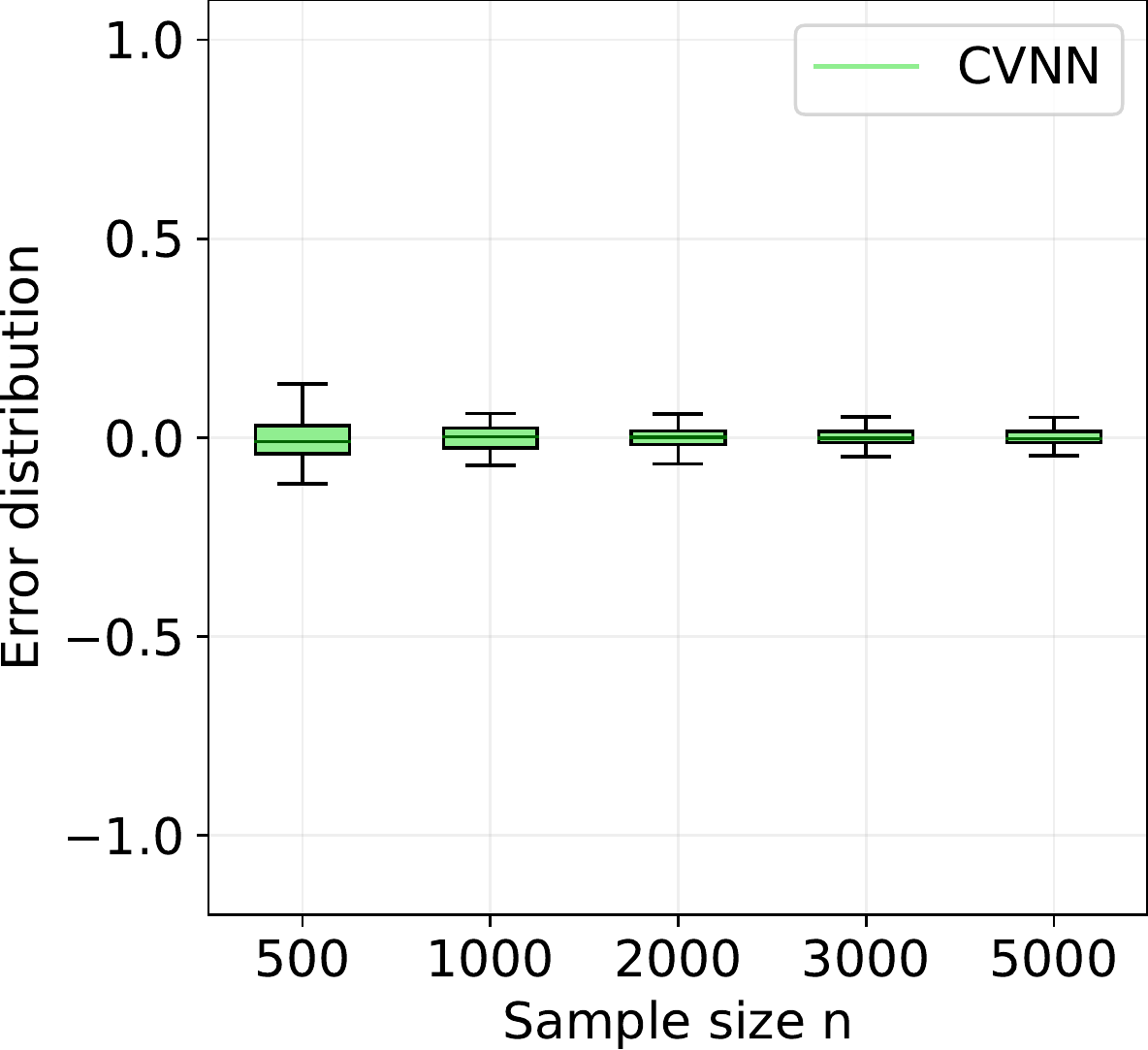}
			\caption{``Up-In'', CVNN}
			\label{fig:ui_cv_heston}
		\end{subfigure}
		\hfill
		\begin{subfigure}[h]{0.242\textwidth}
			\centering
			\includegraphics[width=\textwidth]{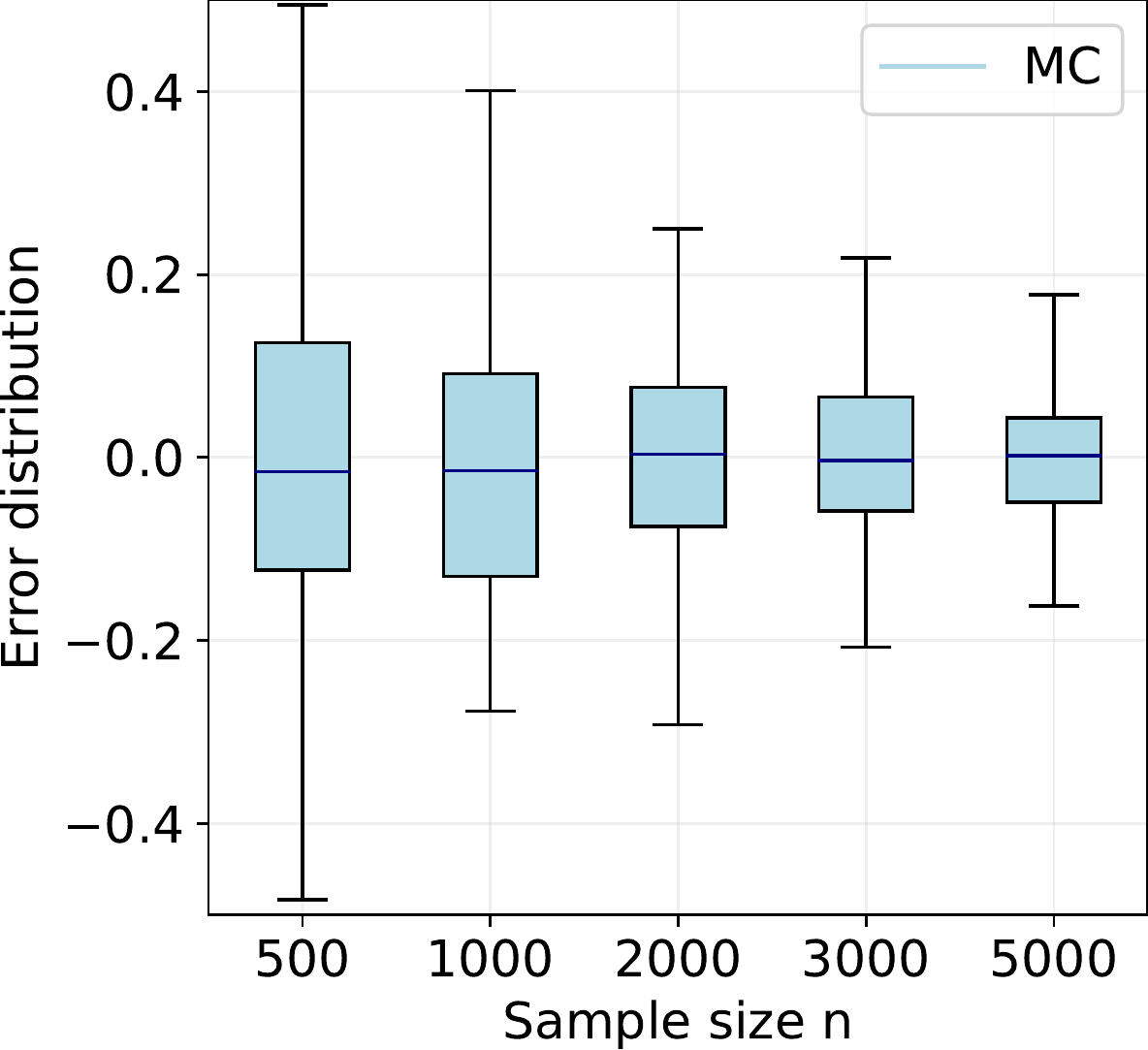}
			\caption{``Up-Out'', MC}
			\label{fig:uo_mc_heston}
		\end{subfigure} 
		\begin{subfigure}[h]{0.242\textwidth}
			\centering
			\includegraphics[width=\textwidth]{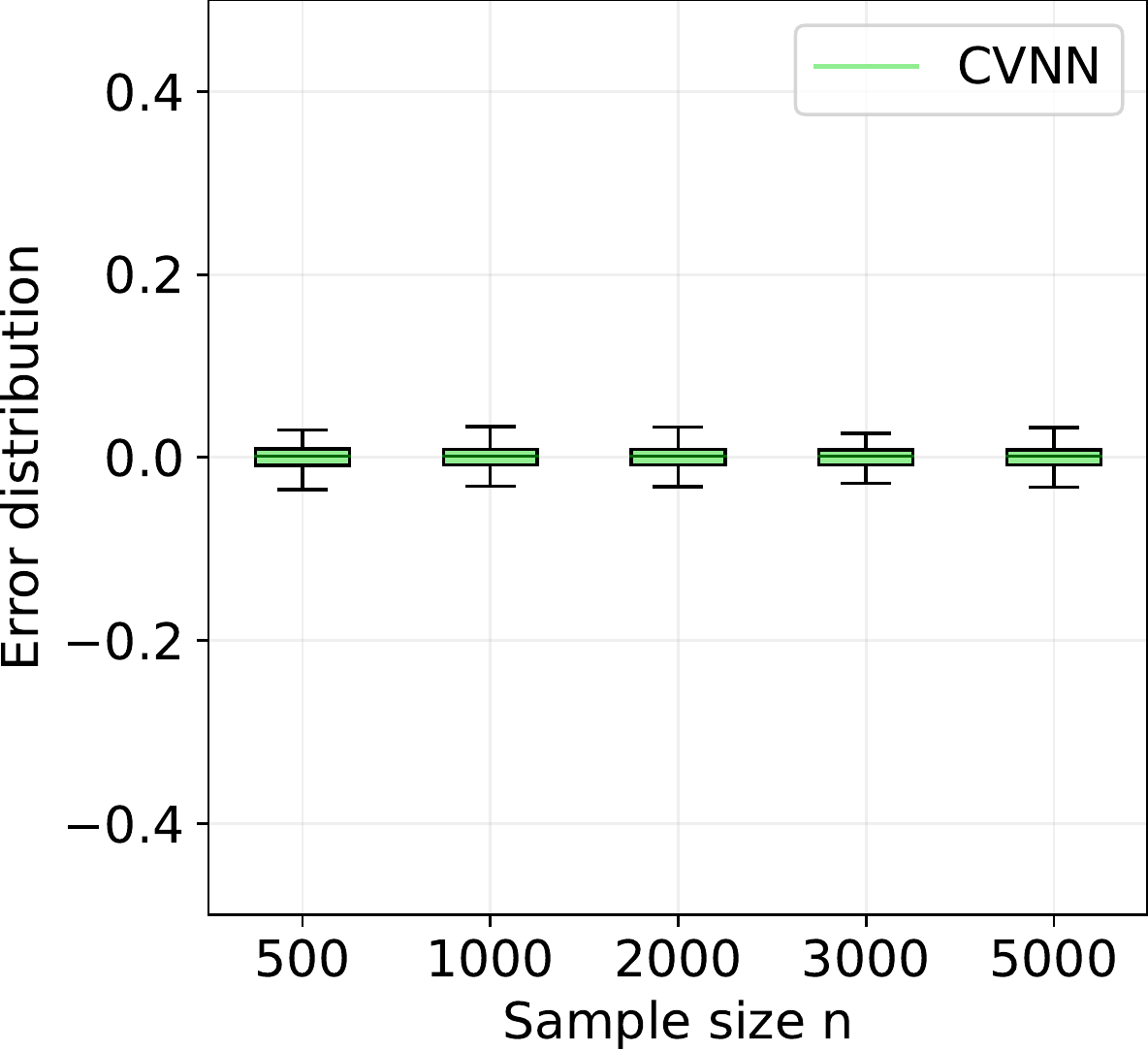}
			\caption{``Up-Out'', CVNN}
			\label{fig:uo_cv_heston}
			\end{subfigure}
		\caption{Boxplots of barrier option pricing with Heston Model with spot price $S_0=100$, strike $K=S_0$, barrier price $H=130$, maturity $T=2$ months, risk-free rate $r=0.1$, initial volatility $v_0=0.1$, long-run average variance $\theta=0.02$, rate of mean reversion $\kappa=4$, instantaneous correlation $\rho=0.8$ and volatility of volatility $\xi=0.9$. The boxplots are obtained over $100$ replications.}
		\label{fig:barrier_options_heston} 
	\end{figure}

\mypar{Results.}
Figure~\ref{fig:barrier_options} shows the error distribution of the different Monte Carlo estimates (naive MC and CVNN) for the pricing of Barrier call options ``up-in'' and ``up-out'' in the Black--Scholes model. The boxplots are computed over $100$ independent replications and the true values of the options are approximated using the Python package \textsf{QuantLib}. Similarly, Figure~\ref{fig:barrier_options_heston} gathers the results for the Heston model. The variance is greatly reduced when using the control neighbors estimate compared to the standard Monte Carlo approach.


\end{document}